\documentclass[11pt,a4paper,openright]{book}

\usepackage{nomencl}
\usepackage{comment}
\usepackage{caption} 
\usepackage{epigraph}
\usepackage{rotating}
\usepackage[section]{placeins}
\usepackage{breakcites}
\usepackage{pst-node}
\usepackage{tikz-cd}

\usepackage{lipsum}  

\usepackage{longtable}

\newcommand{\NewAppendix}[1]{
    \clearpage

    \pagenumbering{arabic}

    \renewcommand{\thepage}{\thesection - \arabic{page}}

}

\RequirePackage{fancyhdr}
\RequirePackage{ccaption}
\RequirePackage{xspace}
\RequirePackage{graphicx} 
\RequirePackage{ifpdf}
\RequirePackage{ifthen}

\usepackage{subcaption}
\usepackage{tabularx,booktabs}  
\usepackage{array} 
\usepackage{multicol} 
\usepackage{multirow}

\usepackage{ccaption}
\captionnamefont{\scshape}
\captionstyle{}
\makeatletter
\renewcommand{\fnum@figure}[1]{\quad\small\textsc{\figurename~\thefigure}:}
\renewcommand{\fnum@table}[1]{\quad\small\textsc{\tablename~\thetable}:}
\renewcommand{\@makecaption}[2]{%
\vskip\abovecaptionskip
\sbox\@tempboxa{#1: #2}%
\ifdim \wd\@tempboxa >\hsize
\def\baselinestretch{1}\@normalsize
#1: #2\par
\def\baselinestretch{1.5}\@normalsize
\else
\global \@minipagefalse
\hb@xt@\hsize{\hfil\box\@tempboxa\hfil}%
\fi
\vskip\belowcaptionskip}
\makeatother

\renewcommand{\baselinestretch}{1}

\newenvironment{itemise}{\begin{itemize}}{\end{itemize}}

\usepackage[grey,times]{quotchap} 

\usepackage[nottoc,notlof,notlot]{tocbibind}  

\newcommand\lrsq[1]{\left[#1\right]}
\newcommand\lr[1]{\left(#1\right)}

\usepackage[pdftitle={PhD Thesis Title}, pdfauthor={Juan Silverio Martínez Baena}, pdfsubject={Thesis},  
 bookmarks,
 bookmarksopen = true,
 bookmarksnumbered = true,
 breaklinks = true,
 linktocpage,
 hyperindex = false,
 hyperfigures,
 pagebackref,
 colorlinks=true,
 linkcolor=blue,
 urlcolor=blue,
 citecolor=red,
 anchorcolor=green
 ]{hyperref}
\usepackage{hyperref}

\usepackage{listings}
\usepackage{color}
\usepackage{tabularx} 
\newcolumntype{A}{>{\arraybackslash}X}
\newcolumntype{L}[1]{>{\raggedright\arraybackslash}p{#1}}
\newcolumntype{C}[1]{>{\centering\arraybackslash}p{#1}}
\newcolumntype{R}[1]{>{\raggedleft\arraybackslash}p{#1}}

\definecolor{dkgreen}{rgb}{0,0.6,0}
\definecolor{gray}{rgb}{0.5,0.5,0.5}
\definecolor{mauve}{rgb}{0.58,0,0.82}

\usepackage[toc,page]{appendix}
\usepackage{graphicx}
\usepackage{amsthm}
\usepackage{amsfonts}
\usepackage{amsmath,amssymb}
\usepackage{multirow}
\usepackage{lscape}
\usepackage{textcomp}
\usepackage[countmax]{subfloat}
\usepackage[linesnumbered,ruled,vlined]{algorithm2e}
\usepackage{mathrsfs}

\theoremstyle{plain}

\theoremstyle{definition}

\usepackage{setspace}

\theoremstyle{plain}
\newtheorem{theorem}{Theorem}[chapter]
\newtheorem{corollary}[theorem]{Corollary}
\newtheorem{lemma}[theorem]{Lemma}
\newtheorem{proposition}[theorem]{Proposition}

\theoremstyle{definition}
\newtheorem{definition}[theorem]{Definition}
\theoremstyle{remark}
\newtheorem{remark}[theorem]{Remark}
\usepackage{mathtools}
\mathtoolsset{showonlyrefs,showmanualtags}
\usepackage{enumitem}
\usepackage{dsfont}
\usepackage{xcolor}
\definecolor{greyforboxes}{RGB}{235,235,235}

\usepackage{tikz}                            
\usepackage[framemethod=tikz]{mdframed}      
\mdfsetup{splittopskip=\topskip, linewidth=0.5}
            
\newenvironment{boxmathdark}[1]
    {\begin{mdframed}[topline=false,rightline=false, bottomline=false,
                      linewidth=2mm, linecolor=black, backgroundcolor=greyforboxes, 
                      innertopmargin = 8pt, innerbottommargin = 8pt,
                      skipabove = 8pt, skipbelow = 8pt]}
    {\end{mdframed}}

\begin{document}

\frontmatter

\pagenumbering{roman}  
\hypersetup{urlcolor=black}{}
\thispagestyle{empty}%
    \null
    \begin{center}
            \huge\sc\expandafter{A priori estimates of stable and finite Morse index solutions to elliptic equations that arise in Physics} 
    \end{center}
    \vskip1cm%
    \begin{center}
        \Large\textit{By}\\
        \vskip0.25cm%
        \LARGE\textbf{\href{mailto:Firstname.Lastname@uon.edu.au}{Juan Silverio Martínez Baena} \\\small{MSc Fisymat(UGR), BSc Physics \& BSc Mathematics (UGR)}}

    \end{center}
    \vfill
    \begin{center}
             \vskip0.2cm
         \Large\textit{Thesis supervised by }\\
         \vskip0.2cm
            \Large\textbf{Dr. Salvador Villegas Barranco.}\\
            \vskip0.2cm
        \textit{Submitted in fulfilment of the requirements for the}
            \vskip0.2cm
          \Large\textbf{PhD Program in Physics and Mathematics}\\
	    \vskip1cm%
           \includegraphics[width=0.3 \columnwidth]{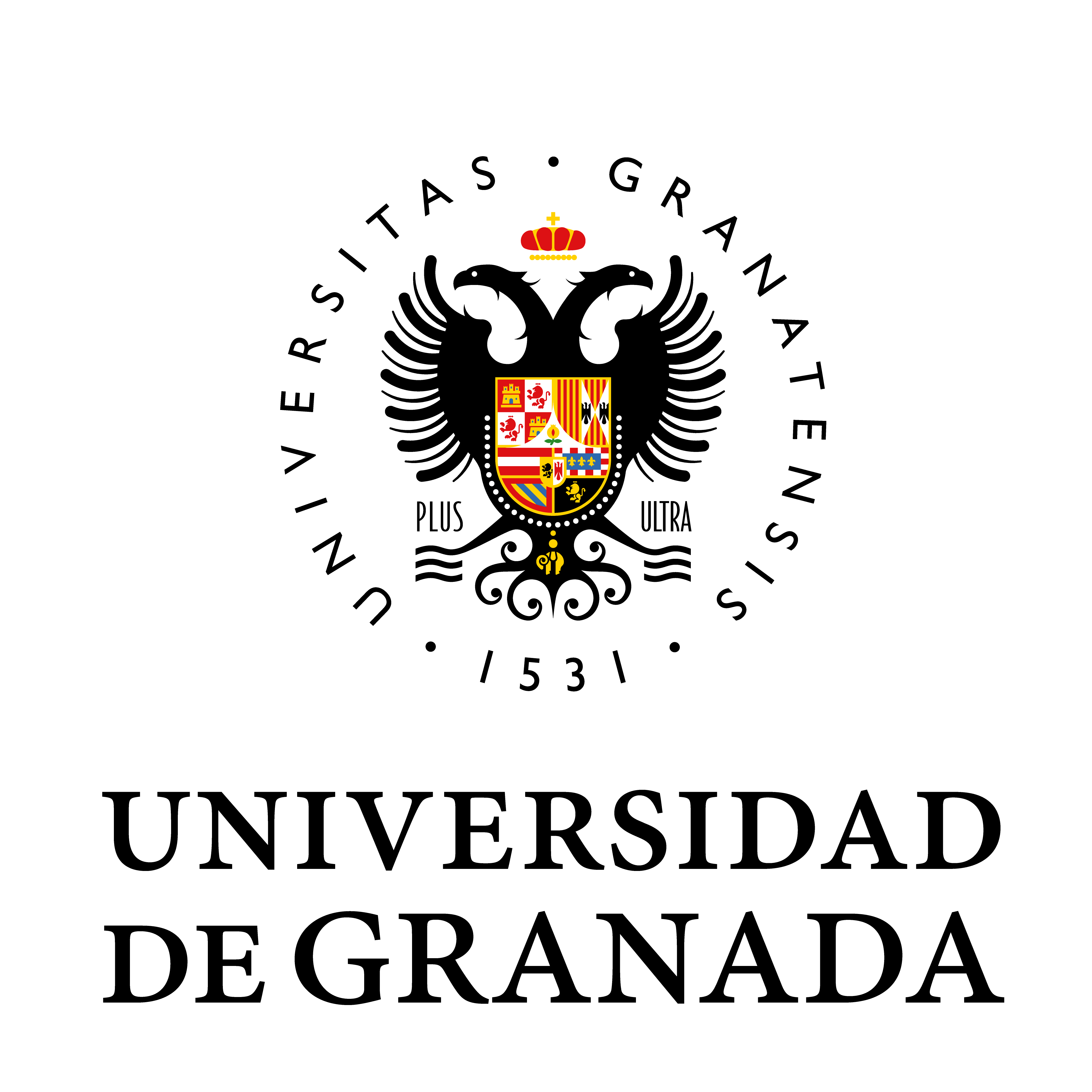}\\
		Departament of Mathematical Analysis, \\
		\href{https://www.ugr.es/}{Faculty of Sciences, University of Granada.}\\
		Avda. de Fuentenueva, s/n C.P. 18071 Granada, Spain.\\
           \today\\
        \vspace{1cm}
    \end{center}
    \vfill
    \newpage
    \thispagestyle{empty}
    \mbox{}
\hypersetup{urlcolor=blue}{}

\onehalfspacing

\chapter*{Agradecimientos}

En primer lugar, quiero expresar mi agradecimiento a mi director, el Catedrático Salvador Villegas Barranco, por el apoyo y la disposición que me ha demostrado durante estos años, ya que sin su ayuda no habría sido posible desarrollar este trabajo. También agradecer al que es mi tutor, el Catedrático David Arcoya, que finalmente decidió darme la oportunidad de contribuir a la producción científica del proyecto del cual en 2020 era investigador principal, mediante mi beca FPI. También quiero agradecer al ex-director del Departamento de Análisis Matemático, el Prof. Juan Carlos Cabello Piñar, por el apoyo que me ha brindado en todo lo que respecta a mi relación con el departamento y sobre todo guiándome al inicio de mi experiencia como docente. De igual manera, agradecezco al Catedrático José A. Carmona de la Universidad de Almería por su gentileza y su inestimable labor gestionando la financiación del proyecto para hacer factible mis actividades académicas externas en congresos y me consta que las de mucha más gente. Por último, agradecer a la Universidad de Granda por hacer posible esta labor investigadora así como al Ministerio de Ciencia, Innovación y Universidades que mediante la ayuda FPI BES-2019-090682 otorgada al proyecto PGC2018-096422-B-I00 ha hecho posible mi sustento durante esta intensa etapa junto con el proyecto actual PID2021-122122NBI00. 

En segundo lugar también quiero agradecer a la Universidad Federal de Sao Carlos UFSCar, en particular al Prof. Franciso Odair, por su incalculable generosidad, no sólo aceptando mi propuesta de estancia en Brasil, sino por el modo en que me acogió, dando lugar a agradecer también a su mujer e hijos por la gran hospitalidad y el buen trato que me dieron. Sin olvidar al grupo de compañeros que hicieron mi estancia productiva tanto a nivel profesional como vital, el maestro, Gustavo, Igor, Eudes, Braulio (incluyendo experiencias menos productivas), Maria Rosilene, Juan Pablo (con interminable tardes de discusión matemática), sin olvidar Hayde y Fabio y sus adorables hijos. Un recuerdo especial para Rodiak y Diana que se convirtieron casi en mis hermanos mayores brasileños aún siendo peruanos.

En tercer lugar, agradecer a la gente local, mis compañeros de Física teórica José M. P. Poyatos, Chema, y nuestras discusiones sobre simetrías y grupos de Lie (entre otras conversaciones y frases célebres) así como a Gerardo G. M. que con su ímpetu me empuja a mantener la motivación por la física teórica y la investigación en general. También especialmente a mi compañero y amigo Alejandro J. C., alias ``el Himene", con el que he compartido muchos ratos de brainstorming y reflexiones sobre la vida y que ha influido en gran medida en parte de esta tesis. 
No quiero olvidar a mis compañeros del IMAG, Fidel (del que agradezco su pasión por la física, no tan común en un matemático puro y la profunda vision de las ideas que hemos compartido), David, Jorge, Diego,...

En cuarto lugar quiero agradecer a mis amigos, Paco, Roberto y Pedro, que han sido un gran apoyo y con los cuales he podido contar para dejar el estrés atrás apoyándome cada fin de semana hablando de lo que me preocupaba y del sentido que le queremos dar a nuestra vida.  

Finalmente quiero agradecer a mi familia, empezando por mi amigo Jesús ya que aunque no nos veamos todo el tiempo que nos gustaría, sigue manteniéndose estoicamente presente. Mis hermanos Rosana y Gonzalo, que son las dos personas, cada uno a su manera, que mejor me podrán entender en muchos aspectos importantes de la vida. Y a mis padres, que con sus imperfecciones me han dado todo lo que tengo a día de hoy.

En especial quiero agradecer a mi compañera de vida estos últimos 6 años, Claudia, por aguantarme en mis peores momentos y continuar aún así iluminando mi camino en la oscuridad y a mi padre que me ha estado apoyando en toda mi carrera académica, ya que sin su sacrificio en momentos clave, no habría llegado a donde estoy.

\chapter*{}
\begin{center}
\textit{{\large A mi padre, por ir más allá...}}
\end{center}


\chapter*{Resumen en español}

El estudio de las propiedades cualitativas (existencia, unicidad/multiplicidad, regularidad, etc) de soluciones de la ecuación de Poisson no lineal $-\Delta u =f(x,u)$ con condición de Dirichlet en un dominio dado (EPnLD), es un problema clásico que admite una cantidad vasta de literatura y que sigue siendo objeto de investigación puntera a nivel internacional junto con las distintas generalizaciones a operadores no lineales, no locales y con condiciones distintas. 

En el estado actual de la comunidad investigadora de matemáticas, en particular, la comunidad de ecuaciones en derivadas parciales, el estudio de la regularidad de soluciones ha destacado por su interés y en ocasiones, por su dificultad. El concepto de estabilidad de solución se ha utilizado en las dos últimas décadas como herramienta para garantizar la regularidad de las soluciones de ecuaciones de tipo elíptico. Una solución es estable cuando la segunda variación del funcional Energía asociado a la ecuación, es semidefinida positiva. Por otro lado, también se ha utilizado como medida de la estabilidad de una solución, el conocido como índice de Morse que viene dado, grosso modo, por la dimensión del mayor subespacio de funciones test en el que la segunda variación de la energía es definida negativa.

Todo ese marco teórico ha dado sus frutos, como han demostrado los autores de la publicación \cite{cabre-figalli-rosoton-serra} cuyo trabajo acabó por demostrar, al mismo tiempo, la regularidad Hölder de las soluciones estables de (EPnLD) en el rango óptimo de dimensiones $N\leq9$, así como la conjetura de Brezis y Vázquez sobre la regularidad de la solución extremal. Estos resultados no han sido generalizados aún a la situación de soluciones de índice de Morse finito (obviamente de las de índice de Morse infinito no se espera nada). Aún así existen resultados parciales como \cite{bahri-lions}
para nolinealidades subcríticas y \cite{figalli-zhang} para nolinealidades supercríticas, que son capaces de acotar las soluciones mediante una constante que esencialmente sólo depende de la dimensión y del índice de Morse. 

El trabajo de esta tesis se ha enmarcado en gran parte en este ámbito. En el Capítulo \ref{ch-counterexample}, se desarrolla un contraejemplo para la que hemos llamado ``Conjetura de Brezis-Vázquez extendida", esto es, la que concierne a soluciones de índice de Morse finito. En particular, demostramos que la acotación del índice de Morse radial, a pesar de lo que uno pueda desear, no aporta la suficiente restricción como para evitar las soluciones singulares. Encontramos una sucesión de soluciones a problemas del tipo (EPnLD) cuyo índice de Morse radial es necesariamente igual a $1$ y el cociente de normas $||\cdot||_{p}/||\cdot||_{q}$ diverge para todo valor $1\leq q<p\leq\infty$ con $p>N/(N-2)$ (la restricción $p>N/(N-2)$ es un hecho inesperado y quizás solamente técnico) en dimensiones $3\leq N\leq 9$.

En el Capítulo \ref{ch-4} se establece la regularidad de soluciones radiales de una versión no autónoma, la ecuación de Hardy-Hénon: $-\Delta u=|x|^\alpha f(u)$, determinando el rango óptimo de dimensiones $2\leq N< 10+4\alpha$ con $\alpha>-2$.

También es tratado en el Capítulo \ref{ch-5} el tema de la existencia y la multiplicidad de soluciones de la ecuación $-\Delta u=g(u)-h(x)f(u)$ donde $h\in L^\infty(\Omega)$ es no negativa, no cero y satisface $h>0  \mbox{ a.e. en } \Omega \backslash \Omega_0  \mbox{ con } \Omega_0=\text{interior}(\{x\in \Omega\, /\, h(x)=0\}),  \ \ {\rm y} \ \ \, |\Omega_0|>0
$. La propiedad de anulación de $h$ hace que este problema haya sido menos estudiado que otros en la literatura. En este trabajo, se demuestra la existencia de una solución positiva y una negativa, así como de hasta cuatro soluciones no triviales distintas, en el marco en que $g$ es asintóticamente lineal y bajo un conjunto apropiado de hipótesis sobre la no linealidad $f$. El factor clave es la posición relativa de la pendiente asintótica $\lambda=g'(0)$ respecto al espectro del operador autoadjunto asociado a una forma cuadrática definida en el dominio $\Omega_0$, un aspecto ya señalado en \cite{AlamaLi1992}. Las herramientas principales utilizadas son el método de sub y supersoluciones en forma integral para encontrar las dos soluciones de signo constante, el Lema de Paso de Montaña para obtener una tercera, y los grupos críticos asociados a todas ellas para identificar y distinguir una cuarta solución (véase el Teorema \ref{th:chang-multiplicity} o \cite[Teorema 3.5]{chang}).

Finalmente, en el Capítulo \ref{ch-upcoming-research} se describen algunas de las ideas que podrían desarrollarse en investigaciones futuras. Nuestros resultados pueden generalizarse en diversas direcciones. En particular, la validez del resultado presentado en el Capítulo \ref{ch-counterexample} como posible contraejemplo de la versión ampliada al índice de Morse completo de la conjetura de Brezis–Vázquez permanece como una cuestión abierta. Asimismo, la posible aplicabilidad de las técnicas desarrolladas en \cite{cabre-figalli-rosoton-serra} al estudio de la regularidad de soluciones estables de la ecuación de Hardy–Hénon no autónoma en dominios convexos también merecería ser explorada.

Los conceptos de estabilidad tratados en el desarrollo de esta tesis, se pueden aplicar al estudio de la estabilidad en teorías de campos de la física teórica. La convergencia de la serie perturbativa y la diferencia entre el orden lineal y el orden superior son algunos de los problemas que surgen en este tema. Hoy en día, cualquier teoría física que se considere seria debe superar pruebas rigurosas de estabilidad. Sin embargo, aún existe cierta carencia de procedimientos precisos que garanticen la detección de todas las formas indeseables de inestabilidad incompatibles con nociones básicas que debe cumplir cualquier teoría razonable desde el punto de vista físico. Sería un avance significativo poder formular definiciones y teoremas precisos que permitan controlar inestabilidades como el acoplamiento fuerte. Un ejemplo representativo es el caso de la gravedad cúbica de Einstein, tal como se expone en \cite{beltran-jimenez2021}.

\clearpage\thispagestyle{empty}

\chapter*{Abstract}

The study of qualitative properties, such as existence, uniqueness/multiplicity, and regularity of solutions to the nonlinear Poisson equation $-\Delta u =f(x,u)$ subject to Dirichlet boundary conditions on a given domain (NPDP), constitutes a classical problem with a vast body of literature. It continues to be an active area of research at the international level, alongside its various generalizations to nonlinear, nonlocal operators and alternative boundary conditions.

Within the current landscape of the mathematical research community, particularly in the field of partial differential equations, the analysis of regularity properties of solutions has garnered significant attention, both for its theoretical interest and the inherent difficulties it presents. Over the last two decades, the notion of solution stability has emerged as a powerful tool to guarantee regularity in elliptic equations. A solution is said to be stable if the second variation of the associated energy functional is nonnegative definite. Alternatively, the stability of a solution can be quantified via its Morse index, which, roughly speaking, is defined as the dimension of the largest subspace of test functions on which the second variation is negative definite.

This theoretical framework has yielded profound results. A notable example is the work of Cabré, Figalli, Ros-Oton, and Serra \cite{cabre-figalli-rosoton-serra}, who simultaneously established the Hölder continuity of stable solutions to (NPDP) in the optimal range of dimensions $N\leq9$, as well as definitively solved a conjecture by Brezis and Vázquez concerning the regularity of the extremal solution. These results have not yet been generalized to solutions with finite Morse index (as opposed to infinite Morse index, for which no such regularity is expected). Nevertheless, there exist partial results: for example, \cite{bahri-lions} for subcritical nonlinearities, and \cite{figalli-zhang} for supercritical ones. These works manage to obtain a priori bounds for the solutions in terms of constants that depend only on the dimension and the Morse index.

A substantial portion of this thesis is situated within this research direction. In Chapter \ref{ch-counterexample}, we construct a counterexample to what we refer to as the “Extended Brezis–Vázquez Conjecture” namely, the conjecture regarding solutions with finite Morse index. In particular, we demonstrate that bounding the radial Morse index does not, contrary to possible expectations, suffice to preclude the occurrence of singular solutions. We construct a sequence of solutions to problems of the form (NPDP) whose radial Morse index is necessarily equal to one, and for which the quotient $||\cdot||{p}/||\cdot||{q}$ diverges for every pair $1\leq q<p\leq\infty$ with $p>N/(N-2)$ (the constraint $p>N/(N-2)$ appears to be unexpected and possibly of technical nature), in dimensions $3\leq N\leq 9$.

In Chapter \ref{ch-4}, we establish regularity results for radial solutions of a non-autonomous version of the equation namely, the Hardy–Hénon equation: $-\Delta u=|x|^\alpha f(u)$ determining the optimal range of dimensions $2\leq N\leq 10+4\alpha$, assuming $\alpha>-2$. The case $\alpha\leq-2$ is much more pathological and there are even some Liouville results for general weak solutions (see \cite{dancer-du-guo}).

In Chapter \ref{ch-5}, we also address the problem of existence and multiplicity of solutions to the equation $-\Delta u=g(u)-h(x)f(u)$ where $h\in L^\infty(\Omega)$ is a non-negative, nontrivial function satisfying
$h>0  \mbox{ a.e. in } \Omega \backslash \Omega_0  \mbox{ with } \Omega_0=\mbox{interior}(\{x\in \Omega\, /\, h(x)=0\}),  \ \ {\rm and} \ \ \, |\Omega_0|>0
$. Due to the vanishing property of $h$, this class of problems has been less addressed in the literature. In this work, we prove the existence of both a positive and a negative solution, as well as up to four distinct nontrivial solutions, under the assumptions that $g$ is asymptotically linear and that the nonlinearity $f$ satisfies a suitable set of conditions. A crucial point is the relative position of the asymptotic slope $\lambda=g'(0)$ with respect to the spectrum of the self-adjoint operator associated to a quadratic form defined on the subdomain  $\Omega_0$, a key aspect already highlighted in \cite{AlamaLi1992}. The main tools employed include the sub- and supersolution method in its integral formulation to obtain two sign-constant solutions, the Mountain Pass Lemma to find a third, and the critical groups associated to all of them to identify and distinguish a fourth solution (see Theorem \ref{th:chang-multiplicity} or \cite[Theorem 3.5]{chang}).

Finally, Chapter \ref{ch-upcoming-research} outlines several avenues for future research. The results obtained herein admit various generalizations. In particular, whether the counterexample in Chapter \ref{ch-counterexample} genuinely invalidates the extended version of the Brezis–Vázquez conjecture for arbitrary finite Morse index remains an open question. Additionally, the potential applicability of the techniques developed in \cite{cabre-figalli-rosoton-serra} to the study of the regularity of stable solutions of the non-autonomous Hardy–Hénon equation on convex domains is a promising direction worth further exploration.

The notions of stability explored in the present thesis may also find applications in the study of stability phenomena in theoretical physics and field theory. Issues such as the convergence of the perturbative series and the discrepancy between linear and higher-order approximations naturally arise in this context. At present, any physically viable theory must withstand stringent stability criteria. Nonetheless, a lack of precise methods for ruling out all physically undesirable instabilities persists. Significant progress would be achieved by formulating rigorous definitions and theorems capable of detecting and ruling out such instabilities—such as strong coupling. A representative case is that of Einstein cubic gravity, as discussed in \cite{beltran-jimenez2021}.

\clearpage\thispagestyle{empty}
\onehalfspacing

\mainmatter


\begin{savequote}[0.75\linewidth]
\textit{`See, in all our searching, the only thing we’ve found that makes the emptiness bearable...is each other''.}\\
Film: Contact (1997), based on Carl Sagan’s novel.
\end{savequote}

\chapter*{Preface}

In the relatively short history of humankind, we have always been looking for answers to the world we see around us. The big questions like where we come from, how everything started, where do we go, or why something exists instead of nothing seem to arise already in the very early stages of our development as a species and being humble, we must say that we are almost equally far from any answer\footnote{In the novel ``The Hitchhiker's Guide to the Galaxy" the answer to this question is 42 but the author Douglas Adams chose that number because a quick poll among his friends revealed that it was very boring.}. With the passage of time, only a few ideas have deeply penetrated our mind beyond artificial (but necessary) progressive specialization and compartmentalization of knowledge into distinct disciplinary domains, shaping the way we think today. The binary way of thinking, the fact that every idea has a contrary in opposition, concepts of good and evil, love and hate, heaven and hell, God and Devil are good examples. The concept of God is possibly very related to the concept of infinite, which is in opposition to finiteness or in opposition to what we can measure, see with our eyes, or even count. Natural numbers seem to have come into place very early; they are also in this range of ideas. All of these ideas seem to help us to describe and understand our life experience so that we can live with some peace of mind. Once one is able to know all the necessary information on a given aspect of life, the mind achieves a state of calm. In that situation, we sometimes say that we have reached the much-desired stability. There is a whole range of different uses of the word stability, from the point that we have described to the mathematical concept that is going to be developed in this thesis. Each of them seems to share intuitively a common pattern that points in the direction of resting in some minimum of energy, understood in the broadest possible sense of the word energy, protected from disturbances or perturbations. The idea of stability seems to me to be one of those ideas as well.


From the way geometry and weight distribution influence the design and construction of an aircraft, and the prediction of weather conditions that pilots must consider, to the discovery of gravitational waves resulting from the merger of massive black holes in the past decade—not to mention the orbits of planets and the balance in the struggle for species survival in nature—the concept of stability plays one of the broadest roles, maybe even in a more notorious way, in science. Most of the time, it defines a criterion for any methodical description or mathematical model of an aspect of reality to be considered acceptable, plausible, and taken seriously.

The common framework from which every scientific particular definition of stability for a system arises can be expressed in terms of an optimization principle: There exists a quantity, sometimes referred to as energy, whose values largely determine the properties of the system and its behavior over time. Any of the different possible properties of the system are called configurations of the system. The central question is then posed: If a system is in a specific configuration with a certain energy and this configuration is slightly perturbed, does this perturbation drastically alter the energy and thus the behavior of the system?, or does the system remain largely unchanged and eventually return to its unperturbed state? If the stability condition holds for the unperturbed configuration, then in some sense that configuration should minimize the energy making the perturbed system underprivileged in terms of efficiency will return to the stable configuration, perhaps after undergoing a transient relaxation process. The majority of this thesis is devoted to the study of qualitative properties of stable solutions of elliptic partial differential equations, which often arise when considering stationary configurations or equilibria in dynamical models of nature coming from Physics.

\newpage

\begin{center}
    \Large{List of symbols}
\end{center}

\begin{longtable}{@{}ll@{}}
\textbf{Symbol} & \textbf{Description} \\
\hline
$N$ & Dimension of the ambient space or manifold \\
$\Omega$ & Smoothly bounded domain in $\mathbb{R}^N$ (open and connected set) \\
$|\Omega|$ & Lebesgue/Haussdorf measure of $\Omega$\\
$B_1$, $\omega_N$ & Unit ball in $\mathbb{R}^N$ with $|B_1|=\omega_N/N$\\
$\mathscr{H}_{N-1}$ & ($N-1$)-dimensional Haussdorf measure\\
$2^* = \frac{2N}{N-2}$ & Critical Sobolev embedding exponent (for $N > 2$) \\
$\lambda_k$ ($\lambda_k(L)$) & $k$-th eigenvalue of the Laplacian (of $L$) with Dirichlet condition \\
$u$ & Solution of an equation\\
$J$, $E(u)$ &A general functional, Energy functional associated to an equation \\
$C_q(J, u)$ & Critical group at a solution $u$ for a functional $J$, in degree $q$ \\
$c,C,C', \kappa, K,K'$ & Constants (with subindices for depence on parameters)\\
$\varphi$ & Test function (in $C^1_c(\Omega)$ or $C_c^\infty(\Omega)$) \\
$d E$, $d^2 E$ & First and second variation of a functional $E$ \\
$H$ & Mean curvature of a hypersurface \\
$\nu$ & Unit normal vector to a hypersurface \\
$\mathcal{SC}$ & Simons Cone\\
$\Phi$ & Field (field theory, theoretical physics)\\
$S(\phi)$ & Action functional in field theory \\
$\mathcal{L}$ & Lagrangian density, $\mathcal{L}(\Phi, \partial \Phi,x)$ \\
$\mathcal{H}$ & Hamiltonian functional (total energy in mechanics or field theory) \\
$\mu, \nu, \rho, \sigma$ & Spacetime indices (range from $0$ to $N-1$ spacetime) \\
$a, b, c$ & Indices for internal symmetries\\
$A,B$ & General multi-index (e.g., for derivatives or fields) \\
$\lbrace x^\mu:\mu=1,\dots,N\rbrace$ & Coordinate system\\
$D_\mu$ & Covariant derivative with respect to coordinate $x^\mu$ \\
$\Gamma^\rho_{\mu\nu}$ & Christoffel symbols of the Levi-Civita connection \\
$T^{\mu\nu}$, $T^\mu\,_\nu$, $T_{\mu\nu}$ & Stress-energy-momentum tensor \\
$M$ & General manifold (smooth, possibly with a metric $g$) \\
$\sqrt{|g|}$ & Invariant measure for metric $g$\\
$(\mathcal{M},g)$, ${M}_g$ & Spacetime $\mathcal{M}$ with Lorentzian metric $g$\\
$\mathbb{R}^{1,3}$ & Minkowski spacetime\\
$G$ & Lie group\\
$\mathfrak{g}$ & Lie algebra of $G$

\end{longtable}

\tableofcontents

\chapter{Introduction}\label{ch-intro}

The differential equation of greatest importance in this thesis is the semilinear Poisson equation \eqref{main_problem}, with Dirichlet boundary conditions imposed on a smoothly bounded domain $\Omega\subset\mathbb{R}^N$.
\begin{equation}\label{main_problem}
\left\lbrace
\begin{array}{rcll}
-\Delta u &=&  f(u)&\text{ in }\Omega\\
u&=&0&\text{ on }\partial \Omega
\end{array}\right.\tag{P1}
\end{equation}
where $\Delta:=\sum_{i=1}^N\partial^2/\partial x_i^2$ is the Laplace operator and the one variable function $f$ is the so-called nonlinearity, which can have a variety of different regularity hypothesis depending on the problem. We will refer to the first equation in \eqref{main_problem} as \textit{the equation} of \eqref{main_problem} and the second equation as \textit{the (Dirichlet)-boundary condition} of \eqref{main_problem}. A solution of the whole Dirichlet problem will be simply referred to as a solution of \eqref{main_problem}. We will assume that $f\in C^1(\mathbb{R})$ if nothing specific is said.\\
Rivers of ink have been spilled over this equation, which continues to play a central role in the study of partial differential equations. The primary goal is to establish qualitative properties of its possible solutions, such as existence, uniqueness or multiplicity, and regularity. These properties critically depend on the nature of the nonlinearity and the domain under consideration. Consequently, it is generally necessary to impose additional assumptions on the properties (regularity, geometry, symmetries, etc.) of both the domain and the nonlinearity in order to obtain meaningful results. One of the fundamental ideas and advices in partial differential equations is to try first radially symmetric solutions. This is widely used in Chapters \ref{ch-4} and \ref{ch-5}. Often, studying this kind of solutions is enough to discover the optimal coefficients that are relevant to make a result true in the most general possible setting.

\section{Organisation of the Thesis}
First of all, in Section \ref{sec-background} the notation and some basic definitions and properties are presented; all of them are now of common knowledge in the partial differential equations community but essential for a correct organization and presentation of the new contributions. 

In Chapter \ref{ch-lit-rev} we review the most relevant and recent progress that has been made in the literature concerning the topics of regularity of stable and finite Morse index solutions in Section \ref{sec-regularity-literature} and about multiplicity of solutions to semilinear elliptic problems involving degenerate nonlinearities in Section \ref{sec-multiplicity-literature}.

Chapter \ref{ch-counterexample} is devoted to show the first novel contribution of the thesis \cite{paper1}, which is a counterexample to the possible regularity of the class of finite radial Morse index solutions to semilinear elliptic equations of type \eqref{main_problem}. We construct a counterexample to what can be referred to as an extension of the Brezis-Vázquez conjecture regarding the boundedness of stable $u\in W_0^{1,2}(\Omega)$ solutions of \eqref{main_problem} for nonlinearities that are positive, nondecreasing, convex and superlinear at $+\infty$ in dimensions $N\leq 9$. The extension is meant to assume that the solutions have finite Morse index rather than being stable. Specifically, we construct a sequence of solutions to problems of the type \eqref{main_problem} with a radial Morse index equal to one, whose quotient of the $L^\infty$-norm over its $L^1$-norm is unbounded. From this fact it is deduced that there exits no universal constant $C=C(N)$ depending only on the dimension such that
\begin{equation}
\frac{||u||_{L^\infty(B_1)}}{||u||_{L^1(B_1)}}\leq C
\end{equation}
for solutions that have only finite radial Morse index. If the finiteness of the full Morse index is enough to avoid this kind of counterexamples without any additional assumption remains a (seemly hard) open problem, even though there exist some partial classical results and some other quite recent ones that have been achieved (see Section \ref{sec-background}).

Chapter \ref{ch-4} develop the new results we got about regularity of stable solutions to a non autonomous semilinear elliptic version, the so-called Hardy-Hénon equation, that is available in \cite{arxiv-paper-2}. 

Chapter \ref{ch-5} deals with the existence and multiplicity of solutions for an even more general version of a semilinear elliptic equation that is still a preprint, having some work in progress \cite{arxiv-paper-multiplicity}.

In Chapter \ref{ch-upcoming-research} some of the upcoming research roadmaps or research avenues that arise from the studies of the thesis are described and put in context.

\section{Background}\label{sec-background}

In this section, we introduce the notions that are fundamental for the coherent development of the novel contributions achieved, as well as establish the notation. This text is not intended to be self-contained; instead, it presents only the most relevant definitions, propositions, and theorems for our work, with references provided for additional details.

\subsection{Stability and Morse index in semilinear elliptic equations}


Let $\Omega\subset \mathbb{R}^N$ be a bounded domain. It is well known (see, for instance, \cite{gilbarg-trudinger}) that the Dirichlet problem

\begin{equation}\label{laplacian_eigenvalue_equation}
\left\lbrace
\begin{array}{rcll}
-\Delta u &=&  \lambda u &\text{ in }\Omega\\
u&=&0&\text{ on }\partial \Omega
\end{array}\right.
\end{equation}
determines a sequence of \textit{eigenvalues}, denoted by $\lbrace\lambda_i\rbrace_{i\in\mathbb{N}}$, ordered according to the \textit{min-max principle} such that $0<\lambda_1<\lambda_2 <\lambda_3<\dots$, up to multiplicity, along with a corresponding sequence of finite-dimensional eigenspaces $\lbrace\Lambda_i\rbrace_{i\in\mathbb{N}}$. In particular, the first eigenvalue $\lambda_1$ is simple, and its associated eigenfunction $\phi_1\in\Lambda_1$ is positive in $\Omega$.  When $f\in C^1(\mathbb{R})$, we can Linearize the equation \eqref{main_problem} yielding the operator $-\Delta(\,\cdot\,) -f'(u)(\,\cdot\, ) $ and the Dirichlet eigenvalue problem for this operator is expressed as

\begin{equation}\label{main_eigenvalue_problem}
\left\lbrace
\begin{array}{rcll}
-\Delta v &=&  f'(u)v+\lambda v &\text{ in }\Omega\\
u&=&0&\text{ on }\partial \Omega
\end{array}\right.
\end{equation}
The min-max characterization of the first eigenvalue for this problem is
\begin{equation}\label{first-eigenvalue-characterization}
    \lambda_1(-\Delta-f'(u),\Omega)=\inf_{\varphi\in H^1_0(\Omega),||\varphi||_{L^2(\Omega)}=1}\left\lbrace\int_\Omega|\nabla \varphi|^2 dx-\int_\Omega f'(u)\varphi^2 dx\right\rbrace
\end{equation}
Clearly, there is no guaranty that the first eigenvalue is positive. However, it will be useful for the exposition a few lines below to assume that
\begin{equation}\label{stability-assumption}
    \lambda_1(-\Delta-f'(u),\Omega)\geq 0.
\end{equation}
On the other hand, integrating by parts equation \eqref{main_problem}, we obtain
\begin{equation}\label{weak-solutions}
    \int_\Omega\nabla u \nabla\varphi\, dx=\int_\Omega f(u)\varphi\,  dx\ \ \ \ \ \forall\varphi\in C^1_c(\Omega).
\end{equation}
For this equation to be well-defined, we must assume that $u\in W^{1,1}(\Omega)$ and  $f(u)\in L^1_{loc}(\Omega)$.
\begin{definition}[Weak solution]
    A function $u\in W^{1,1}(\Omega)$ that satisfies \eqref{weak-solutions} is called weak solution of the equation of \eqref{main_problem}. If $u\in H^1_0(\Omega)$, $f(u)\mathsf{dist}(\cdot,\partial\Omega)\in L^1(\Omega)$ and $u$ satisfies \eqref{weak-solutions} for all $\varphi\in C^1(\overline{\Omega})\cap C^2_0(\Omega)$ then is called a weak solution of the problem \eqref{main_problem}.
\end{definition}
After another integration by parts, equation \eqref{weak-solutions} yields
\begin{equation}\label{eq:L1-weak-solutions}
    -\int_\Omega u\, \Delta\varphi\, dx=\int_\Omega f(u)\varphi\,  dx\ \ \ \ \ \forall\varphi\in C^2_c(\Omega).
\end{equation}
\begin{definition}[$L^1$-weak solution]
    A function $u\in L^1(\Omega)$, such that $f(u)\in L^1_{loc}(\Omega)$ and that satisfies \eqref{eq:L1-weak-solutions} is called a $L^1$-weak (or very weak) solution of the equation of \eqref{main_problem}. If $f(u)\,\mathsf{dist}(\cdot,\partial\Omega)\in L^1(\Omega)$ and \eqref{eq:L1-weak-solutions} is satisfied by $u\in L^1(\Omega)$ for all $\varphi\in C^1(\overline{\Omega})\cap C^2_0(\Omega)$ then $u\in L^1(\Omega)$ is called a $L^1$-weak solution of the problem \eqref{main_problem}.
\end{definition}
\begin{remark}\ 
\begin{enumerate}[label=(\roman*)]
    \item From the equation \eqref{eq:L1-weak-solutions} one merely need $u\in L^1_{loc}(\Omega)$. However, it is convenient to maintain the assumption $u\in L^1(\Omega)$ in order to avoid \textit{complete blow-up} phenomena (see \cite[Section 3.1.1]{dupaigne}).

    \item Every classical solution, $u\in C^2(\Omega)$ of the equation of \eqref{main_problem} or $u\in C^2(\Omega)\cap C(\overline\Omega)$ of the problem \eqref{main_problem}, is also a solution in both weak and $L^1$-weak senses.
\end{enumerate}

\end{remark}
Take $\epsilon>0$ and a function $\varphi\in C_c^1(\Omega)$ and let $u$ be a weak solution of the equation of \eqref{main_problem}. Now define the so-called \textit{Energy functional} as follows
\begin{equation}\label{energy-functional}
    E_\Omega(u)=\int_\Omega\left\lbrace\frac{\left|\nabla u\right|^2}{2}-F(u)\right\rbrace dx,
\end{equation}
where $F(u)=\int_0^uf(s)ds$. The behavior of the energy under the perturbation $u+\epsilon\varphi$ of the solution $u$ at first order would be measured by the \textit{first variation} of the energy:
\[\lim_{\epsilon\to 0}\frac{E(u+\epsilon\varphi)-E(u)}{\epsilon}\equiv\left.\frac{d }{d\epsilon}E(u+\epsilon\varphi)\right|_{\epsilon=0}=:dE_u(\varphi)=\int_\Omega\left\lbrace\nabla u \nabla\varphi -f(u)\varphi \right\rbrace dx\]
In this sense, the equation of \eqref{main_problem} is the Euler-Lagrange equation of the functional \eqref{energy-functional}, meaning that weak solutions of the equation of \eqref{main_problem} correspond to critical points of the energy functional (analogously with respect to the problem \eqref{main_problem} with the mentioned additional assumptions in the previous definitions). Assuming that $f'(u)\in L^1_{loc}(\Omega)$, the second variation can also be calculated to obtain
\begin{equation}\label{second-variation-energy}
    Q_{u,\Omega}(\varphi):=d^2E_u(\varphi,\varphi)=\int_\Omega \left\lbrace\left|\nabla\varphi\right|^2-f'(u)\varphi^2\right\rbrace dx\ \ \ \ \ \forall\varphi\in C^1_c(\Omega)
\end{equation}
which is a quadratic form that is positive semidefinite if one take into account the assumption \eqref{stability-assumption}.
\begin{definition}[Stable solution]
    Assume $f'(u)\in L^1_{loc}(\Omega)$. An $L^1$-weak solution $u$ of the equation \eqref{main_problem} is called \textit{stable}\footnote{Some authors prefer to call this solutions semi-stable, and reserve the word stable for the strict inequality in \eqref{def-stability}.} if it satisfies the inequality
    \begin{equation}\label{def-stability}
         Q_{u,\Omega}(\varphi)\geq 0\ \ \ \ \ \forall\varphi\in C^1_c(\Omega)
    \end{equation}
\end{definition}

\begin{definition}[Local stability]
    A solution is locally stable if for any $x\in \Omega$ there exists an open neighborhood $\omega_x$ of $x$ such that $Q_{u,\omega_x}[\varphi]\geq 0$ for all $\varphi \in C_c^1(\omega_x)$.
\end{definition}
Stability is a property that has been proved extremely powerful for many purposes. Stable solutions exhibit a large amount of well-behaved properties (see the monograph \cite{dupaigne} for an extensive treatment) while, at the same time, it is such a natural concept, that we can say that stable solutions are the most expected ones in non chaotic systems in nature. One of the benefits that we will use is the following approximation property:
\begin{proposition}[Theorem 3.2.1, \cite{dupaigne}]\label{prop-approximation-by-stable-solutions}
    Let $u$ be a stable weak solution of \eqref{main_problem} and assume $f\in C^1(\mathbb{R})$, $f\geq 0$ and $f$ is convex. Then there exists a sequence $\lbrace u_\epsilon\rbrace$ of classical solutions to the problem
    \begin{equation}
\left\lbrace
\begin{array}{rcll}
-\Delta u_\epsilon &=& (1-\epsilon) f( u_\epsilon) &\text{ in }\Omega\\
u_\epsilon&=&0&\text{ on }\partial \Omega
\end{array}\right.
\end{equation}
such that $u=\lim_{\epsilon\to 0}u_\epsilon$ in $H^1_0(\Omega)$.
\end{proposition}
Since its first introduction by M. Morse in the early 20th century in \cite{morse}, it is now fairly well established the interest in finite Morse index solutions.
\begin{definition}[Morse index]
    Let $u$ be a solution of \eqref{main_problem} with $f'(u)\in L^1_{loc}(\Omega')$ for a subdomain $\Omega'\subset\Omega$. We say that $u$ has Morse index $\mathsf{ind}(u,\Omega')=k\in \mathbb{N}$ in $\Omega'$ if $k$ is the maximal dimension of a subspace $X_k\subset C^1_c(\Omega')$ such that,
\begin{equation*}
Q_{u,\Omega'}(\varphi) < 0\ \ \ \ \ \ \forall \varphi\in X_k\setminus\lbrace 0 \rbrace
\end{equation*} 
If we restrict the set $C^1_c(\Omega')$ to its intersection with the space of radial functions, say $C^1_{c,rad}(\Omega')$, we get the definition of \textit{radial Morse index} in  $\Omega'$ denoted by $\mathsf{ind_r}(u,\Omega')$. We will omit the reference to the domain in the case $\Omega'=\Omega$.
\end{definition}
\begin{remark}
It is immediate from the definition that,
\begin{itemize}
\item[$i)$] $\mathsf{ind}(u)=0\ \text{in}\ \Omega'\iff u$ is stable in $\Omega'$.
\item[$ii)$] $\mathsf{ind_r}(u,\Omega')\leq \mathsf{ind}(u,\Omega')$.
\end{itemize}
\end{remark}
In general, finite Morse index also guarantees  a number of desirable properties from which many advantages can be gained. A particular example, that we will exploit, is the following:

\begin{proposition}[Chapter 1, \cite{dupaigne}]\label{dupaigne_proposition} Let $\Omega$ be a bounded domain and $u\in H^1_0(\Omega)$ be a  solution of \eqref{main_problem} with $f\in C^1(\mathbb{R})$. Then:
\begin{enumerate}\label{morse_index_properties}

\item \label{local_stability}
If $\mathsf{ind}(u)< +\infty$, then $u$ is locally stable.

\item \label{finite_morse_index_for_regular_solutions}
If $u\in C^2(\overline{\Omega})$, then $\mathsf{ind}(u)< +\infty$ and it is equal to the number of negative eigenvalues of the linearized operator $-\Delta -f'(u)$ (with Dirichlet boundary conditions).

\item $\mathsf{ind_r}(u)=0\iff \mathsf{ind}(u)=0$. \label{radial_morse_index_0_iff_morse_index_0}

\end{enumerate}

\end{proposition}
The Proposition \ref{morse_index_properties} \eqref{radial_morse_index_0_iff_morse_index_0} is a consequence of the well-known increasing behavior of the sequence of eigenvalues of the linearized operator $(-\Delta -f'(u))$ and the fact that its first eigenvalue is simple and for $\Omega=B_1$ the corresponding eigenfunction is radial.

For the sake of our study regularity of stable solutions it will be necessary to recall two classical results. The first one is about the \textit{mass} (the $L^1$-norm) of the solution for nonlinearities that go to $+\infty$ faster than a linearly with slope $\lambda_1$.  
\begin{proposition}[Proposition B.1, \cite{cabre-figalli-rosoton-serra}]
\label{prop:L1}
Let $\Omega\subset \mathbb{R}^N$ be a bounded domain of class {$C^1$},
and let $u\in C(\overline\Omega)\cap C^2(\Omega)$ solve \eqref{main_problem} for some $f: \mathbb{R}\to [0,+\infty)$ satisfying
\begin{equation}
\label{eq:superlin f}
f(t)\geq C_1 t-C_2\qquad \mbox{for all }\,t \geq 0,\quad \text{with $C_1>\lambda_1$ and $C_2\geq 0$,}
\end{equation}
where $\lambda_1=\lambda_1(\Omega)>0$ is the {first eigenvalue of the Laplacian in $\Omega$ with Dirichlet homogeneous boundary condition}.
Then there exists a constant $C$, depending only on $C_1$, $C_2$, and $\Omega$, such that
$$
\|u\|_{L^1(\Omega)}\leq C.
$$
\end{proposition}
\noindent The second one tells us that for almost every admissible nonlinearity, boundedness of the solution implies its smoothness.
\begin{theorem}[Elliptic regularity, \cite{gilbarg-trudinger}]
    Let $u:\Omega\to\mathbb{R}$ be a (very) weak solution of \eqref{main_problem}.
    \begin{enumerate}
        \item  If $1<p<+\infty$, $f(u)\in L^p(\Omega)$, then $u\in W^{2,p}(\Omega)$ (Calderon-Zygmund estimates).
        \item If $0<\alpha<1$, $f(u)\in C^{k,\alpha}(\Omega)$, then $u\in C^{k+2,\alpha}(\Omega)$ (Schauder estimates).
    \end{enumerate}
    Consecuently, by the standard Sobolev embedding theorems \ref{th:sobolev-embeddings}, if $u\in L^\infty(\Omega)$ and $f\in C^\infty(\mathbb{R})$, then $u\in C^\infty(\Omega)$.
\end{theorem}

\subsection{Dynamical stability vs elliptic stability}
The most common intuition behind the concept of stability regarded as the ability of a system to return to the equilibrium (i.e., the configuration of minimum energy, after a small perturbation), is not directly related with the Definition \ref{def-stability}, since it is made from a dynamical point of view, as we are accustomed to seeing in our daily experience. The formal definition of dynamical stability, expressed from an elliptic perspective could be stated as follows (see \cite[Section 1.4]{dupaigne}).

\begin{definition}
Let $\Omega\in\mathbb{R}^N$ be a smooth bounded domain. A solution $u\in C^2(\overline{\Omega})$ of \eqref{main_problem} is said \textit{asymptotically stable} if there exists $\epsilon>0$ such that for all $u_0\in C^2(\overline{\Omega})$ with $||u-u_0||_{L^\infty(\Omega)}<\epsilon$, the solution $v\in C^2(\overline{\Omega}\times [0,T])$ of,
\begin{equation}\label{eq:dyn_stab}
\left\lbrace
\begin{array}{rcll}
\frac{\partial v}{\partial t}-\Delta v &=&  f(v)&\text{ in }\Omega\times (0,T)\\
v&=&0&\text{ on }\partial \Omega\times (0,T)\\
v(x,0)&=&u_0(x)    &\ x\in\Omega
\end{array}\right.
\end{equation}
is defined for all times $T>0$ and,
\begin{equation}
\lim_{t\rightarrow +\infty}||v(t)-u||_{L^\infty(\Omega)}=0
\end{equation}
\end{definition}
An \textit{equilibrium} of \eqref{eq:dyn_stab} is just an stationary solution, i.e., a solution that does not depend on time $\partial_tv=0$. The previous definition can be better understood in the context of \textit{Lyapunov stability}, where $u$ can be regarded as any equilibrium of the parabolic equation \eqref{eq:dyn_stab}. An equilibrium is called just \textit{Lyapunov stable} if for every $\epsilon >0$ there exists a radius $\delta>0$ such that $||v(\cdot,t)-u||_{L^\infty(\Omega)}<\epsilon$ for all $t\geq 0$.

 In the dynamical case, the perturbations will evolve over time. The idea of the last definition focuses on shifting the initial conditions slightly enough, while maintaining the equation of motion unchanged. Besides the perturbation is given only at a specific instant in time, it fully accounts for the nonlinear features of the equation. 

It is interesting to compare this nonlinear notion of stability with the linearized stability that comes from \eqref{stability-assumption}. In other words, it is legitimate to ask if any stationary solution which is stable under the elliptic definition is also asymptotically stable and vice versa.

\begin{proposition}[Proposition 1.4.1 and 1.4.2, \cite{dupaigne}]
Assume $\Omega$ smoothly bounded domain of $\mathbb{R}^N$ and $\lambda_1(-\Delta-f'(u);\Omega)\not= 0$. Then, the solution $u$ of \eqref{main_problem} is stable, if and only if, it is asymptotically stable. 
\end{proposition}
  
In the case $\lambda_1(-\Delta-f'(u);\Omega)= 0$, no conclusive answer can be stated and there are positive and negative examples.  

\begin{enumerate}
    \item Example (\cite[Proposition 1.4.3]{dupaigne}): Let $\Omega\subset\mathbb{R}^N$ be a smooth bounded domain such that the first eigenvalue $\lambda_1(-\Delta,\Omega)>0$ and $\phi_1$ its first eigenfunction, $\epsilon>0$, $1<p<2^*-1$ and consider the following equation
    \begin{equation}\label{linear-vs-nonlinear1}\tag{A1}
        \left\{
        \begin{aligned}
        \frac{\partial v}{\partial t} - \Delta v &= \lambda_1 v+|v|^{p-1}, 
           && \text{in } \Omega \times (0,T),\\
                v       &= 0, 
           && \text{on }\partial\Omega \times (0,T).\\
        v(x,0)       &= \pm \epsilon \,\phi_1(x), 
           && \text{for }x\in \Omega.
        \end{aligned}
        \right.
    \end{equation}
    Then, $u=0$ is a solution of \eqref{main_problem} with $f(u)=\lambda_1 u+|u|^{p-1}$ such that $\lambda_1(-\Delta-f'(u);\Omega)= 0$. The energy functional $E_\Omega:H_0^1\rightarrow\mathbb{R}$ associated with the elliptic equation is
    \begin{equation*}
        E_\Omega (u)=\int_\Omega\frac{\left|\nabla u\right|^2}{2}dx-\frac{\lambda_1}{2}\int_\Omega u^2dx-\frac{1}{p+1}\int_\Omega|u|^{p+1} dx,
    \end{equation*}
    The solution of \eqref{linear-vs-nonlinear1} satisfies $E(v)\leq E(\pm\epsilon\,\phi_1)<0$ while $E(0)=0$ so $u=0$ is not asymptotically stable.

    \item Example (\cite[Proposition 1.4.4]{dupaigne}): With the same assumptions of the previous example but with $p>1$ and consider the following equation
    \begin{equation}\label{linear-vs-nonlinear2}\tag{A2}
        \left\{
        \begin{aligned}
        \frac{\partial v}{\partial t} - \Delta v &= \lambda_1 v-|v|^{p-1}v, 
           && \text{in } \Omega \times (0,T),\\
                v       &= 0, 
           && \text{on }\partial\Omega \times (0,T).\\
        v(x,0)       &= u_0(x), 
           && \text{for }x\in \Omega.
        \end{aligned}
        \right.
    \end{equation}
    for $||u_0||_{\infty}$ small enough. Then, $u=0$ is a solution to \eqref{main_problem} with $f(u)=\lambda_1 u-|u|^{p-1}$ such that $\lambda_1(-\Delta-f'(u);\Omega)= 0$. In addition, one can check that given $\phi_1$ an eigenfunction associated to $\lambda_1$, $\underline{u}_n=-\phi_1/n$ and $\overline{u}_n=\phi_1/n$ are a sub and a supersolution respectively such that $||\underline{u}_n-\overline{u}_n||_{L^\infty(\Omega)}\rightarrow 0$ and $\underline{u}_n<v<\overline{u}_n$ uniformly in time. Thus, $u=0$ is asymptotically stable.
\end{enumerate}

\subsection{Infinite dimensional Morse theory}\label{subsec-morse-theory}
This section is dedicated to establish the notation and state the results that will be used in the development of Chapter \ref{ch-5} of this thesis. It focuses on infinite-dimensional Morse theory applied to the search for critical points of functionals, for which we follow the book of K.C. Chang \cite{chang}. This theory serves as a bridge between several distinct branches of mathematics: Algebra and Topology, or more precisely, Algebraic Topology, for the study of singular homology in its abstract form, and Functional Analysis and Differential Equations, which provide the conditions under which this homology can be applied to the study of nontrivial ``holes" (i.e., critical points) in the level sets of the functionals under consideration.\\
The main mathematical object to present is the \textit{critical group} of a functional in a certain point of its domain. In the following, few algebraic notions are reviewed to arrive at its definition.
Let $\Delta_q\subset\mathbb{R}^q$ be the q-dimensional \textit{simplex} (in dimension $3$ a tetrahedron), i.e.,
    $$\Delta_q:=\left\lbrace\sum_{j=0}^q \lambda_je_j|\lambda_j\geq 0,\sum_{j=0}^q\lambda_j=1\right\rbrace$$
and let $X$ be a topological space. Any continuous map $\Delta_q\to X$ is called a singular q-simplex. The set of all singular q-simplexes is denoted by $\Sigma_q$.
Given an Abelian group $G$, one consider formal linear combinations, $\sigma=\sum g_i\sigma_i$ with $g_i\in G$ and $\sigma_i\in \Sigma_q$, these are called \textit{singular $q$-chains} and we will denote by $C_q(X,G)$ the set of all singular q-chains. Given a \textit{topological pair} $(X,Y)$ (i. e., $Y\subset X$ as a topological subspace), the set 
$$C_q(X,Y,G):=C_q(X,G)/C_q(Y,G),$$
is called the \textit{singular q-relative chain module}. A fundamental operator for defining homology is the \textit{boundary operator}
$$
\begin{array}{rcll}
\partial :C_q(X,G)&\to&  C_{q-1}(X,G)\\
\sigma&\mapsto&\sum_{j=0}^q (-1)^j\sigma^{(j)},
\end{array}
$$
where $\sigma^{(j)}$ stands for the singular $(q-1)$-simplex generated by omitting the $j^{th}$-vector in the linear combination. Is is not difficult to see that this operator satisfies $\partial^2=0$ and that can be induced to the quotient of the singular $q$-relative chain module, i.e., such that the diagram
\[ \begin{tikzcd}
C_q(X,G) \arrow{r}{\pi_q} \arrow[swap]{d}{\partial} & C_q(X,Y,G) \arrow{d}{\overline{\partial}} \\%
C_{q-1}(X,G) \arrow{r}{\pi_{q-1}}& C_{q-1}(X,Y,G)
\end{tikzcd}
\]
is commutative ($\pi$ is the natural projection). Here we can stress the main object under study in homology:
\begin{definition}[Homology groups]
    Let $Z_q(X,Y,G):=ker(\overline{\partial})$ and $B_q(X,Y,G):=Im(\overline{\partial})$. \textit{The singular $q$-relative homology module} is
    \begin{equation}
        H_q(X,Y,G):=Z_q(X,Y,G)/B_q(X,Y,G).
    \end{equation}
    The rank of $H_q(X,Y,G)$ is called the \textit{singular $q$-Betti number}. In the case $Y=\emptyset$ we write $H_q(X,G)$. 
\end{definition}
\begin{remark}
    In the case $G$ is a field $Q$, we have that $\mbox{rank}\, H_q(X,Y,Q)=\dim\, H_q(X,Y,Q)$. The quantity $\chi(X,Y,Q)=\sum_{q=0}^\infty (-1)^q\dim\, H_q(X,Y,Q)$ or in particular when $Y=\emptyset$ $\chi(X,G)$, is called the \textit{Euler characteristic}.
\end{remark}
For $\mathbb{S}^N$ being the $N$-dimensional sphere and for the $N$-dimensional unit ball $B_1^N$ the following homology groups are known:
\begin{equation}\label{homology-sphere}
H_q(\mathbb{S}^N,G)\cong\left\lbrace
\begin{array}{llrl}
0 &q\not=N\ \   q,N\geq 1, \\
G &q=N \ \ q,N\geq 1 & \vee & q=0, N\geq 1,\\
G^2 &q=N=0.&& 
\end{array}\right.
\end{equation}
\begin{equation}\label{homology-ball}
H_q(B_1^N,\mathbb{S}^{N-1},G)\cong\left\lbrace
\begin{array}{lr}
0 &q\not=N,\\
G &q=N  .
\end{array}\right.
\end{equation}
Now consider a Banach manifold $M$\footnote{Technically, one needs a Finsler structure on $M$ in order to define a pseudo-gradient vector field that allows us to use some deformation lemmas to detect critical points (see \cite[Ch. I, Sec. 3]{chang}).} and a $C^1$-functional $J:M\to\mathbb{R}$. Let $c\in \mathbb{R}$ and denote its $c$-sublevel set by $J_c:=\lbrace u\in M|J(u)\leq c\rbrace$. Let $K:=\lbrace u\in M|dJ_u(\varphi)=0\ \ \forall\varphi\in C^\infty(M)\rbrace$ be the set of critical points of $J$, $u\in K$ and let $c=J(u)$.
\begin{definition}[Critical group]
We call the $q^{th}$-critical group of $J$ at $u$ to the homology group
\begin{equation}\label{def-critical-group}
    C_q^G(J,u):=H_q(J_c\cap U_u,(J_c\setminus\lbrace u\rbrace)\cap U_u,G),
\end{equation}
where $U_u$ is any neighborhood of $u$ such that $K\cap(J_c\cap U_u)=\lbrace u\rbrace$.
\end{definition}
If the Banach manifold is a Hilbert Riemannian manifold $(M,<,>)$ and $u\in K$ is a \textit{nondegenerate critical point} (i.e. $dJ_u$ has a bounded inverse) then $d^2J_u$ is a self-adjoint operator which posses a resolution of the identity. Morse index can also be defined in this more abstract setting as the dimension of the negative eigenspace of $d^2J_u$ corresponding to the spectral decomposing. Now we have the following key lemma.
\begin{lemma}[Morse Lemma 4.1, \cite{chang}]
    Assume that $J\in C^2(M,\mathbb{R})$, and $u\in M$ is a nondegenerate critical point, then there exists a neighborhood $U_u$ of $u$ and a local diffeomorphism $\Psi:U_u\to T_u(M)$ with $\Psi(u)=0$, such that
    \begin{equation}
        J\circ\Psi^{-1}(\xi)=J(u)+\frac{1}{2}<d^2J_u(\xi),\xi>\ \ \ \ \forall\xi\in\Psi(U_u).
    \end{equation}
\end{lemma}
With this lemma we see that the second variation of the functional gives almost all the information about de functional near some nondegenerate critical point and thus, allows us to compute critical groups of nondegenerate critical points via its Morse index.
\begin{theorem}[Theorem 4.1, \cite{chang}]
    Suppose that $J\in C^2(M,\mathbb{R})$ and $u\in M$ is a nondegenerate critical point of $J$ with Morse index $k$, then
    \begin{equation}
        C_q(J,u)=\left\lbrace
\begin{array}{lr}
G &q=k,\\
0 &q\not=k.  
\end{array}\right.
    \end{equation}
\end{theorem}
In particular we already know that local minimizers are stable (hence its Morse index is zero) which is now also reflected in the following lemma for critical points arising from the sub-supersolution method. We will denote $2^*-1=({N+2})/({N-2})$ where $2^*=N/(N-2)$ is the critical Sobolev exponent for the embedding $W^{1,2}(\mathbb{R}^N)\hookrightarrow L^p(\mathbb{R}^N)$. 
\begin{lemma}[Lemma 2.1, \cite{chang}]\label{lemma-critical-group-minimal-solution}
Let $\underline{u}<\overline{u}$ a pair of strict sub- and supersolutions of \eqref{main_problem} for a possibly non autonomous non linearity $f\in C(\overline{\Omega}\times\mathbb{R},\mathbb{R})$ with the growth condition $|f(x,t)|\leq C(1+|t|^p)$ with $1<p<2^*-1$ for all $t\in \mathbb{R}$. Then, there exists a local minimum $u\in C^1_0(\overline{\Omega})$ of $\widetilde{J}=J|_{C^1_0(\overline{\Omega})}$ such that $\underline{u}<u<\overline{u}$. Moreover, if $u$ is isolated, then
\begin{equation}
    C_q(\widetilde{J},u)=\left\lbrace
\begin{array}{lr}
G &q=0,\\
0 &q\not=0.  
\end{array}\right.
\end{equation}
\end{lemma}
On the other hand it is also known (see \cite{fang-ghoussoub}) that the Morse index of a Mountain Pass type solution is less or equal than one. This is reflected in the following result.
\begin{theorem}[Theorem 1.6, \cite{chang}]\label{th:chang-mountain-pass}
    Assume that $J\in C^2(M)$, that $u$ is a mountain pass solution  of \eqref{main_problem} and that $d^2J_u$ is a Fredholm operator such that if $0\in ker(d^2J_u)$ then $\dim \,ker(d^2J_u)=1$. Then 
    \begin{equation}
        C_q(J,u)=\left\lbrace
        \begin{array}{lr}
            G &q=1,\\
            0 &q\not=1 . 
        \end{array}\right.
    \end{equation}
\end{theorem}
All these characterizations of critical groups under certain assumptions can be applied to find and distinguish different critical points of the energy functional of semilinear elliptic equations that allow us to achieve some existence and multiplicity results. The following result, makes use of everything mentioned above and is in the line of what we present in Chapter \ref{ch-5}.
\begin{theorem}[Theorem 3.5, \cite{chang}]\label{th:chang-multiplicity}
    Assume that the nonlinearity of \eqref{main_problem} satisfies
    \begin{itemize}
        \item[$i)$] $\lim_{|t|\to\infty}\frac{f(t)}{t}<\lambda_1$, the first eigenvalue of $-\Delta$ with Dirichlet boundary condition.
        \item[$ii)$] $f(0)=0$ and $f\in C^1(\mathbb{R})$.
    \end{itemize}
    and $\lambda:=f'(0)$. Then
    \begin{enumerate}
        \item For $\lambda>\lambda_1$, \eqref{main_problem} has two nontrivial solutions.
        \item For $\lambda>\lambda_2$ or $\lambda=\lambda_2$ and $\frac{f(t)}{t}\geq \lambda_2$ in a neighborhood $U$ of $t=0$, \eqref{main_problem} has three nontrivial solutions.
        \item For $\lambda>\lambda_2$, assume that if $\lambda\in\sigma(-\Delta)$ either
        $$\frac{f(t)}{t}\geq \lambda\ \ \ \ \text{or}\ \ \ \ \frac{f(t)}{t}\leq \lambda$$
        holds for $t\not=0$ in a neighborhood of $t=0$, then \eqref{main_problem} has at least four nontrivial solutions.
    \end{enumerate}
\end{theorem}
As one can see in the previous theorem and will become clearer later, the relative position of the asymptotic parameter $\lambda$ and $\lambda_1$ with the spectrum of a self adjoint operator, the classical Laplacian $\sigma(-\Delta, \Omega)$, is deeply connected to the multiplicity of the problem \eqref{main_problem}. In Chapter \ref{ch-5} we will introduce significant variations or generalizations of the theorem above.

\section{Motivation}


There are several strong reasons to study the properties of stable solutions, both from a purely mathematical perspective and from the viewpoint of their application to physical theories that aim to be taken seriously as plausible descriptions of some aspect of reality. From a mathematical standpoint, as discussed in \cite{dupaigne}, stability is a property that can be exploited to derive a wealth of bounds and estimates, providing precise control over the behavior of solutions, particularly their regularity. This part will be developed in the first place, in two steps according to two famous conjectures that are still open in a sense that will be precised. Secondly, some of the physical context in which these equations and concepts have been employed will be exposed. Despite this thesis being mainly focused on pure mathematical contributions in which the usual tendency is to skip this part as quickly as possible or to avoid any further discussion, this part is intended to have a precise treatment that allows the reader to attain a good intuition of the physical phenomena involved. Furthermore, in the Section \ref{subsec:stability-in-field-theories} a linearized dynamical stability setting is employed to study the physical plausibility of different field theories in physics.\\

In the long-standing problem of determining when weak solutions of a boundary value problem are actually classical solutions, or at least limits of classical solutions, the concept of stability has been proven to be a valuable tool. 
In Proposition \ref{prop-approximation-by-stable-solutions} at least for convex nonlinearities, it was shown that weak solutions falling inside the family of stable solutions are, in a sense, always limit of classical solutions. However, the reciprocal is not true in general as we see in the next section.

\subsection{The Brezis-V\'azquez conjecture}\label{sec:brezis-vazquez-conjecture}
Let us show the analysis of the following slight modification of the problem \eqref{main_problem}.
\begin{equation}
\label{eq:lambda}
\left\{
\begin{array}{cl}
-\Delta u=\lambda f(u) & \text{in }\Omega,\\
u>0 & \text{in }\Omega,\\
u=0 & \text{on }\partial\Omega,
\end{array}
\right.
\end{equation}
where $f:[0,+\infty)\to \mathbb{R}$ satisfy $f(0)>0$ and nondecreasing, convex, and superlinear at $+\infty$ (in the sense that $\lim_{t\to+\infty}\frac{f(t)}{t}=+\infty$). \\
It is now well known that the ($L^\infty$-norm of the) branch of classical stable solutions of this $\lambda$-family of semilinear elliptic problems forms a continuum $\lbrace u_\lambda:\lambda\in(0,\lambda^\star)\rbrace$ up to some critical value $\lambda^\star>0$ and in the limit $\lambda\uparrow\lambda^\star$ the problem has a unique weak solution $u_\star$ which is called the {\em extremal solution}. In the literature, this problem is usually referred to as the ``Gelfand problem'', or a ``Gelfand-type problem''. It was first presented by Barenblatt in a volume edited by Gelfand \cite{gelfand}, and was motivated by problems occurring in combustion\footnote{Originally, Barenblatt introduced problem \eqref{eq:lambda} for the exponential nonlinearity $f(u)=e^u$ (arising as an approximation of a certain empirical law). Nowadays, the terminology of Gelfand or Gelfand-type problem applies to all $f$ satisfying the assumptions above.}. 
Later, it was studied by a series of authors; see for instance  \cite{brezis,dupaigne} for a complete account on this topic.

The basic results concerning \eqref{eq:lambda} can be summarized as follows
(see for instance \cite[Theorem 1 and Remark 1]{brezis} or the book \cite{dupaigne} by Dupaigne):

\begin{theorem}[\cite{brezis,dupaigne}]\label{thm:lambda star}
There exists a constant $\lambda^\star \in (0,+\infty)$ such that:
\begin{itemize}
\item[(i)] For every $\lambda \in (0,\lambda^\star)$ there is a unique $W^{1,2}_0(\Omega)$ stable solution $u_\lambda$ of \eqref{eq:lambda}. 
Also, $u_\lambda$ is a classical solution and $u_\lambda<u_{\lambda'}$ for $\lambda<\lambda'$.

\item[(ii)] For every $\lambda >\lambda^\star$ there is no classical solution.

\item[(iii)] For $\lambda=\lambda^\star$ there exists a unique $L^1$-weak solution $u^\star$ which is called the  extremal solution of  \eqref{eq:lambda}
and satisfies $u_\lambda\uparrow u^\star$ as $\lambda\uparrow \lambda^\star$.
\end{itemize}
\end{theorem}
On the other hand, the uniqueness of weak solution for $\lambda=\lambda^\star$ was proved by Martel \cite{martel}.

Brezis noticed a strong dependence of the behavior of the extremal solution of \eqref{eq:lambda} on the dimension. In \cite[Open problem 1]{brezis}, he asked the following:
 
\vspace{3mm}

\noindent {\it Is there something ``sacred'' about dimension 10? More precisely, is it possible in ``low'' dimensions to construct some $f$ (and some $\Omega$) for which the extremal solution $u^\star$ is unbounded? Alternatively, can one prove in ``low'' dimension that $u^\star$ is
smooth for every $f$ and every $\Omega$? }

\vspace{3mm}

\begin{remark}
    The key observation is that Brezis's question in low dimensions can be interpreted as an a priori bound for the family of stable solutions $\{u_\lambda\}_{\lambda<\lambda^\star}$. Consequently, analyzing the regularity of extremal solutions is essentially equivalent to establish a priori estimates for stable classical solutions.
\end{remark}

Framed in this way, the problem can be connected to the following long-standing conjecture proposed by Brezis and Vázquez in \cite{brezis-vazquez}.

\vspace{3mm}

\noindent  \textbf{Brezis-V\'azquez conjecture:} {\it Let $u\in W_0^{1,2} (\Omega)$ be a stable weak solution to \eqref{main_problem}. Assume that f is positive, non-decreasing, convex, and superlinear at $+\infty$, and let $N\leq 9$. Then u is bounded. }

This conjecture was unsolved for decades until a relatively recent publication \cite{cabre-figalli-rosoton-serra} finally gave a conclusive answer. This will be analyzed in detail in Section \ref{sec-regularity-literature}. \\

\subsection{The De Giorgi's conjecture}
Often, the interest of the community in solving a problem increases as time goes by without a solution. That is sometimes enough to become a true motivation for mathematicians that take it as a challenge. This is one of the features, certainly not the main one during decades, of the so-called De Giorgi conjecture, a conjecture posed on \cite{de-giorgi-original-conjecture} by Ennio de Giorgi that states the following
\vspace{3mm}

\noindent
{\em 
{\rm {\bf De Giorgi's Conjecture:}}  
Let $u:\mathbb{R}^{N}\to (-1,1)$ be a solution of the Allen-Cahn equation 
\begin{equation}    \label{eq-allen-cahn}
    -\Delta u=u -u^3\quad \mbox{in $\mathbb{R}^{N}$}.
\end{equation}
satisfying the monotonicity condition: $\frac{\partial u}{\partial x_N}>0$ in $\mathbb{R}^N$. Then, $u$ is a $1$-d solution (or equivalently, all level sets $\{u=s\}$ of $u$ are
hyperplanes), at least if $N\leq 8$.
}\\

\begin{figure}
    \centering
    \includegraphics[width=0.5\linewidth]{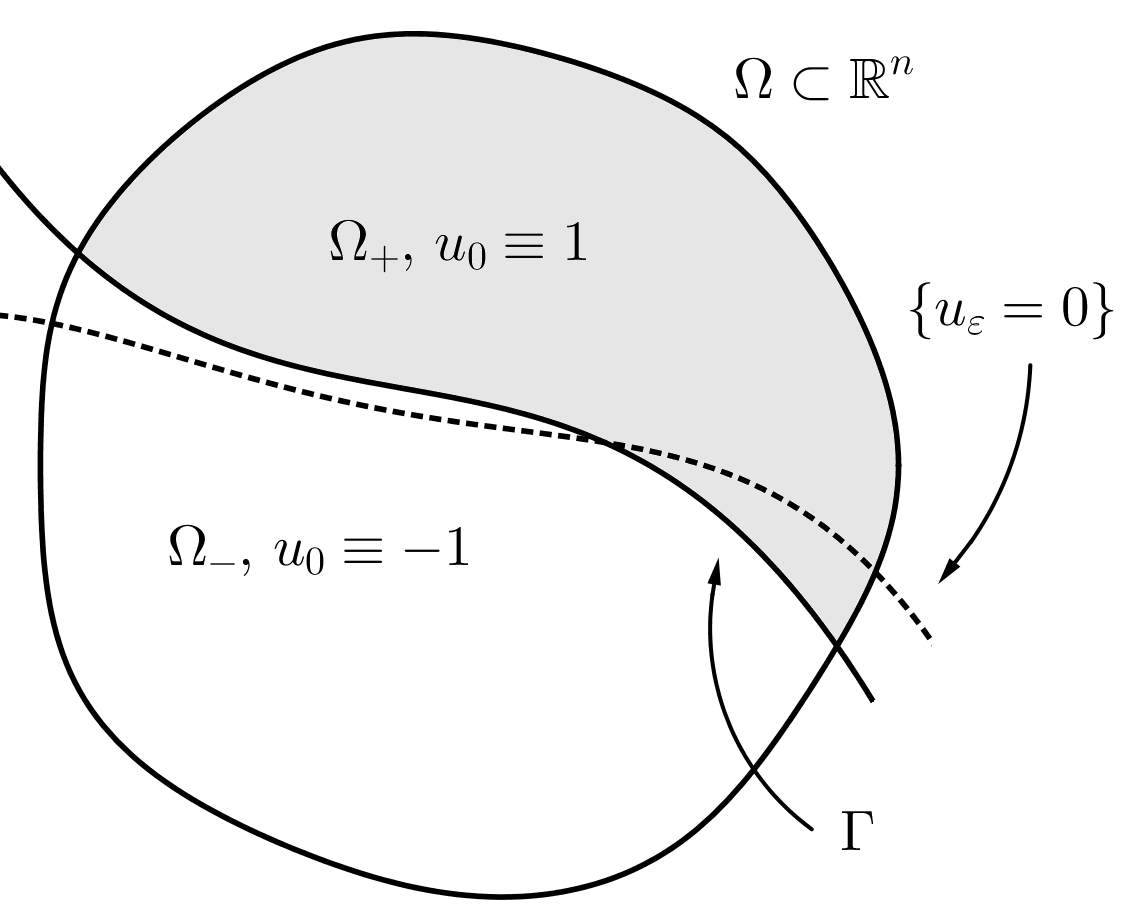}
    \caption{$\Gamma^-$-convergence to a characteristic function limiting interface between two substances described by the Cahn-Hilliard model (figure from \cite{Cabré2018})}
    \label{fig:interface-modica-mortola}
\end{figure}

There exists a deep connection that is still not fully understood between the Allen Cahn equation, this conjecture in particular, and the theory of minimal surfaces in which the concept of stability seems to play a significant role. Let us describe this assertion in a bit more detail.
Let us rescale the problem and pose it on a bounded domain, i.e., let $u_\epsilon(x):=u\left(\frac{x}{\epsilon}\right)$ for all $x\in\Omega$. The limit $\Omega=\mathbb{R}^N$ can be recovered by taking $\epsilon\to 0^+$. The energy for the equation \eqref{eq-allen-cahn} rescaled in a domain $\Omega$ reads 
\begin{equation}
    E_\epsilon(u,\Omega)=\int_{\Omega}\left\lbrace\frac{\epsilon}{2}|\nabla u|^2+\frac1 \epsilon W(u)\right\rbrace dx
\end{equation}
where $W(u)=\frac{1}{4}(1-u^2)^2$. Usually, this energy models diffusive processes between two substances that are contained in a region of space $\Omega$ (see Figure \ref{fig:interface-modica-mortola} and \cite[Sec. 3.4]{Elliott1989}). 
\begin{definition}[Minimizer of the Allen-Chan equation]\label{definit:MinimizAC}
A function  $u:\mathbb{R}^{n}\to (-1,1)$ is said to be a {\it minimizer} of \eqref{eq-allen-cahn}
when
$$E_{\epsilon=1} (u,B_R) \leq E_{\epsilon=1} (v,B_R)$$
for every open ball $B_R$ and functions $v: \overline{B_R}\to \mathbb{R}$ such that $v\equiv u$ on 
$\partial B_R$.
\end{definition}
Here, the variational approach to the theory of minimal surfaces comes surprisingly into play, as long as, it turns out that the interface of the two substances to tends to prefer to minimize its area.  The relevant functional is for minimal surfaces is in general the following.
\begin{definition}[Perimeter]
Let $\Omega \subset \mathbb{R}^N$ be an open set, regular enough. For a given open ball $B_R$ we define the {\it perimeter of} $\Omega$ {\it in} $B_R$ as 
\begin{equation*}
P(\Omega, B_R) := \mathscr{H}_{N-1} (\partial \Omega \cap B_R ),
\end{equation*}
where $\mathscr{H}_{N-1}$ denotes the $(N-1)$-dimensional Hausdorff measure.
\end{definition}

It is a known fact (see chapter 10 of \cite{giusti}) that the first and second variations of the perimeter (this time perturbations are taken via normal perturbation fields of compact support $\Omega_\epsilon=\Phi_\epsilon(\Omega)$ where $\Phi_\epsilon=I+\epsilon\,\varphi\,\nu$, $\varphi\in C^\infty_c(B_R)$ is a scalar function and $\nu$ is the normal vector to the surface $\partial \Omega$) are given by
\begin{align}
\left.\frac{d}{d\epsilon} P(\Omega_\epsilon, B_R) \right|_{\epsilon=0} 
&=\int_{\partial \Omega} H\,\varphi \,d\mathscr{H}_{N-1},
\label{eq-first-variation-perimeter}\\
\left.\frac{d^2}{d\epsilon^2}  P(\Omega_\epsilon, B_R) \right|_{\epsilon=0}
&=\int_{\partial \Omega}\left\{|\nabla_T\varphi|^2-(c^2-H^2)\varphi^2
\right\}d\mathscr{H}_{N-1},
\label{eq-second-variation-perimeter}
\end{align}
where $H=H(x)$ is the {\it mean curvature} of $\partial \Omega$ at $x$ and $c^2=c^2(x)$ is the sum of the squares 
of the $N-1$ principal curvatures $k_1, \dots, k_{N-1}$ of $\partial \Omega$ at $x$. To be more precise,
$$
{H}(x)= k_1 + \dots + k_{N-1} \quad \mbox{ and } \quad c^2= k_1^2 + \dots + k_{N-1}^2 .
$$
In \eqref{eq-second-variation-perimeter}, $\nabla_T$ is the tangential gradient to the surface $\partial \Omega$, given by
\begin{equation}\label{def:tangentialgradient}
\nabla_T\,\varphi = \nabla \varphi - (\nabla \varphi \cdot \nu)\, \nu
\end{equation}
for any function $\varphi$ defined in a neighborhood of $\partial \Omega$. Thus, a \textit{minimal surface} is a surface with zero mean curvature $H=0$ on $\partial\Omega$ and is called \textit{stable minimal surface} if \eqref{eq-second-variation-perimeter} is non negative for all balls $B_R$.
\begin{theorem}[\cite{modica-mortola}]
    Let $X$ be the space of Lebesgue-measurable functions such that $\int_\Omega u=M<|\Omega|$, equipped with the $L^1(\Omega)$ norm. Define the functional
    \begin{equation}
        \hat E_\epsilon(u)=
            \left\{
                \begin{array}{cl}
                E_\epsilon(u) & \text{if  }u\in W^{1,2}(\Omega)\cap X\\
                +\infty & \text{else}
                \end{array}
            \right.
    \end{equation}
    Then, in the variational sense of $\Gamma^-(X)$-convergence (see the definition in \cite{modica-mortola}) 
    \begin{equation}
        \hat E_\epsilon(u_\epsilon,\Omega)\xrightarrow[\epsilon\to 0]{\Gamma^-(X)}c_W\,P(\partial\lbrace u^{\Gamma^-}_0=1\rbrace,\Omega),
    \end{equation}
    where $c_W=\int_0^1\sqrt{2W(s)}\,ds$ and the limiting function $u^{\Gamma^-}_0$ is nothing but the an indicatrix function inside $\Omega$ that separates two subregions corresponding to the two substances by a minimal hypersurface, the interface (see Figure \ref{fig:interface-modica-mortola}).
\end{theorem}
That said, it is intuitive that stable solutions to \eqref{eq-allen-cahn} can be linked to stable minimal surfaces in some way, maybe even before taking the asymptotic limit $\epsilon\to 0^+$. There are two results that let emphasize this assertion
\begin{theorem}[\cite{simons}]\label{th-simons-cone}
Let $\Omega \subset \mathbb{R}^N$ be an open set such that $\partial \Omega$ is a stable minimal cone and $\partial \Omega\setminus\{0\}$ is regular. 
Thus, we are assuming $H=0$ and
\begin{equation*}
\int_{\partial \Omega}\left\{|\nabla_T\,\varphi|^2-c^2\varphi^2\right\}d\mathscr{H}_{N-1}\geq 0,
\label{eq:1-2v2}
\end{equation*}
for every $\varphi\in C^1(\partial \Omega)$ with 
compact support outside the origin.
If $3 \leq N\leq 7$, then $\partial \Omega$ is a hyperplane.
\end{theorem}
\begin{theorem}[\cite{savin}]\label{th-savin}
Assume that $N\leq 7$ and that $u$ is a minimizer of \eqref{eq-allen-cahn} in~$\mathbb{R}^N$. 
Then, $u$ is a $1$-d solution.
\end{theorem}
The interesting point is that if the solution satisfies the additional assumption (not appearing in the original conjecture of De Giorgi)
\begin{equation}\label{eq-de-Giorgi-limits-assumption}
\lim_{x_N\rightarrow\pm\infty}u(x', x_N)=\pm 1
\quad \text{ for all }x'\in \mathbb{R}^{N-1}.
 \end{equation}
then (together with the gain of one more dimension by the monotonicity of the solution) the Theorem \ref{th-savin} leads to a positive answer to the De Giorgi conjecture (see \cite[Theorem 2.4]{Cabré2018}). The proof of Theorem \ref{th-savin} for its part, uses Theorem \ref{th-simons-cone}. In this way we see that stability of minimal cones seems to be a key ingredient.

\subsubsection{The saddle-shaped solution vanishing on the Simons cone}\label{subsubsec:saddlesha}

Let $m \ge 1$, $x=(x_1,\dots, x_{2m})\in\mathbb{R}^{2m}$, and $s$, $t$ denote the two
radial variables
\begin{equation}\label{coor}
s =   \sqrt{x_1^2+...+x_m^2} \quad \mbox{ and } \quad t  =  \sqrt{x_{m+1}^2+...+x_{2m}^2}.
\end{equation}
The so-called \textit{Simons cone} is given by 
\begin{equation}\label{def-simons-cone}
    \mathcal{SC}=\{s=t\}=\partial \Omega, \quad\text{ where }
\Omega =\{s>t\}.
\end{equation}
\begin{definition}[Saddle-shaped solution]\label{def:saddlesolution} 
We  say  that $u:\mathbb{R}^{2m}\rightarrow\mathbb{R}$ is a {\it saddle-shaped
solution}  (or simply a saddle solution) of the Allen-Cahn equation 
\begin{equation}\label{eq2m}
-\Delta u= u- u^3 \quad {\rm in }\;\mathbb{R}^{2m},
\end{equation}
whenever $u$ is a
solution of
\eqref{eq2m} and, with $s$ and $t$ defined by \eqref{coor},
\renewcommand{\labelenumi}{$($\alph{enumi}$)$}
\begin{enumerate}
\item $u$ depends only on the variables $s$ and $t$. We write
$u=u(s,t)$; 
\item $u>0$ in $E:=\{s>t\}$; 
\item  $u(s,t)=-u(t,s)$ in $\mathbb{R}^{2m}$.
\end{enumerate}
\end{definition}

\begin{remark}
Notice that if $u$ is a saddle-shaped solution, then we have $u=0$ on $\mathcal{SC}$.
\end{remark}
While the existence of a saddle-shaped solution is easily established, its uniqueness is more delicate.
\begin{theorem}[\cite{cabre}] \label{unique}
For every even dimension  $2m\geq 2$,
there exists a unique saddle-shaped solution $u$ of \eqref{eq2m}.
\end{theorem}
Due to the minimality of the Simons cone when $2m \ge 8$ (and also because of the minimizer from \cite{LWW} referred to after Theorem \ref{th-savin}), the saddle-shaped solution is expected to be a minimizer when $2m \ge 8$:
\vspace*{3mm}

\noindent\textbf{Open problem: }Is the saddle-shaped solution a minimizer of \eqref{eq2m} in $\mathbb{R}^8$, or at least in higher even dimensions?.
\vspace*{3mm}

No progress has been made on this open problem except for the following result. It establishes stability (in this case in the sense of \eqref{def-stability} that is of course weaker than minimality) for $2m \ge 14$.
\begin{theorem}[\cite{cabre}]\label{stab} 
If $2m\geq 14$, the saddle-shaped solution $u$ of \eqref{eq2m} is stable in $\mathbb{R}^{2m}$.
\end{theorem}
The study and characterization of stable solutions to the Allen Cahn equation and stable minimal surfaces in $\mathbb{R}^N$ for low dimension $N$, could be a major step forward in finally solving De Giorgi's conjecture in its general form.

\subsection{Physical motivation: Stability of Field Theories in physics}\label{subsec:stability-in-field-theories}

In the modern approach of theoretical physics a \textit{field theory} can be understood as a triple $(\mathcal{M}_g,V_G,\lbrace\Phi_i\rbrace_{i\in I})$ where $\mathcal{M}$ is a (commonly $4$-dimensional) differentiable manifold  equipped with a Lorentzian-signature metric tensor $g$ on its tangent space $T\mathcal{M}$, called \textit{spacetime}, $V$ is a (possibly infinite dimensional) \textit{representation} vector space of a Lie group $G$ that realizes the symmetries of the field theory and $\lbrace\Phi_i:\mathcal{M}\to V\rbrace_{i\in I}$ is a set of differentiable \textit{fields} that represent the content in matter and forces of the theory.

In general, fields over $\mathcal{M}$ have precise tensorial transformation properties according to linear representations of a sort of semi-direct product $\operatorname{Diff}\mathcal{M}\ltimes G$ of the group of diffeomorphisms of the spacetime, that represent the symmetries under general change of coordinates or reference frames; with a semisimple Lie group (through its Lie algebra) that encode the internal symmetries of the theory. By the so-called \textit{gauge principle} the action of the group locally in $\mathcal{M}$ generates the interactions or forces (see for example, the classical paper \cite{gauge-principle} and Appendix \ref{app:noether} to read more about the celebrated gauge principle). 

Given a basis $\lbrace e_A\rbrace_{A\in I}$ of $V$, the components of the fields will be labeled by a multi-index $A$: $\Phi_i=\sum_{A}\Phi^A_i e_{A}$. $A$ may contain \textit{spacetime indices}  that determine the covariant or tensorial behavior of the field under the change of reference frame, as well as \textit{internal indices},  which specifies the internal symmetry transformations. We shall use the usual conventions: 
\begin{itemize}
\item Spacetime indices labeled by Greek letters from the middle of the alphabet $\mu,\nu,\rho,\sigma\dots$ and internal indices labeled with lowercase Latin letters of the from the beginning of the alphabet $a,b,\dots$
\item  Signature $(-1,+1,\dots,+1)$ on the Lorentzian metric, lower and upper indices and 
$$(x_0\equiv t,x_1,\dots,x_{\dim \mathcal{M}-1})^{T}\overset{\cdot}{=}x^\mu \ \ \mu\in \lbrace0,1,\dots,\dim \mathcal{M}\rbrace.$$
\item Einstein summation convention.

\item Torsion free or Levi-Civita connection on $\mathcal{M}$ with covariant derivative denoted by $D_\mu=\partial_\mu +\Gamma_\mu$ where $\Gamma_\mu$  stands for the operator whose components $\Gamma_{\mu\nu}^\rho$ are the so-called Christoffel symbols. 
\end{itemize}

The fields dynamics is governed by the \textit{action functional}
\begin{equation}\label{eq:action-def}
    \mathcal{S}(\Phi):=\int_\mathcal{M} d^4x \sqrt{|g|}\mathcal{L}(\Phi),
\end{equation}
i.e., it contains all the information of the interacting fields over the spacetime. See that the action is nothing but the integral of the \textit{Lagrangian density} which for zero-curvature (flat) spacetimes already encodes all the information about the interaction of the fields. This optimization/variational framework has proven to be a very successful description of almost all known physical phenomenology at the most fundamental level.

Specifically in particle physics, the spacetime under consideration is always the Minkowski spacetime $\mathbb{R}^{1,3}$, the four-dimensional flat spacetime since gravity effects can be neglected due to the weakness of the gravity force compared to the other three fundamental forces at the available energy scales in experiments. In that case, the group of symmetry is reduced to $ISO(1,3)\times G$ where $ISO(1,3)$ stands for the Poincaré group of translations, spatial rotations and boosts in Minkowski spacetime. Particle physicists usually identify \textit{theories} and  \textit{Lagrangians} forgetting about the metric part $\sqrt{|g|}$ of the invariant volume form. In such a field theory one can clearly differentiate the two types of derivatives: $\partial_t$ the time derivative and $\partial_{x_i}$ the spatial derivatives. The time/space splitting allows us to conceive time evolution of physical quantities or observables through a differential equation of the form $\mathcal{D}(\Phi)=0$ for a certain differential operator $\mathcal{D}$. Even in a curved spacetime, under \textit{global hyperbolicity}\footnote{A spacetime \( (\mathcal{M}, g) \) is said to be \emph{globally hyperbolic} if there exists a Cauchy hypersurface \( \Sigma \subset \mathcal{M} \), i.e., a spacelike hypersurface intersected exactly once by every inextensible causal curve.} one can construct a globally well behaved evolution in time and derive those equations. They are called \textit{equations of motion} and are, of course, characterized by the theory. More precisely, they are the Euler-Lagrange equations of the action functional, so they describe its critical points, a method that has been called \textit{Principle of least action} and is where functional analysis and partial differential equations become important for physics.

In a general field theory over a spacetime, the Lagrangian may contain a bunch of different types of scalar terms constructed from contraction of all the indices of different fields. Some of them appear most of the time and have a clear interpretation. 
We shall distinguish
\begin{enumerate}
    \item \textit{Kinetic terms}: Quadratic terms in the first derivatives of the fields components. They represent the kinematic energy of a field,
    \begin{equation}\label{eq:kinetic-terms}
        \frac{1}{2}D_\mu\Phi^aD^\mu\Phi_a \subset \mathcal{L},
    \end{equation}
    in the case of a scalar field $u$ this is equal to $-(\partial_t u)^2+|\nabla u|^2$ (sometimes is preferred to call kinetic term to only the time derivative part).
    \item \textit{Mass terms}: Quadratic terms in the fields components. They give the inertia of the fields,
    \begin{equation}\label{eq:mass-terms}
        m_\Phi \Phi_A\Phi^A\equiv m_\Phi \Phi^2\subset\mathcal{L},\ \ \ \ m_\Phi>0,
    \end{equation}
    where $m_\Phi>0$ is the mass of the field. Sometimes $m_\Phi=m_\Phi(c_i,\Psi)$ depend on some scalar arrangement of other fields or parameters but still is interpreted as an effective mass (see Figure \ref{fig:higgs-mechanism})
    \item \textit{Interaction terms}: These are contractions of different fields or its derivatives of any order. The \textit{minimal coupling principle} states that the interaction between two fields $\Phi_1=\Phi$ and $\Phi_2=\Psi$ should be represented by its simplest quadratic contraction:
    \begin{equation}\label{eq:interaction-terms}
        c_{\Phi\Psi}\Phi^A\Psi_A\subset\mathcal{L},\ \ \ \ c_{\Phi\Psi}\in \mathbb{R}
    \end{equation}
    where $c_{\Phi\Psi}$ is the coupling constant between the fields $\Phi$ and $\Psi$.
    \item \textit{Curvature terms}: Given a field $\Phi^A$, they are combinations of the form
    \begin{equation}
        F_{\mu\nu}^a:=\partial_\mu \Phi^a_\nu-\partial_\nu\Phi^a\mu+f^a_{bc}\Phi^b_\mu,\Phi^c_\nu
    \end{equation}
    where $f^a_{bc}$ stand for the structure constants of the underlying Lie algebra. 
    \item \textit{Potential terms}: Higher order contractions of the fields.
    \begin{equation}
        c\Phi^2+\lambda \Phi^4\subset\mathcal{L},\ \ \ \ (c<0)
    \end{equation}
\end{enumerate}
To measure the distribution and flow of the energy and momentum in the spacetime one divide the action in two terms
\begin{align}
    \mathcal{S}(g,\Phi)=\mathcal{S}_{\text{gravity}}(g)+
    \mathcal{S}_\text{matter}(g,\Phi),
\end{align}
perturb the metric $g_{\mu\nu}+\epsilon h_{\mu\nu}$ and take the first variation of the matter part of the action. It is a task for the so-called \textit{stress-energy tensor} $T_{\mu\nu}$ which is determined by,
\begin{equation}\label{def:stress-energy-tensor}
    d\mathcal{S}_{\text{matter}}(h):=\left.\frac{d}{d\epsilon}\mathcal{S}_{\text{matter}}(g+\epsilon h)\right|_{\epsilon=0}=\int_\mathcal{M}\sqrt{-g}\,T^{\mu\nu}h_{\mu\nu}\,d^4x
\end{equation}
In Minkowski spacetime $x\in\mathbb{R}^{1,3}$ the above definition connects\footnote{Even though this canonical version coming from Noether theorem, $T_{}^\mu{}_\nu $, could be in principle non symmetric, in the celebrated paper \cite{belinfante} the author describes a method to symmetrize it adding total derivatives (that are terms that maintain the action invariant).} with the canonical definition of energy-momentum tensor that arise as the conserved current associate by the \textit{Noether theorem} (see Theorem \ref{th:noether-theorem} in Appendix \ref{app:noether}) to the symmetry under spacetime constant translations (i.e., $x_\mu\mapsto x_\mu+v_\mu$ with $v\in\mathbb{R}^{1,3}$):
\begin{equation}\label{def:canonical-stress-energy-tensor}
    (T_{}^\mu{}_\nu)_{can} := \frac{\partial \mathcal{L}}{\partial \partial_\mu \Phi^A}\partial_\nu\Phi^A - \mathcal{L} \delta^\mu{}_{\nu}\,,
\end{equation}
To get physical intuition on it, one can take the $2$-dimensional submanifold given by the local coordinates $x^\mu$ and $x^\nu$ with values in a neighborhood of $x$, then $T_{\mu\nu},\ \mu,\nu\in\lbrace1,2,3\rbrace$ measures the stress (force over area) over the spatial submanifold. On the other hand, $T_{0\mu},\ \mu\in\lbrace1,2,3\rbrace$ measures the flux of energy and momentum in direction $x^\mu$. Finally, and most importantly, $T_{00}$ measures the flux of the energy in the temporal direction, i.e., the energy density.  This energy density is going to play as a snitch for instabilities as we will see below. 
There are currently two theories that are both independently the two fundamental pillars of physics. Both are field theories that can be described by its action functional.
\begin{enumerate}
    \item \textit{General Relativity}: Published in 1915, due almost solely to A. Einstein,
    \begin{equation}
        \mathcal{S}_{GR}=\frac{1}{16\pi G}\int_\mathcal{M}d^4x \sqrt{|g|}\,(R-2\Lambda)
    \end{equation}
    where $R$ stands for the Ricci scalar curvature of $\mathcal{M}$. The Euler-Lagrange equations of this action are the \textit{Einstein field equations}.
    \begin{equation}\label{eq:general-relativity}
        R_{\mu\nu}-\frac{1}{2}g_{\mu\nu}R+\Lambda=-\frac{8\pi G}{c^4} T_{\mu\nu}
    \end{equation}
    where $\Lambda$ is the so-called \textit{cosmological constant} that measures a constant  energy density of spacetime\footnote{It is usually employed to account for the so-called \textit{Dark Energy} that makes the universe exhibit an accelerating expansion compatible with the observations.}
    \item \textit{Quantum mechanics}: Due to many different physicists of the early of the 20th century like M. Planck, A. Einstein, N. Bohr, E. Schrödinger, B. Heisenberg, L. de Broglie, P. Dirac, among them. In the non relativistic case:
    \begin{equation}
        \mathcal{S_{QM}}=\int_{(0,T)\times\Omega} dtd^3x \frac{i\hbar}{2} \left( \Psi^* \frac{\partial \Psi}{\partial t} - \Psi \frac{\partial \Psi^*}{\partial t} \right)
- \frac{\hbar^2}{2m} \nabla \Psi^* \cdot \nabla \Psi
- V(\vec{x}, t) \Psi^* \Psi
    \end{equation}
    where $\Psi:\Omega\to \mathbb{C}$ is the so-called \textit{wave function} whose modulus squared $|\Psi(x)|^2$ describes the probability amplitude of finding a quantum particle at the position $x\in \Omega$. The famous \textit{Schrödinger equation} can be derived by the action:
    \begin{equation}
    i\hbar \frac{\partial}{\partial t} \Psi = \left( -\frac{\hbar^2}{2m} \nabla^2 + V(\vec{x}, t) \right) \Psi
    \end{equation}
\end{enumerate}
Additionally, it is worth saying that the actual field theory that is currently the pillar of particle physics is the result of a great effort by many physicists in the first half and early second half of the 20th century to make a theory in which quantum mechanics and \textit{special relativity} (that doesn't take into account gravity) are compatible and in which little by little had to describe a greater number of new particles and phenomena that were being discovered. That theory is now called \textit{Standard Model} (SM) and is the most accurate description of nature achieved in the history. The action of SM can be described in the same way as before
\begin{equation}
    \mathcal{S}_{SM}=\int _{M}d^4x\, \mathcal{L}_{SM}
\end{equation}
where
\begin{equation}\label{eq:standard-model-lagrangian}
\begin{split}
    \mathcal{L}_{\text{SM}} = & -\frac{1}{4} G_{\mu\nu}^a G^{\mu\nu}_a - \frac{1}{4} W_{\mu\nu}^b W^{\mu\nu}_b - \frac{1}{4} B_{\mu\nu} B^{\mu\nu} + \sum_{\Psi} \bar{\Psi} i \gamma^\mu D_\mu \Psi \\
&+ (D^\mu \Phi)^\dagger (D_\mu \Phi) - \mu^2 \Phi^\dagger \Phi + \lambda (\Phi^\dagger \Phi)^2 - \sum_{f} y_f \bar{\Psi}_L \Phi \Psi_R + \text{h.c.}
\end{split}
\end{equation}
With a little patience one can identify the different types of terms that have been mentioned above.\\
All these theories give us predictions of the time evolution of physical observables through the equations of motion once a set of initial conditions is introduced. Therefore, a notion of dynamical stability with respect to the initial conditions can be applied.\\
Given any time $t_0>0$ and a field $\Phi(t,\vec x)$ one take a perturbation of the field $\Phi_{\epsilon}(x,t)=\Phi(t,\vec x)+\epsilon \varphi(t,\vec x)$. The space in which the perturbations live depends strongly on the equation of motion under consideration and often the search for the appropriate space of perturbations is part of the problem. Then, one performs the second variation of the action.
\begin{equation*}
    d^2\mathcal{S}^\Phi(\varphi,\varphi)=\left.\frac{d^2 }{d\epsilon^2}\mathcal{S}(\Phi+\epsilon\varphi)\right|_{\epsilon=0}
\end{equation*}
The result will give us the dynamics of the evolution of the perturbations that will determine the stability of the theory.
\begin{remark}
    Note that the equations for linear perturbations can also be obtained from the perturbative expansion of the Lagrangian
\begin{equation}
    {\mathcal{L}}(\Phi+\epsilon\varphi)={\mathcal{L}}(\Phi)+\sum_{n=1}^{\infty}\epsilon^n{\mathcal{L}}_{\Phi}^{(n)}(\varphi),
\end{equation}
where ${\mathcal{L}}^{(n)}_{\Phi}(\varphi)$ is of order $\epsilon^n$ and depends on the solution $\Phi$ around which the perturbative expansion is considered. If the order under consideration is $\epsilon$ then it is called lineal or linearized stability.
\end{remark}
To make a correspondence with the convexity of the energy functional as in \eqref{def-stability} one has to perform a \textit{Legendre transformation} to go from the Lagrangian density to the so-called \textit{Hamiltonian} density which coincides with the energy functional for field theories at the classical (non quantum) level and with the usual energy functional \eqref{energy-functional} that we have been using for semilinear elliptic equations when stationary states are considered. Then, the dynamics can be expressed in terms of the fields and their canonical momenta (without velocities $\partial_t\Phi$). More precisely
\begin{equation}\label{eq:def-hamiltonian}
    \mathcal{H}(\Phi,\Pi)=\Pi_a\,\frac{\partial\Phi^a}{\partial t}\,-\mathcal{L}
\end{equation}
where the \textit{canonical momenta} are defined in components as
\begin{equation}\label{def:canonical-momenta}
    \Pi_a=\frac{\partial \mathcal{L}}{\partial (\partial_t\Phi^a)}
\end{equation}

In order for that transformation to be well defined, it is essential to have the global hyperbolicity of the spacetime, as well as, the capacity to invert the relation \eqref{def:canonical-momenta}. By the inverse function theorem this implies that
\begin{equation}
    \mathcal{K}_{ab}=\frac{\partial^2 \mathcal{L}}{\partial (\partial_t\Phi^a)\partial (\partial_t\Phi^b)}\not=0 \ \ \ \ \ \text{in }\mathcal{M}
\end{equation}
The matrix $(\mathcal{K}_{ab})$ is called the \textit{kinetic matrix} and governs the ability of the fields configurations to properly propagate and evolve in time\footnote{In fact, in quantum field theory, one defines the \textit{kinetic operator} $\mathcal{K}$ as the linearized operator $\mathcal{K}=d\mathcal{D}$ of the equation of motion $\mathcal{D}(\Phi)=0$. It would appear as $\Phi^\dagger\mathcal{K}\Phi$ after integrating by parts the kinetic energy term $D_\mu\Phi^\dagger D^\mu\Phi$. The kinetic matrix woul be just the temporal part of the matrix representation of the kinetic operator in spacetime coordinates. It allows to define the so-called \textit{propagator} of the field $\mathcal{D}_F(x,y)$, that is which is nothing but the Green's function/operator of the kinetic operator, i.e., the integral kernel corresponding to the inverse operator $\Phi(x)=\int _{\mathcal{M}}\mathcal{D}_F(x,y)\mathcal{K}( \Phi)(y)\,d^4y$.}.
\begin{remark}
    If $(\mathcal{M},g)$ is the Minkowski spacetime $\mathbb{R}^{1,3}$, with canonical kinetic terms, i.e., $\mathcal{K}=\mathbb{I}$, then the energy density $T_{00}$ defined in \eqref{def:stress-energy-tensor} and the Hamiltonian \eqref{eq:def-hamiltonian} coincides: $$T_{00}=\mathcal{H}.$$
\end{remark}

\subsubsection{Types of Instabilities in Field Theories}

We now discuss the main types of instabilities encountered in theoretical physics. To analyze the fundamental types of instabilities that may arise, consider a scalar Lagrangian quadratic (see Theorem \ref{ostrogradsky} in Appendix \ref{app:ostrogradsky}) on the linear scalar field governing perturbations, on an isotropic background solution ($\partial_{x_i}\Phi=0\ \ \forall i \in\lbrace1,2,3\rbrace$) leading to a scalar perturbation Lagrangian of the form:
\begin{equation}\label{eq:EffectiveScalarLag}
\mathcal{L}^{(1)}=\frac{1}{2} (\partial_t\varphi)^2- \frac{b}{2a} (\vec{\nabla}\varphi)^2-\frac{\mu}{2a} \varphi^2,
\end{equation}
where the coefficients $a$, $b$, $\mu$ will determine the relevant information about the stability of the background solution. The equations of motion are 
\begin{equation}
\partial_{tt}\varphi-\frac{b}{a}\Delta\varphi+\frac{\mu}{a}\varphi=0.
\end{equation}
Let us  we express the solution in Fourier expansion 
\begin{align*}
    \varphi (x)=\int\varphi_k(x)  \,d^3k
\end{align*}
 where
\begin{equation}\label{eq:SolutionKleinGordonInstabilities}
\varphi_k=A_k \,e^{-i\left(\sqrt{\frac{b}{a}}\vec{k}\cdot\vec{x}-\omega t\right)}+A^\ast_k e^{i\left(\sqrt{\frac{b}{a}}\vec{k}\cdot\vec{x}-\omega t\right)}
\end{equation}
with frequency $\omega=+\sqrt{\frac{b}{a}|\vec{k}|^2+\frac{\mu}{a}}$. We can have the following types of instabilities:

\begin{enumerate}

\item{\textbf{Gradient Instability ($a<0,b>0$ or $a>0,b<0$)}}

We would say that the spatial gradient term has the wrong sign:
\[
\mathcal{L} \supset -\frac{1}{2} c_\varphi^2 (\nabla \varphi)^2 \quad \text{with } c_\varphi^2 < 0
\]
where $c_\varphi=\sqrt{b/a}$. Clearly $c_\varphi$ cannot be interpret as the velocity of propagation in this case. The frequency $\omega$ is purely imaginary, triggering an exponential growth of modes with an arbitrary fast growth rate as $\vec k$ increases. They have been found in some cosmological scalar-field models (see, for example \cite{gradient-instability-reference}).

\item{\textbf{Ghosts ($a<0$ and $b<0$)}}

The kinetic term has the wrong sign:
\[
\mathcal{L} \supset -\frac{1}{2} (\partial_\mu \varphi)^2
\]
The Hamiltonian becomes unbounded from below. If the theory is interacting this kind of field could trigger other fields to achieve arbitrary large energy modes by interacting with the ghost field (that turn into arbitrary large negative energy modes). That is why is also called a \textit{run-away} instability and in quantum mechanics leads to unitarity violation. They could also be found in some higher-derivative gravity or modified gravity theories \cite{ghost-condensate-reference}.

\item{\textbf{Tachyonic Instability ($a,b>0$, $\mu<0$)}}

A field has a wrong sign mass term in the Lagrangian:
\[
\mathcal{L} \supset -\frac{1}{2} \mu\,\varphi^2 \quad \text{with } \mu > 0
\]
Small fluctuations grow exponentially due to non zero imaginary part of $\omega$
 \begin{equation}\label{eq:Tachyonicmodes}
\varphi_k=A_k e^{-i\sqrt{\frac{b}{a}}\vec{k}\cdot\vec{x}}e^{-|\omega| t}+A^\ast_k e^{i\sqrt{\frac{b}{a}}\vec{k}\cdot\vec{x}}e^{|\omega| t} \end{equation}

A prominent instance of this phenomenon is the well-known Higgs mechanism within the Standard Model (see \eqref{eq:standard-model-lagrangian} above)  the term $-\mu \Phi^\dagger\Phi$ acts as a tachyonic mass term that disappears after the \textit{spontaneous symmetry breaking} (see Figure \ref{fig:higgs-mechanism}).
\begin{center}
    \begin{figure}
        \centering
        \includegraphics[width=0.5\linewidth]{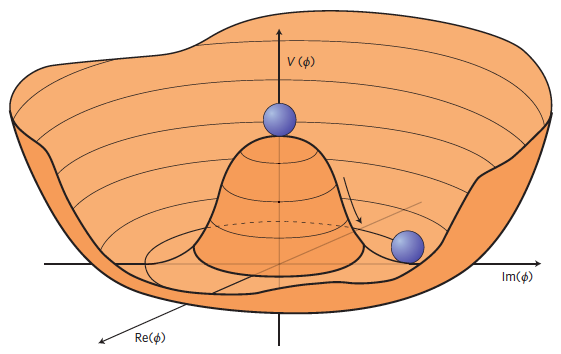}
            \caption{Higgs boson potential and face transition after the cooling of the universe universe (obtained from https://cds.cern.ch/record/2012465/files/higgspotential.png.). }
        \label{fig:higgs-mechanism}
    \end{figure}
\end{center}
\end{enumerate}

\begin{remark}
Far from being a problem, for some specific processes instabilities are essential to maintain the dynamics of the theory. As we have seen, tachyonic instability triggers spontaneous symmetry breaking and vacuum decay is key to cosmological phase transitions. Even ghost condensation appears in speculative modified gravity theories (see \cite{ghost-condensate-reference}).
\end{remark}

In physics literature one can also find other concepts of stability that are more phenomenological but not so related with the partial differential equations theory. 

\begin{itemize}

\item{Vacuum Instability.} Field theories (both classical and quantum), instabilities are sometimes a signal that the expected (for instance, suggested by the observed phenomenology in the experiment) solution of minimal energy, also called \textit{vacuum} or \textit{ground state} is not the true absolute minimum but possibly just a local one.
\begin{remark}
    This is a key difference with the elliptic or stationary stability theory. A local minimum of the energy is by definition stable in elliptic theory. However, a perturbed solution that evolves in time can carry exponentially growing perturbations that can throw the solution away of the local minimum.
\end{remark}
By the process of \textit{quantum tunneling}, a particle can overcome a potential barrier that is higher than its kinetic energy (due to the spread beyond the barrier of its wave function obeying Schrödinger equation). In this way, it can transition to a deeper, even more stable, vacuum state (see for an example, Coleman-de Luccia tunneling \cite{Coleman-tunneling}).

\item{Jeans Instability (Gravitational Collapse).} Describes the phenomenon in which matter distribution under gravity becomes unstable. For instance, when gravity overcomes internal pressure in a star, it can lead to a unstoppable collapse that ends up in the formation of a black hole (see the celebrated \cite{oppenhaimer} for more details).

\item Other dynamical instabilities occurs when the system is stable under small perturbations but unstable under large deviations like oscillons, soliton decay, bubble nucleation, etc.
\end{itemize}

As a resume of all possible uses of the concept of stability one can check the Table \ref{tb:resume-types-of-stability}.

\begin{sidewaystable}
\centering
\begin{tabular}{p{0.15\textwidth}p{0.25\textwidth}p{0.25\textwidth}p{0.32\textwidth}}
\toprule
\textbf{Feature} & \textbf{Elliptic PDEs (Static)} & \textbf{Parabolic PDEs / Dynamical Systems} & \textbf{Field Theories in Physics} \\
\midrule
\textbf{Equation type} & $-\Delta u = f(u)$ & $\partial_t u = \Delta u + f(u)$ & Euler-Lagrange from action $S[\Phi] = \int \mathcal{L}(\Phi) \, d^4x$: 
\begin{center}
    $\partial_{tt} \varphi-\Delta\Phi-m\Phi=0$
\end{center}\\
\midrule
\textbf{Object of study} & Steady solutions $u(x)$ & Evolutionary solutions $u(t,x)$ & Field configurations $\Phi(t,x)$ and small fluctuations \\
\midrule
\textbf{Stability notion} & Energy second variation $\delta^2 E \geq 0$ (linearized stability) & Lyapunov or asymptotic stability of trajectories & Small size of perturbations: Positivity of kinetic and gradient terms\\
\midrule
\textbf{Linearized problem} & Determine the spectrum of the linearized operator:\begin{center} $L_u(v):=(-\Delta - f'(u))(v)=\lambda v$ \end{center} & Linearization around steady state  $u_0$ and asymptotic behavior of the eigenvalues: \begin{center} $\partial_t v + L_{u_0}v=0;\,L_{u_0}v=\lambda v$ \end{center} &  Linear terms ($o(\epsilon)$) for perturbations of the Lagrangian:
$${\mathcal{L}}(\Phi+\epsilon\varphi)={\mathcal{L}}(\Phi)+\epsilon{\mathcal{L}}_{\Phi}^{(1)}(\varphi)$$ \\
\midrule
\textbf{Instability signals} & Negative eigenvalues of $L$ & Positive sign of the real part of eigenvalue of $-L_{u_0}$ $\Rightarrow$ exponential growth of perturbations (hence instability) & Ghosts (wrong sign kinetic term), gradient or tachyonic instabilities \\
\midrule
\textbf{Linear vs Nonlinear} & Linear elliptic stability $\wedge$ $\ \lambda_1(-\Delta-f'(u))\not=0$ $\Rightarrow$ Nonlinear asymptotic stability  & Asymptotic stability $\Rightarrow$ Lyapunov stability & Stability at first (linear) order $(\epsilon{\mathcal{L}}_{\Phi}^{(1)}(\varphi))$ imply stability of the field theory only if the higher orders are negligible (validity of the perturbative expansion).\\
\bottomrule
\end{tabular}
\caption{Comparison of stability notions in elliptic PDEs, parabolic PDEs/dynamical systems, and field theories in physics.}\label{tb:resume-types-of-stability}
\end{sidewaystable}


\begin{savequote}[0.75\linewidth]
\textit{`The regularity theory in many cases is very subtle and involves a delicate
machinery''.} Michael Struwe (\cite[Chapter I]{struwe})
\end{savequote}
\chapter{Literature Review}\label{ch-lit-rev}
In this chapter, it is presented the state of the art using stability and Morse index boundedness as a tool for achieving regularity of solutions and the study of multiplicity of semilinear elliptic equations involving degenerate terms on the nonlinearity; to give the appropriate context for a better contextualization of our work in Chapters \ref{ch-counterexample}, \ref{ch-4} and \ref{ch-5}.

\section{Regularity of stable and finite Morse index solutions}\label{sec-regularity-literature}

Extremal solutions of \eqref{eq:lambda} are a priori, merely in $L^1(\Omega)$. It is then natural to ask whether extremal solutions do belong to the natural energy space $W^{1,2}_0(\Omega)$.  This question was posed by Brezis and V\'azquez in \cite{brezis-vazquez}:

\vspace{3mm}

\noindent {\it Does there exist some $f$ and $\Omega$ for which the extremal solution $u^\star$ is a weak solution\footnote{In the sense of Theorem \ref{thm:lambda star}(iii).} not in $W^{1,2}_0(\Omega)$?}

\vspace{3mm}

Concerning this problem, it has been proved that $u^\star$ belong to the energy space $W^{1,2}_0(\Omega)$ when $N\leq 5$ by Nedev \cite{nedev}, for every $N$ when $\Omega$ is convex also by Nedev in an unpublished preprint, and when $N=6$ by Villegas \cite{villegas-extremal}. The story ends up with the collaboration of four authors in a recent celebrated paper that finally establish the regularity of almost all stable solutions to autonomous semilinear elliptic equations. It is worth mentioning here, as a result, that our contributions are made for the regularity of finite (radial)-Morse index solutions of the autonomous case, and for stable solutions in a non autonomous setting.

\subsection{Smoothness of stable solutions up to dimension 9}

In 2020 in \cite{cabre-figalli-rosoton-serra}, the authors proved the following theorem.

\begin{theorem}[Theorem 1.2, \cite{cabre-figalli-rosoton-serra}]
\label{thm:L1-CalphaW12}
Let $B_1$ denote the unit ball of $\mathbb{R}^N$. 
Assume that $u\in C^2(B_1)$ is a stable 
solution of
$$
-\Delta u=f(u) \quad \text{in }B_1,
$$
with  $f:\mathbb{R} \to\mathbb{R}$ locally Lipschitz and nonnegative.

Then
\begin{equation}
\label{eq:W12g L1 int}
\|\nabla u\|_{L^{2+\gamma}({B}_{1/2})} \le C\|u\|_{L^1(B_1)},
\end{equation}
where $\gamma>0$ and $C$  are dimensional constants.
In addition, if $N \leq 9$ then
\begin{equation}
\label{eq:Ca L1 int}
\|u\| _{C^\alpha(\overline{B}_{1/2})}\leq C\|u\|_{L^1(B_1)},
\end{equation}
where $\alpha>0$ and $C$  are dimensional constants.
\end{theorem}
The idea is to extract meaningful information about $u$ from stability inequality \eqref{def-stability} by choosing a suitable test function $\varphi$. In previous results, including works by Crandall-Rabinowitz \cite{CrandallRabinowitz1975} and Nedev \cite{nedev}, the common choice was $\varphi=h(u)$, where $h$ is a carefully selected function that depends on the nonlinearity $f$. 

A different approach was taken in \cite{Cabre2010} to obtain $L^\infty$ bounds for dimensions $N\leq4$, where the test function was chosen as $\varphi=|\nabla u|\xi(u)$, with a specific selection of $\xi$ depending directly on the solution $u$.

In the case of \cite{cabre-figalli-rosoton-serra}, a central idea in the proof is to employ a test function of the form
\begin{equation}\label{eq:test-function-CFRS}
    \varphi=(x\cdot \nabla u)|x|^{(2-N)/2}\zeta,
\end{equation}
where $\zeta$ is a smooth cut-off function satisfying $0\leq \zeta\leq 1$, equal to $1$ on $B_\rho$ and supported in $B_{3\rho/2}$. This construction allows them to derive the following key estimate: there exists a constant $C$ depending only on the dimension such that

\begin{align}
\label{eq:2n10}
(N-2)(10-N)\int_{B_\rho}|x|^{-N}|x\cdot \nabla u|^2\,dx \leq C\rho^{2-N}\int_{B_{3\rho/2}\setminus B_\rho}|\nabla u|^2\,dx &\\\nonumber
\forall\  0<\rho<{\textstyle\frac23}&
\end{align}
From this inequality, one immediately observes that significant nontrivial information can be obtained when $3 \leq N\leq 9$. Since we may always reduce to $N\geq 3$ by augmenting the number of variables if necessary, the critical dimension assumption here is $N\leq 9$. Therefore, when $N\leq 9$, inequality \eqref{eq:2n10} implies that the radial component of the gradient of $u$ inside a ball is controlled by the total gradient in a surrounding annulus. 
\begin{remark}
    Nevertheless, it is important to highlight that \eqref{eq:2n10} alone does not provide an $L^\infty$ bound for general solutions $u$ of $-\Delta u = f(u)$ when $N \leq 9$. For a counterexample one can consider functions $u$ in $\mathbb{R}^3$ depending on only two variables.
\end{remark}
If one could establish that for stable solutions the radial derivative $x\cdot \nabla u$ and the full gradient $\nabla u$ are comparable in the $L^2$-norm at every scale, then one could control the right-hand side of \eqref{eq:2n10} by
$\int_{B_{3\rho/2}\setminus B_\rho}|x|^{-N}|x\cdot \nabla u|^2\,dx$. This would lead to the estimate
$$
\int_{B_\rho}|x|^{-N}|x\cdot \nabla u|^2\,dx \leq C\int_{B_{3\rho/2}\setminus B_\rho}|x|^{-N}|x\cdot \nabla u|^2\,dx,
$$
and through an iteration argument combined with an appropriate covering, we would ultimately deduce that $u \in C^{\alpha}(B_1)$.
The main arguments to prove the Theorem \ref{thm:L1-CalphaW12} can be summarized as follows: they demonstrate that the radial and total gradients of $u$ are comparable in $L^2$ (under some sort of doubling condition on the Dirichlet energy of $u$ over balls that is not actually a restriction) and use a compactness argument supported by the following sequence of a priori estimates:
\begin{enumerate}
\item[$i)$] Curvature-type bounds for the level sets of $u$, derived using the test function $\varphi = |\nabla u|\eta$ in the stability condition:
Let $u\in C^2(B_1)$ be a stable solution of $-\Delta u=f(u)$ in $B_1\subset \mathbb{R}^N$,  with $f$ locally Lipschitz.    Then, for all $\eta \in C^{0,1}_c(B_1)$ the following is satisfied
\[
\int_{B_1}  \mathcal{A}^2 \eta ^2dx \le \int_{B_1} |\nabla u|^2 |\nabla \eta|^2dx,
\]
where the quantity $\mathcal A$ controls the second fundamental form of the level sets of $u$ defined as:
\begin{equation}\label{defAAA}
\mathcal{A} :=  
\begin{cases} \left( \sum_{ij} u_{ij}^2  - \sum_{i} \left(\sum_{j} u_{ij} \frac{u_j}{|\nabla u|} \right)^2  \right)^{1/2} \quad  \quad &\mbox{if  } \nabla u\neq 0
\\
0 &\mbox{if  } \nabla u=0.
\end{cases}
\end{equation}
\item[$ii)$] An improved $L^{2+\gamma}$ integrability estimate of the gradient, obtained via ($i$) and a suitable energy estimate on each level set of $u$, namely for a.e. $t\in \mathbb{R}$,
\begin{equation}\label{eq:level-set-energy-estimate}
\int_{\lbrace u=t\rbrace\cap B_{3/4}} |\nabla u|^2 d\mathscr{H}^{n-1}  \leq  C;
\end{equation}
\item[$iii)$] a general compactness result for superharmonic functions; see Lemma \ref{strongconvergence}
\item[$iv)$] and the nonexistence of nontrivial $0$-homogeneous superharmonic functions.
\end{enumerate}
Combining Theorem \ref{thm:L1-CalphaW12} with the moving plane method this already proves smoothness of the solutions for convex domains and dimension $N\leq 9$. In 2024, in \cite{figalli-zhang}, the authors were able to achieve the regularity result for finite Morse index solutions, using the moving plane/convexity argument, for supercritical nonlinearities $f(u)\sim u^{\frac{N+2}{N-2}+\epsilon}$ (see Chapter \ref{ch-counterexample} after Proposition \ref{dupaigne_proposition}). In their paper, the result for stable solutions is crucially exploited inside each small enough ball\footnote{Remember that finite Morse index solutions are locally stable.} of an appropriate covering of the domain together with some compactness convergence of finite Morse index sequences of solutions.

For an arbitrary $C^3$ domain, the regularity near/up to the boundary is addressed through the estimates of the following theorem.
\begin{theorem}[Theorem 1.5, \cite{cabre-figalli-rosoton-serra}]\label{thm:globalCalpha}
Let  $\Omega\subset \mathbb{R}^N$ be a bounded domain of class $C^3$. Assume that $f: \mathbb{R} \to \mathbb{R}$ is nonnegative, nondecreasing, and convex. 
Let $u\in C^0(\overline\Omega) \cap C^2(\Omega)$ be a stable solution of \eqref{main_problem}. Then
\begin{equation}
\label{eq:W12g L1 glob}
\|\nabla u\|_{L^{2+\gamma}(\Omega)}\le C\|u\|_{L^1(\Omega)},
\end{equation}
where $\gamma>0$ is a dimensional constant and $C$  depends only on $\Omega$.
In addition, if $N \leq 9$ then
\begin{equation}
\label{eq:C0a L1 glob}
\|u\| _{C^\alpha(\overline\Omega)} \leq C\|u\|_{L^1(\Omega)},
\end{equation}
where $\alpha>0$ is a dimensional constant and $C$ depends only on $\Omega$.
\end{theorem}
This theorem states regularity up to the boundary of the domain of the classical stable solution. Both together and by classical elliptic regularity theory, provide the smoothness of any stable $W^{1,2}_0(\Omega)$ weak solution to semilinear equations \eqref{main_problem}.
\begin{remark}
   To extend the regularity results to any weak stable solution Proposition \ref{prop-approximation-by-stable-solutions} and Proposition \ref{prop:L1} are employed since they allows for an approximation argument by classical stable solutions. 
\end{remark}
Among other implications, the following corollary on the regularity of the extremal solution of \eqref{eq:lambda} closes the Brezis-Vázquez conjecture.
\begin{corollary}[Corollary 1.8, \cite{cabre-figalli-rosoton-serra}]\label{thm:energy-solns}
Let  $\Omega\subset \mathbb{R}^N$ be a bounded domain of class $C^3$. 
Assume that $f: [0,+\infty) \to (0,+\infty)$ is nondecreasing, convex, and superlinear at $+\infty$, and let $u^\star$ denote the extremal solution of \eqref{eq:lambda}.

Then $u^\star \in W^{1,2+\gamma}_0(\Omega)$ for some dimensional exponent $\gamma>0.$
In addition, if $N \leq 9$ then $u^\star$ is bounded and it is therefore a classical solution.
\end{corollary}

Note that, thanks to the superlinearity of $f$, it follows by Proposition \ref{prop:L1} that the $L^1(\Omega)$-norms of the functions $\{u_\lambda\}_{\lambda<\lambda^\star}$ are uniformly bounded by a constant depending only on $f$ and $\Omega$. Hence, by applying Theorem \ref{thm:globalCalpha} to the functions $\{u_\lambda\}_{\lambda<\lambda^\star}$ and letting  $\lambda \uparrow \lambda^\star$, it is immediately deduced that extremal solutions are always $W^{1,2}$ (actually even $W^{1,2+\gamma}$) in every dimension, and that they are universally bounded (and hence smooth) in dimension $N\leq 9$. 
This result answer the Brezis-V\'azquez conjecture and hence also the questions stated above (after Theorem \ref{thm:lambda star} and at the beginning of Section \ref{sec-regularity-literature}) in the relevant family of nonlinearities and smoothly bounded domains. 
\subsubsection{The case  \texorpdfstring{$\mathbf{N\geq 10}$}{Lg}.}  \label{sec-case-N-greater-10}
In view of the results described in the previous sections, it is natural to ask what can one say about stable solutions in dimension $N \geq 10.$ 
Their strategy of proof provide optimal (or perhaps almost optimal) integrability estimates in Morrey spaces in every dimension. Recall that Morrey norms are defined as
\[\|w\|_{M^{p,\beta}(\Omega)}^p:=\sup_{y\in\overline\Omega,\ r>0}r^{\beta-N}\int_{\Omega\cap B_r(y)}|w|^p\,dx,\]
for $p\geq1$ and $\beta\in(0,N)$.

\begin{theorem}[Theorem 1.9, \cite{cabre-figalli-rosoton-serra}]\label{thm11}
Let $u\in C^2(B_1)$ be a stable 
solution of
$$
-\Delta u=f(u) \quad \text{in }B_1\subset \mathbb{R}^N,
$$
with  $f:\mathbb{R} \to\mathbb{R}$ locally Lipschitz.
Assume that $N\geq 10$ and  define
\begin{equation}
\label{eq:pn}
p_N:=\left\{
\begin{array}{ll}
+\infty &\text{if }N=10,\\
\frac{2(N-2\sqrt{N-1}-2)}{N-2\sqrt{N-1}-4} &\text{if }N\geq 11.
\end{array}
\right.
\end{equation}

Then
\begin{equation}
\label{eq:p n11}
\Vert u\Vert_{M^{p,2+\frac{4}{p-2}}(B_{1/2})}\leq C \Vert u\Vert_{L^1(B_1)}\quad\ \text{ for every }\ \
p<p_N,
\end{equation}
where  $C$ depends only on $N$ and $p$.

In addition, if {$f$ is nonnegative and nondecreasing,}  $\Omega\subset \mathbb{R}^N$ is a bounded domain of class~$C^3$, and $u\in C^0(\overline\Omega) \cap C^2(\Omega)$ is a stable solution of \eqref{main_problem}. Then 
\begin{equation}
\label{eq:p n11 global}
\Vert u\Vert_{M^{p,2+\frac{4}{p-2}}(\Omega)}\leq C \Vert u\Vert_{L^1(\Omega)}\quad\ \text{ for every }\ \
p<p_N,
\end{equation}
for some constant $C$ depending only on $p$ and $\Omega$.
\end{theorem}

It is interesting to observe that the above result is essentially optimal.
To see this we recall that, in dimension $N=10$, the function $u=\log(1/|x|^2)$ is an unbounded $W^{1,2}_0(B_1)$ stable solution in $B_1$
(recall Proposition \ref{prop-approximation-by-stable-solutions}).
Also, as shown in  \cite{brezis-vazquez}, for $N\geq 11$ the function $u(x)=|x|^{-2/(q_N-1)}-1$ is the extremal solution of
\begin{equation}
\left\{ \begin{array}{cl}
 -\Delta u  = \lambda^\star (1+u)^{q_N} & \textrm{in }B_{1}\\
u > 0 & \textrm{in }B_{1}\\
 u  =  0 & \textrm{on }\partial B_{1},
\end{array}
\right.
\label{eq:05}
\end{equation}
with $\lambda^\star=\frac{2}{q_N-1}\big(N-2-\frac{2}{q_N-1}\big)$ and $q_N:=
\frac{N-2\sqrt{N-1}}{N-2\sqrt{N-1}-4}$. 
In particular, it is easy to see that $u \in M^{p,2+\frac{4}{p-2}}(B_{1/2})$ if and only if $p\leq p_N$.
It is an open question whether \eqref{eq:p n11} holds with $p=p_N$ for a general stable solution $u$.

It is worth mentioning that Morrey spaces have an optimal embedding into Lebesgue spaces characterized by the optimal exponent $q_{opt}=q_{opt}(p,\beta)$ for which
\begin{equation}
    ||u||_{L^q(\Omega)}\leq ||\nabla u||_{M^{p,\beta}(\Omega)}\ \ \ \ \forall q\leq q_{opt}
\end{equation}
More precisely, we can ask the following natural question: Take for simplicity the unit ball $B_1\subset\mathbb{R}^N$.\\
{\it Given functions $u$ with support in the unit ball $B_1\subset\mathbb{R}^N$ and with gradient in the Morrey space $M^{p,\lambda}(B_1)$, where $1<p<\beta<N$,
what is the  optimal range of exponents $q$ such that necessarily $u\in L^{q}(B_1)$?}.
There was an answer in \cite{Adams.Xiao.2011-morrey-mal} saying that the optimal exponent was $q_{opt}=\frac{Np}{\beta -p}$. In \cite{Cabre-Charro} the optimal exponent was corrected to be the smaller value $q_{opt}=\frac{\beta p}{\beta -p}$ providing also a counterexample for the previous exponent that depends on $\beta$ being integer or not. 
\begin{itemize}
\item If $\beta$ is an integer, such function $u$ can be taken to be
\begin{equation}\label{counterexample.noninteger.general.statement.thm}
u(x)=\left(|x'|^{-\alpha}-2^\alpha\right)_{+}\,\xi(|x''|) 
\end{equation}
where $x=(x',x'')\in\mathbb{R}^{\beta}\times\mathbb{R}^{n-\beta}$, 
 the parameter $\alpha$ satisfies
\begin{equation}\label{hyp.alpha.counterexamples.intro.statement.thm}
\frac{\beta}{q}\leq\alpha<\frac{\beta-p}{p},
\end{equation}
and $\xi:\mathbb{R}^+\to[0,1]$ is a cutoff function with $\xi\equiv 1$ in $[0,1/2)$ and $\xi\equiv0$ in $\mathbb{R}^+\setminus[0,\sqrt{3}/3).$

\item If $k-1<\beta<k$ for some integer $k\in[2,N]$, the function $u$ can be set to,
\begin{equation}\label{counterexample.noninteger.intro.statement.thm}
u(x) =
\left\{
\begin{array}{rl}
&\left(\textnormal{dist}(x,\mathcal{C}_{n,\beta})^{-\alpha}-4^\alpha\right)_{+}\hspace{62pt}\textnormal{if}\ k=N ,\\
&\left(\textnormal{dist}(x',\mathcal{C}_{k,\beta})^{-\alpha}-4^\alpha\right)_{+}\,\xi(|x''|) \qquad\textnormal{if}\ k<N,
\end{array}
\right.
\end{equation}
where $x=(x',x'')\in\mathbb{R}^{k}\times\mathbb{R}^{N-k}$, $\alpha$ satisfies \eqref{hyp.alpha.counterexamples.intro.statement.thm}, $\xi$ is a cutoff function as above, and $\mathcal{C}_{k,\beta}$ is a set of Hausdorff dimension $k-\beta$ in $\mathbb{R}^k$ given by
\[
\mathcal{C}_{k,\beta}=\{0\}\times C_{\gamma}\subset\mathbb{R}^{k-1}\times[-1/2,1/2],
\]
where $C_{\gamma}\subset[-1/2,1/2]$ is the generalized Cantor set with parameter 
$\gamma=1-2^{1-\frac{1}{k-\beta}}$.
\end{itemize}
\begin{remark}
    The generalized Cantor set $C_{\gamma}$ is obtained from the interval $[-1/2,1/2]$ by removing at iteration $j=1,2,\ldots$ the central interval of length $\gamma\,l_{j-1}$ from each remaining segment of length $l_{j-1}=((1-\gamma)/2)^{j-1}$. The usual Cantor set corresponds to $\gamma=1/3$.
The reason for their choice of $\gamma$ is that
the Hausdorff dimension of $C_{\gamma}$ is $ \frac{-\log{2}}{\log{\frac{1-\gamma}{2}}}=k-\beta\in(0,1)$.
 In particular, letting $\beta$ range from $k-1$ to $k$ yields any fractal dimension between 0 and 1, and \eqref{counterexample.noninteger.intro.statement.thm} somehow interpolates between the integer cases $\beta=k-1$ and $\beta=k$.
\end{remark}

The fact that for non-integer $\beta$'s the counterexample relays on the distance to a fractal set, is an uncomfortable issue that the authors and other researchers thought that would be desirable to fix. Nevertheless, providing a non-fractal type counterexample remains a slippery open problem to this day.

\section{Multiplicity of solutions of equations with vanishing nonlinearity}\label{sec-multiplicity-literature}

This section is devoted to present the current status and the progress made during the last decades on the multiplicity of solutions to non-autonomous problems of the following type

\begin{equation}\label{main-multiplicity-problem}
\left \{ 
\begin{array}{rcll}
-\Delta u & = & {g(u)} - h(x)f(u), & \mbox{in } \Omega, \\
u & = & 0, & \mbox{on } \partial\Omega,\\
\end{array}
\right.\tag{$P3$}
\end{equation}
where $\Omega$ is a bounded smooth domain in $\mathbb{R}^N$,    $h\in L^\infty(\Omega)$ is nonnegative, nonzero  and satisfies
\begin{equation}\label{h}
h>0  \mbox{ a.e. in } \Omega \backslash \Omega_0  \mbox{ with } \Omega_0=\mbox{interior}\,\{x\in \Omega\, /\, h(x)=0\},  \ \ {\rm and} \ \ |\Omega_0|>0,
\end{equation} 
and $g$  and  $f$ are  continuous functions.  

Questions concerning the solvability of equations \eqref{main-multiplicity-problem} in the case where
the coefficient $h$ may vanish were introduced in \cite{KazdanWarner1975} in
the context of the prescribed curvature problem on compact manifolds. In the
particular case where $g(u)=\lambda u$ and $f ( u ) = |u|^{p-2}u$, the bifurcation curve of positive solutions to \eqref{main-multiplicity-problem} has been described in \cite{Ouyang1992} and \cite{delPino1994}, and the existence and nonexistence of changing sign solutions have recently been treated in \cite{delPino-Felmer1995}. In particular, they show that (for homogeneous $f$) \eqref{main-multiplicity-problem} admits at least $N(\lambda):=\#\lbrace j\in\mathbb{N}|\lambda_j<\lambda\rbrace-\#\lbrace j\in\mathbb{N}|\lambda_j<\lambda\rbrace$ pairs of nontrivial solutions. Then, in \cite{at1} it is provided a simplified proof of this result. They also emphasized the relationship between the solvability of certain semilinear equations and the number of bound states of an associated spectral problem that had been previously remarked in \cite{AlamaLi1992} and \cite{HEINZ1995} for nonlinear perturbations of Schrödinger operators in $\mathbb{R}^N$. The eigenvalue problem in question is the following
\begin{equation}\label{eq:eigenvalue-multiplicity-problem}
\left \{ 
\begin{array}{rcll}
-\Delta v - \mu W(x)v& = & \lambda v , & \mbox{in } \Omega, \\
v & = & 0, & \mbox{on } \partial\Omega,\\
\end{array}
\right.\tag{$E_\lambda$}
\end{equation}
with $\mu>0$ an auxiliary parameter and where $W(x)$ will be equal to $h(x)$ or $h(x)f(u)/u$ depending on the needs. Let $\lbrace \lambda_i(\mu)\rbrace_{i\in \mathbb{N}}$ denote the discrete set of eigenvalues of \eqref{eq:eigenvalue-multiplicity-problem} (ordered by min-max).  Note that\footnote{Remember that $\lbrace\lambda_i\rbrace_{i\in\mathbb{N}}$ denote the set of eigenvalues of the Laplacian as in \eqref{laplacian_eigenvalue_equation}.}  $\lambda_1(\mu)\rightarrow \lambda_1$ as $\mu\rightarrow 0^+$  and using \cite[Proposition 2.3]{fht} we also have $\lambda_1(\mu)\rightarrow \lambda_1(\Omega_0)$ as $\mu\rightarrow +\infty$ in the resolvent sense. Apart from other assumptions that we will do in Chapter \ref{ch-5}, we will denote and assume the following:
\begin{equation}
    0\leq \mu_0^\pm:=\lim_{u\rightarrow 0^\pm}{\frac{f(u)}{u}}< \mu_\infty^\pm:=\lim_{\pm\infty}\frac{f(u)}{u}.
\end{equation}
It is useful to keep in mind the model cases $f ( u ) = |u|^{p-2}u$ with $p>1$ and $g(u)=\lambda u$ with $\lambda>0$ for which these assumptions are manifestly satisfied, in order to have a better intuition of the behavior of the problem and the ideas behind the proofs. 

In \cite{at1}, just for $g(u)=\lambda u$, $\lambda>0$, the authors obtained existence of a positive and a negative solution of \eqref{main-multiplicity-problem} when $\lambda_1<\lambda<\lambda_1(\Omega_0)$. In addition, if it is assumed $N(\lambda)\geq 2$, they proved existence of a third nontrivial solution (see \cite[Theorem A]{at1}). In the case $\lambda\geq\lambda_1(\Omega_0)$, existence of a nontrivial solution is obtained if $N(\lambda)\geq 1$ and $\lambda\notin\lbrace\lambda_i\rbrace_{i\in\mathbb{N}}$ (see \cite[Theorem B]{at1}).
For $f ( u ) = |u|^{p-2}u$ and even with an extra term in the equation: 
$$-\Delta u=g(u)+|u|^{p-2}u+k(x)\ \ \ \ \text{in}\ \Omega,$$
where $k\in L^2(\Omega)$, it was proved in \cite[Lemma A.3]{AlamaLi1992} that the problem admits a unique solution provided $\lambda<\lambda_1(\Omega_0)$. Existence of solutions for all $\lambda\in\mathbb{R}$ was obtained in \cite{ARCOYA2020111728}.
More recently in \cite{Arcoya2024} these results were extended to the case of general asymptotically linear nonlinearity $g$ and general superlinear nonlinearity $f$. Their approach is rather more abstract than the previous references. They employ some of the usual infinite dimensional critical point theory or Morse theory, in which critical groups (introduced in Section \ref{subsec-morse-theory}) serve as the fundamental tool that distinguish different critical points and gives the relevant information about multiplicity. However, the problem of finding solutions with fixed sign as in \cite{at1} in the general asymptotically linear setting has not been addressed before up to our knowledge. In Chapter \ref{ch-5} we deal with that problem. As we will see, we find existence of positive and negative solutions and up to four non trivial different solutions. The main tools are the method of sub and supersolutions in integral form, to find the two signed solutions, the mountain pass lemma to find the third one and the critical groups of all of them to find and distinguish the fourth solution (see Theorem \ref{th:chang-multiplicity} or \cite[Theorem 3.5]{chang}).



\begin{savequote}[0.75\linewidth]
	This chapter contains the results of the paper ``A negative result on regularity estimates of finite radial Morse index solutions to elliptic problems." published on the Journal Communications on Pure and Applied Analysis Vol. 24, No. 3, March 2025, pp. 358-364. doi:10.3934/cpaa.2024092 (\cite{paper1}).
	\qauthor{} 
\end{savequote}

\chapter{A counterexample to the extended Brezis-Vazquez conjecture}\label{ch-counterexample}

In this chapter we deal with the following question: does the boundedness of the Morse index lead to boundedness of the solution?. We provide a counterexample for a large family of semilinear elliptic equations showing the existence of sequences of solutions with radial Morse index equal to 1 for which regularity estimates can not be satisfied.\\

\noindent\textbf{Keywords:} Stability, elliptic problems, radial solutions, Morse index, regularity
estimates.
\\
\textbf{MSC2020:} Primary: 35J91; Secondary: 35B35, 35B45.

\section{Introduction}

 Now we focus our attention on the properties of the Morse index of solutions $u:\Omega\to \mathbb{R}$ of the nonlinear elliptic Dirichlet problem \eqref{main_problem} with domain $\Omega=B_1$ the unit ball in $\mathbb{R}^N$ with dimension $N\geq 3$ and the nonlinearity $f\in C^1(\mathbb{R})$ being nonnegative and nondecreasing. Under these assumptions, it is well known by the celebrated Gidas-Ni-Nirenberg symmetry result \cite{gidas-ni-niremberg} that $u$ is radially symmetric and decreasing. 
As we have shown in the literature review, quite recently in \cite{cabre-figalli-rosoton-serra} the authors have shown that for nonlinearities $f\geq 0$ in dimensions $3\leq N\leq 9$, stable solutions to \eqref{main_problem} are as smooth as the regularity of $f$ permits, proving interior estimates:
\begin{align}
||\nabla u||_{L^{2+\gamma}(B_{1/2})} \leq\,& C||u||_{L^1(B_1)}\ \ \ \mbox{for } N\geq 1,\\
||u||_{C^{\alpha}(\overline{B}_{1/2})}  \leq\, & C||u||_{L^1(B_1)}\ \ \ \mbox{for } 1\leq N\leq 9\label{interiorL2gamma},
\end{align}
with $\gamma>0$, $\alpha$ and $C$ dimensional constants and analogous results up to the boundary in the case $f$ nondecreasing and convex.  
Furthermore, in \cite[Th.1.6]{villegas} the author proved sharp pointwise estimates for stable solutions of \eqref{main_problem}:
\begin{theorem}[Theorem 1.6, \cite{villegas}]
Let $N\geq 2$, $f\in C^1(\mathbb{R})$, and $u\in H^1_0(B_1)$ be a stable radial solution of \eqref{main_problem}. Then there exists a constant $M_N$ depending only on $N$ such that:
\begin{itemize}
\item[$i)$] If $N < 10$, then $||u||_{L^\infty(B_1 )}\leq M_N ||u||_{H^1(B_1\backslash\overline{B_{1/2}})}$.\label{result_to_reproduce}
\item[$ii)$] If $N = 10$, then $|u(r)|\leq M_N ||u||_{H^1(B_1\backslash\overline{B_{1/2}})}(|\log (r)| + 1),\ \ \forall r \in (0, 1]$.
\item[$iii)$] If $N > 10$, then $|u(r)|\leq M_N ||u||_{H^1(B_1\backslash\overline{B_{1/2}})}r^{-N/2+\sqrt{N-1}+2}\ \ \ \forall r \in (0, 1]$.
\end{itemize}\label{Salvador_theorem}
\end{theorem}
These results provide a good control on the behavior of solutions under the assumption of stability, which can be a too strong assumption in some situations. Hence, it is reasonable to require results in the same direction for a larger class of solutions. We shall use boundedness of radial Morse index as the criterion to the extent the class of solutions under study.

As was mentioned in the previous chapter, Morse index have been employed to characterize the uniform boundedness of a solution. This idea was proposed for the first time in \cite{bahri-lions}. In that paper, the authors consider a class of subcritical nonlinearities, i.e. nonlinearities with asymptotic behavior $f(t)t^{-1}|t|^{-(p-1)}\to C$ as $t\to \pm\infty$ with $1<p<\frac{N+2}{N-2}$ and prove the equivalence between the boundedness of solutions and their Morse index. In \cite{figalli-zhang} the authors achieve similar results for the supercritical case being able to prove uniform boundedness for convex and bounded domains and a suitable compact subclass of nonlinearities (positive, nondecreasing and convex). They also mention that in the critical case, one can find a counterexample to the boundedness of finite Morse index solutions. More precisely, if one define $U:\mathbb{R}^+\times\Omega\to \mathbb{R}$,
\begin{equation*}
U(\lambda,x)=\left(\frac{\sqrt{\lambda N(N-2)}}{\lambda^2+|x|^2}\right)^{\frac{N-2}{2}}
\end{equation*}
then the family of functions,
\begin{equation*}
u_\lambda(x)=U(\lambda,x)-U(1,x)
\end{equation*}
are Morse index $\mathsf{ind}(u_\lambda)=1$ solutions (see \cite{marchis-ianni-pacella} Sec. 5 for a clear exposition of the details) of the problem \eqref{main_problem} for the critical nonlinearity $f_\lambda(u)=(\lambda+u)^{\frac{N+2}{N-2}}$ that exhibit a singular behavior $\lbrace||u_\lambda||_{L^\infty(B_1)}\rbrace\to\infty$, $\lbrace||u_\lambda||_{L^1(B_1)}\rbrace\to 0$  as $\lambda\to 0^+$. This example shows that one cannot expect to get something similar to Theorem \ref{result_to_reproduce} $i)$ only in terms of the finite Morse index assumption. However, it is still interesting to ask if one can control quotients such as $\frac{||\cdot||_{L^p (B_1)}}{||\cdot||_{L^q(B_1)}}$ for $q<p$ under the same conditions (note that $||u_\lambda||_p/||u_\lambda||_q$ remains bounded for $1<q<p<\infty$). Our main result is again a negative answer using a bounded radial Morse index sequence as a counterexample.

\section{Main Result}

\begin{theorem}\label{main_theorem}
Let $1\leq q < p \leq \infty$ with $p>{N}/({N-2})$. For $3\leq N\leq 9$ there exists a sequence $\{u_n\}_{n\in \mathbb{N}}\subset C^\infty(\overline{B_1})$ of solutions to Dirichlet problems \eqref{main_problem} with nonlinearities $\lbrace f_n\rbrace_{n\in \mathbb{N}}\subset C^\infty(\mathbb{R})$ and radial Morse index $\mathsf{ind_r}(u_n)= 1\ \ \forall n\in\mathbb{N}$, such that:
\begin{equation}\label{main_theorem_equation}
\frac{||u_n||_{L^p (B_1)}}{||u_n||_{L^q(B_1)}}\longrightarrow +\infty\ \ \ \ \ \text{as}\ \ n\rightarrow\infty
\end{equation} 
\end{theorem}

\begin{remark}

It is remarkable the fact that we are not able to find divergent sequences for $q<p\leq {N}/({N-2})$. Despite the suggestive fact that $N/(N-2)$ is the critical Sobolev exponent for the embedding $W^{2,1}\hookrightarrow L^p$ in dimension $N$, the ``raison d'être'' of this kind of bound in $p$ and its optimality remains for us an intriguing open question. We want also to emphasize again that this divergent behavior is contrary to the stability case in which the estimate in the Theorem \ref{Salvador_theorem} $i)$ automatically implies the quotient that appears in \eqref{main_theorem_equation} to be bounded.
\end{remark}

Our main result can be interpreted as the impossibility to find a sort of reciprocal result of Prop. \ref{morse_index_properties}-\textit{(\ref{finite_morse_index_for_regular_solutions})}, i.e., the impossibility to prove certain notion of regularity of a solution only from the boundedness of its Morse index. It also remains as an open problem to find sufficient conditions to ensure such a reciprocal.
%

\section{Proof of the main result}

For the sake of clarity before going into the proof of the main theorem we show the following tree lemmas.

\begin{lemma}\label{lemma1}
Let $\Omega=B_1$, $u(x)=u(|x|)\equiv u(r)$ a solution of \eqref{main_problem}. Assume $\mathsf{ind_r}(u)=0$ in $B_{\delta}$ and in $B_{1}\setminus\overline{B_\delta}$ for some $\delta\in (0,1)$. Then $\mathsf{ind_r}(u)\leq 1$ in $B_1$.
\end{lemma}
\noindent {\bf Proof.}
If on the contrary $\mathsf{ind_r}(u)>1$, we could find a two dimensional subspace $X\subset C^1_{0,rad}(B_1)$ such that $Q_{u,B_1}[\varphi]<0$ for any $\varphi\in X\setminus\lbrace 0 \rbrace$. We choose two linearly independent $\varphi_1,\varphi_2\in X$. Clearly $\varphi_2(\delta)\not= 0$ because otherwise 
$$
Q_{u,B_1}[\varphi_2]=Q_{u,B_\delta}\left[\varphi_2|_{_{B_{\delta}}}\right]+Q_{u,B_1\setminus\overline{B_\delta}}\left[\varphi_2|_{_{B_1\setminus\overline{B_{\delta}}}}\right]\geq 0,
$$
which is immediately contradictory. We can define,
$$\varphi=\frac{\varphi_1(\delta)}{\varphi_2(\delta)}\varphi_2-\varphi_1$$
that vanishes at $r=\delta$. Hence
$$
Q_{u,B_1}[\varphi]=Q_{u,B_\delta}\left[\varphi|_{_{B_{\delta}}}\right]+Q_{u,B_1\setminus\overline{B_\delta}}\left[\varphi|_{_{B_1\setminus\overline{B_{\delta}}}}\right]\geq 0,
$$
which is a contradiction.\qed

\begin{lemma}[Lemma 2.1, \cite{cabre-capella}]\label{lemma1,5}
Let $u(x)=u(|x|)=u(r)$ be a non constant radially symmetric solution of \eqref{main_problem}. Then, $u$ is stable in $B_1\setminus\overline{B_{r_0}}$ if and only if:
\begin{equation*}
\frac{\int_{r_0}^1 r^{N-1}u_r^2\omega'^2\ dr}{\int_{r_0}^1 r^{N-1} u_r^2\frac{\omega^2}{r^2}\ dr}\geq N-1\ \ \ \ \ \ \forall \omega\in C^\infty_{0,rad}((r_0,1))\setminus\lbrace 0 \rbrace
\end{equation*}
where we have denoted $u_r:=\frac{\partial u(r)}{\partial r}$.
\end{lemma}

Actually, our Lemma \ref{lemma1,5} is an easy adaptation and not exactly the same statement of \cite[Lemma 2.1]{cabre-capella} which is slightly more general.

\begin{lemma}\label{lemma2}
Let $0<\Psi\in C^\infty ((0,1])$ such that $\Psi(t)=t, \, t\in (0, \alpha)$, for some $\alpha\in (0,1)$. Define $u:B_1\subset \mathbb{R}^N\rightarrow \mathbb{R}$, $u(x)=u(|x|)\equiv u(r)$ for $r\in [0,1]$ by
\begin{equation}
u(r)=\int_r^1\Psi(s^N)s^{1-N}ds\ \  \ \ \forall r\in[0,1].\label{definition_of_u}
\end{equation}
Then, $u\in C^\infty (\overline{B_1})$ is a solution to \eqref{main_problem} for some $f\in C^\infty(\mathbb{R})$ satisfying:

\begin{itemize}
\item[$I)$] $f\geq 0\ \text{in}\  u(B_1)\iff \Psi\label{psi_i}$ nondecreasing.\label{I)}
\item[$II)$] $f'\geq 0\ \text{in}\  u(B_1)\iff \Psi\label{psi_ii}$ is concave.
\end{itemize}
\end{lemma}
\noindent {\bf Proof.}
Since the radial function $u\in  C^\infty (\overline{B_1}\setminus\{ 0\})$ satisfies $u_r<0$ for every $0<r\leq 1$ we have
\begin{equation}\label{eq:nonlinearity-counterexample}
    f(s):=-(\Delta u)({u}^{-1}(s)), \ s\in[0,u(0)) 
\end{equation}
is well defined and satisfies $f\in C^\infty ([0,u(0)))$. On the other hand, from $\Psi(t)=t, \, t\in (0, \alpha)$, we obtain $u_r(t)=-t$ for $t\in (0,\alpha^{1/N})$. Consequently
$$ u(r)=u(0)-\frac{r^2}{2}, \ r\in [0,\alpha^{1/N}).$$
Therefore $u\in  C^\infty (\overline{B_1})$ and $f(s)=N$ for every $s\in (u(0)-(\alpha^{1/N})^2/2,u(0))$. Finally, extending $f$ to $\mathbb{R}$ by $f(s)=N$, if $s\geq u(0)$ and in a $C^\infty$ way if $s<0$, we conclude that $f\in C^\infty(\mathbb{R})$ and $u\in C^\infty (\overline{B_1})$ is a solution to \eqref{main_problem} as claimed.
On the other hand, under radial symmetry \eqref{main_problem} is written as
\begin{equation}\label{1.1-sphericall}
-u''(r)-\frac{N-1}{r}u'(r)=f(u(r)),\, r\in (0,1].
\end{equation}
Taking the first derivative in \eqref{definition_of_u}, 
\begin{equation}\label{first_derivative_u}
-u'(r)=r^{1-N}\Psi(r^N),\, r\in (0,1].
\end{equation}
Taking the second derivative and using \eqref{1.1-sphericall} we have 
\begin{equation}\label{relation-f-psi}
f(u(r))=N\Psi'(r^N)),\, r\in (0,1],
\end{equation}
so the first part $I)$ is proven. Taking derivatives again in \eqref{relation-f-psi} and using \eqref{first_derivative_u},
$$f'(u)=-N^2r^{2(N-1)}\frac{\Psi''(r^N)}{\Psi(r^N)},\, r\in (0,1]$$
and $II)$ follows.

\

\noindent {\bf Proof of Theorem \ref{main_theorem}.}
Take an arbitrary $r_0\in (0,1)$ and define the function,
\begin{equation*}
\Psi_{r_0}(r)=\left\lbrace\begin{array}{ccc}
r         &      if       &     0< r\leq r_0^N\\
\xi_{r_0}(r)    &      if       &      r_0^N< r\leq 1
\end{array}\right.
\end{equation*}
where $\xi_{r_0}(r)$ is a $C^\infty\left([r_0,1]\right)$ strictly increasing and concave function, such that $\Psi_{r_0}\in C^\infty((0,1])$ and chosen to be bounded by 
\begin{equation}\label{bound_on_Psi}
\xi_{r_0}(r)\leq\kappa_Nr_0^N,\ r\in (r_0^N,1],
\end{equation}
with $\kappa_N=\frac{N}{2\sqrt{N-1}}$. Note that this is always possible, since $N\geq 3$ implies $\kappa_N>1$. 
Now we define a radially symmetric function $u_{r_0}:B_1\rightarrow \mathbb{R}$ as
$$
u_{r_0}(r)=\int_r^1\Psi_{r_0}(s^N)s^{1-N}\,ds\ \  \ \ \forall r\in[0,1].
$$
By the Lemma \ref{lemma2} above $u_{r_0}\in C^\infty (\overline{B_1})$ is a solution to \eqref{main_problem} for some $f_{r_0}\in C^\infty(\mathbb{R})$ satisfying $f_{r_0},f'_{r_0}\geq0$. As a result, by Proposition \ref{morse_index_properties}, $u_{r_0}$ must have finite Morse index and be locally stable. Indeed $u_{r_0}$ has zero radial Morse index in $B_{r_0}$ since by \eqref{psi_ii} and \eqref{relation-f-psi} $f_{r_0}(u_{r_0}(r))=N$ and thus $f_{r_0}'(u_{r_0}(r))=0$ for $r<r_0$. It follows immediately that  $u_{r_0}$ is stable in $B_{r_0}$. Taking into account Lemma \ref{lemma1,5} we can prove stability also in the annulus $B_1\setminus\overline{B_{r_0}}$:
\begin{eqnarray*}
\frac{\int_{r_0}^1 r^{N-1}{(u_{r_0})_r}^2\omega'^2\ dr}{\int_{r_0}^1 r^{N-1}{(u_{r_0})_r}^2\ \frac{\omega^2}{r^2}\ dr}&=&
 \frac{\int_{r_0}^1 r^{N-1}(r^{2-2N}\Psi_{r_0}^2(r^N))\omega'^2\ dr}{\int_{r_0}^1 r^{N-1}(r^{2-2N}\Psi_{r_0}^2(r^N))\frac{\omega^2}{r^2}\ dr}\geq  \frac{\int_{r_0}^1 r^{N-1}(r^{2-2N}r_0^{2N})\omega'^2\ dr}{\int_{r_0}^1 r^{N-1}(r^{2-2N}\kappa_N^2 r_0^{2N})\frac{\omega^2}{r^2}\ dr}\geq\\
& \geq & \frac{r_0^{2N}}{\left(\frac{N}{2\sqrt{N-1}}r_0^N\right)^2}\frac{N^2}{4}
 =  N-1
\end{eqnarray*}
where we have used \eqref{bound_on_Psi} and the generalized Hardy inequality ($0<a<b$),
 $$\int_a^b r^{\alpha+1}\omega'^2\ dr\geq\frac{\alpha^2}{4}\int_a^b r^{\alpha-1}\omega^2\ dr\ \ \ \forall \alpha\in \mathbb{R},\ \forall \omega\in C^1_0([a,b])$$
  had been used with $a=r_0,b=1$ and $\alpha=-N$. Then Lemma \ref{lemma1} implies that the radial Morse index $\mathsf{ind_r}(u_{r_0})\leq 1$. 
  Now, we compute the estimates for its $L^p$ and $L^q$ norms in the case $\frac{N}{N-2}<p<\infty$. Note that, in this case, there is no loss of generality in assuming $N/(N-2)<q<p$ since $L^{q_2} (B_1)$ is continuously embedded in $L^{q_1} (B_1)$ if $q_1<q_2$.
\begin{eqnarray}\nonumber
||u_{r_0}||^p_{L^p(B_1)}&=&C_N\int_0^1r^{N-1}\left|\int_r^1s^{1-N}\Psi_{r_0}(s^N)ds\right|^p
dr\\\nonumber
&\geq &C_N\int_0^{r_0} r^{N-1}\left(\int_r^{r_0}s\,ds\right)^p dr= C_N\int_0^{r_0} r^{N-1}\left( \frac{r_0^2-r^2}{2}\right)^p dr \\\nonumber
&= &C_{N,p}\ r_0^{N+2p}.\\\nonumber\\\nonumber
\end{eqnarray}
\begin{eqnarray}
||u_{r_0}||^q_{L^q(B_1)}
&= & C_N\int_0^{r_0}r^{N-1}\left(		\frac{r_0^2-r^2}{2}+\int_{r_0}^1 s^{1-N}\Psi_{r_0}(s^N)ds\right)^q\,dr\\\nonumber
&+&C_N\int_{r_0}^1r^{N-1}\left(\int_r^1s^{1-N}\Psi_{r_0}(s^N)ds\right)^q dr\\\nonumber
&\leq& C_N\int_0^{r_0}r^{N-1}\left(		\frac{r_0^2-r^2}{2}+\kappa_n r_0^N\int_{r_0}^1 s^{1-N}ds\right)^q\,dr\\\nonumber
&+&C_N\int_{r_0}^1r^{N-1}\left(\kappa_n r_0^N \int_r^1s^{1-N}ds\right)^q dr\\\nonumber
&\leq& C_N\int_0^{r_0}r^{N-1}\left(		\frac{r_0^2-r^2}{2}+\kappa_n r_0^N \frac{r_0^{2-N}}{N-2}\right)^q\,dr\\\nonumber
&+&C_N\int_{r_0}^1r^{N-1}\left(\kappa_n r_0^N \frac{r^{2-N}}{N-2}\right)^q dr\\\nonumber
&\leq&C_{N,q}r_0^{2q}\int_0^{r_0}r^{N-1}dr+C_{N,q}r_0^{Nq}\frac{r_0^{N+(2-N)q}}{-N+(N-2)q}=C'_{N,q}r_0^{N+2q}. \\\nonumber 
\end{eqnarray}
Combining these inequalities we obtain
\begin{equation}
\frac{||u_{r_0}||_{L^q(B_1)}}{||u_{r_0}||_{L^p(B_1)}}\leq C_{N,p,q}r_0^{N\left(\frac{1}{q}-\frac{1}{p}\right)}\label{quotient_ending_proof}
\end{equation}
with $C_{N,p,q}$ a dimensional constant depending on $q$ and $p$. Since $q<p$ the quotient \eqref{quotient_ending_proof} goes to $0$ as $r_0\rightarrow 0$. This proves \eqref{main_theorem_equation} in the range $\frac{N}{N-2}<p<\infty$. If $p=\infty$ we see at once that
\begin{eqnarray*}\nonumber
||u||_{L^\infty (B_1)}=u(0)=\int_0^1 s^{1-N}\Psi(s^N)\,ds\geq \int_0^{r_0} s^{1-N}\Psi(s^N)\,ds=\frac{r_0^2}{2}\\\nonumber
\end{eqnarray*}
Now we get,
\begin{equation}\label{Linfty_Lq_estimate}
\frac{||u_{r_0}||_{L^\infty(B_1)}}{||u_{r_0}||_{L^q(B_1)}}\geq C_{N,q}\,r_0^{-\frac{N}{q}}
\end{equation}
 which is unbounded as $r_0\rightarrow 0$.
 
 Finally, by Proposition \ref{dupaigne_proposition} \eqref{radial_morse_index_0_iff_morse_index_0}, if $\mathsf{ind_r}(u_{r_0})=0$ in $B_1$ then $u$ would be stable and, taking $q=1$,\eqref{interiorL2gamma} would imply the quotient $\Vert u_{r_0} \Vert_{\infty} / \Vert u_{r_0}\Vert_1$ to be uniformly bounded by a dimensional constant, fact that we have just proven impossible for small $r_0>0$. This, together with $\mathsf{ind_r}(u_{r_0})\leq 1$ deduced above, shows that $\mathsf{ind_r}(u_{r_0})=1$ for sufficiently small $r_0>0$ and the theorem follows.
\qed

 \medskip




\begin{savequote}[0.75\linewidth]
	This chapter contains parts of the paper ``A priori estimates of stable solutions of the general Hardy-Hénon equation in the ball" that is available in: Electron. J. Qual. Theory Differ. Equ. 2025, No. 71, 1-11. DOI: https://doi.org/10.14232/ejqtde.2025.1.71
)
	\qauthor{} 
\end{savequote}
\chapter{Regularity of stable solutions of the general Hardy-Hénon equation on the ball}\label{ch-4}

This chapter is devoted to the study of stable radial solutions
$u\in H^1(B_1)$ of $-\Delta u=\vert x\vert^\alpha f(u) \mbox{ in } B_1\setminus\lbrace0\rbrace$, where $f\in C^1(\mathbb{R})$ is a general nonlinearity, $\alpha>-2$  and $B_1$
is the unit ball of $\mathbb{R}^N$, $N>1$. We establish the boundedness of such solutions for dimensions $2\leq N<10+4\alpha$ and sharp pointwise
estimates in the case $N\geq10+4\alpha$. 
In addition, we provide, for this range of dimensions, a large family of
stable radially decreasing unbounded $H^1(B_1)$ solutions.
\\
\\
\textbf{Keywords:} Hardy-Hénon equation, regularity estimates, stable solutions.
\\
\textbf{MSC2020:} 35B45, 35B35, 35J61

\section{Introduction}\label{section.intro}

We deal with the semi-stability of radial solutions $u\in H^1(B_1)$ of
 \begin{equation} \label{non-autonomous-equation}
-\Delta u  =  |x|^\alpha{f(u)} \qquad \mbox{in } B_1\setminus\lbrace0\rbrace\tag{P2},
\end{equation}
where $x=(x_1,\dots,x_N)$ and $B_1$ the unit ball of $\mathbb{R}^N$ with measure $\omega_N/N$ for $N>1$. In the whole paper we assume $\alpha>-2$, and $f\in C^1(\mathbb{R})$. We deal with the regularity of stable radially symmetric energy solutions $u\in H^1(B_1)$ of \eqref{non-autonomous-equation}. Denoting by $r:=|x|$ standard arguments show that these radial solutions $u=u(|x|)\equiv u(r)$ are continuous in $r\in (0,1]$ and can be viewed as classical solutions. We denote as usual $u_r:=\frac{<x,\nabla u>}{|x|}$ and constants by capital letters $C,K,C',K'\dots$ that may vary throughout the discussion and computations, with subscripts that determine the dependence of the constants on other parameters.

 The equation \eqref{non-autonomous-equation} corresponds to the Euler-Lagrange equation of the energy functional
 \begin{equation}
     E(u):=\int_{B_1}\left(|\nabla u|^2-|x|^\alpha F(u)\right)dx,
 \end{equation}
where $F(t)=\int_0^tf(s)ds$.
Since $f\in C^1(\mathbb{R})$ and $u$ is continuous, then $f'(u)$ is continuous. It is meaningful to call a solution $u\in H^1(B_1)$ stable if the second variation of the energy is nonnegative definite, i.e.,
\begin{equation}\label{Hénon-stability-inequality}
    d^2E(\varphi,\varphi)=\int_{B_1}\left(|\nabla\varphi|^2-|x|^\alpha f'(u)\varphi^2\right) dx\geq 0,
\end{equation}
for all $\varphi\in C^1(B_1)$ with compact support in $B_1\setminus\lbrace 0 \rbrace$. 

Equation \eqref{non-autonomous-equation} can be regarded as a generalized version of the so-called Hardy-Hénon equation which corresponds to power-type nonlinearity $f(u)=|u|^{p-1}u$. That equation is in particular called the Lane-Emden equation in the autonomous case $\alpha=0$. Another nomenclature has been used for $f(u)=e^u$, the Hénon-Gelfand problem. 

The study of this problem is preceded by the study of the so-called Hénon-equation whose origins goes back to the celebrated work of S. Chandrasekhar \cite[Chapter IV, pp. 87-88]{chandrasekhar1957stellar} on the stellar structure in polytropic equilibrium in astrophysics. 

The pioneering contribution in the autonomous case $\alpha =0$ is due to \cite{gidas-spruck} where the authors prove that for power-type nonlinearity $f(u)=|u|^{p-1}u$ and posed in the whole space $\mathbb{R}^N$, there is no positive solution for $p\in (1,p_S(N))$, where $p_S(N):=2^*-1$ if $N\geq 3$ and $p_S=+\infty$ if $N\leq 2$ is the classical critical Sobolev exponent.

For $p=p_{S}(N)$, the same equation is known to have (up to translation and rescaling) a unique positive solution, which is radial and explicit (see \cite{caffarelli-gidas-spruck}).

Let now $p_{JL}(N)>p_{S}(N)$ denote the so-called Joseph-Lundgren exponent:
$$
p_{JL}(N)=\left\{
\begin{aligned}
+\infty&\quad\text{if $N\leq 10$},\\
\frac{(N-2)^2-{4N}+8\sqrt{N-1} }{(N-2)(N-10)}&\quad\text{if $N\ge 11$}.
\end{aligned}
\right.
$$
This exponent can be characterized as follows: for $p\ge p_{S}(n)$, the explicit function $u_{s}(x)=C_{p,N}\vert x\vert^{-\frac2{p-1}}$ for an appropriate constant $C_{p,N}$ depending only on $p$ and $N$, is a singular solution of the Lane-Emden equation, which is stable if and only if $p\geq p_{JL}(N)$.
It was proved in \cite{farina-classification} that Lane-Emden equation has no nontrivial finite Morse index solution whenever $1<p<p_{JL}(N)$, $p\neq p_{S}(N)$. Through the application of some blow-up analysis technique, such Liouville-type theorems imply interior regularity for solutions of a large class of semilinear elliptic equations: they are known to be equivalent to universal estimates for solutions of
\begin{equation} \label{general}
-Lu = f(x,u,\nabla u)\quad\text{in $\Omega$,}
\end{equation}
where $L$ is a uniformly elliptic operator with smooth coefficients, the nonlinearity $f$ scales like $\vert u\vert^{p-1}u$ for large values of $u$, and $\Omega$ is an open set of $\mathbb{R}^N$. For precise statements, see the work \cite{P-Q-S} in the subcritical setting, as well as its adaptation to the supercritical case by  \cite{davila-dupaigne-farina}.

For the Hénon-Hardy equation (i.e. $\alpha\not=0$, $f(u)=|u|^{p-1}u$), there are also an extensive literature. In \cite{dancer-du-guo}, the authors show the following:

\begin{theorem}(Theorem 1.2, \cite{dancer-du-guo}). Let $u$ be a stable solution of the Hénon-Hardy equation with $\alpha>-2$, $\Omega=\mathbb{R}^N$, and one of the following hypothesis is true: $ 2\leq N\leq 10+4\alpha$ and $ p>1$; or $ N>10+4\alpha$ and $ 1<p<p_{JL}(N,\alpha)$; where we have defined
    \[
    p_{JL}(N,\alpha)=\left\{
    \begin{array}{lr}
    +\infty\quad & \text{if } 2\leq N \leq 10+4\alpha,\\
    \frac{(N-2)^2-2(\alpha+2)(\alpha+N)+2\sqrt{(\alpha+2)^3(\alpha+2N-2)}}{(N-2)(N-4\alpha-10)}&\quad\text{if $N> 10+4\alpha$}.
    \end{array}
    \right.
    \]
Then $u\equiv0$. On the contrary, for $p\geq p_{JL}(N,\alpha)$ the equation admits a family of radially symmetric positive stable solutions.
\end{theorem}

 These results were later generalized in \cite{WangYe2012} for sign-changing solutions where the authors also proved that assuming $0\in\Omega$, no weak solution exists for $\alpha\leq-2$. As they remark in their proof, the result remains valid for any positive, convex and nondecreasing nonlinearity $f$. This indicates that the case $\alpha\leq -2$ carries additional pathological behavior, highlighting the necessity of avoiding that values for $\alpha$. For $\alpha>-2$, $\Omega=\mathbb{R}^N$ they also prove the following two results: 
\begin{enumerate}
    \item[$i)$] If $2\leq N<10+4\alpha$, there is no weak stable solution of the Gelfand-Hénon equation.
    \item[$ii)$] If $2< N<10+4\alpha^-$ where $\alpha^-=\min\lbrace\alpha,0\rbrace$, then any weak solution of the Gelfand-Hénon equation has infinite Morse index.
\end{enumerate}
\begin{remark}
    The dimension range $2\leq N<10+4\alpha$ is optimal, since for $N\geq 10+4\alpha$ $$u(x)=-(2+\alpha)\log|x|+\log[(2+\alpha)(N-2)],$$ is a radial stable weak solution to the Gelfand-Hénon equation. It is remarkable that essentially it coincides with what the inequality \eqref{eq:2n10+4alpha} of the observation made in Section \ref{non-autonomus-future-research} suggests.
\end{remark}

For analogous recent results in the fractional/non local setting we refer to \cite{dupaigne-non-local}, \cite{FazlyHuYang2021}, \cite{hasegawa-ikoma-kawakami}, \cite{hasegawa-ikoma-kawakami-part-II}, \cite{HyderYang2021} and even with a Hardy potential term $|x|^{-2s}u$ on the right-hand side \cite{kim2022finite} or concerning systems \cite{dai2024liouville}, \cite{duong2024liouville} and references therein.

In our case, we are mainly concerned with the regularity of stable radial solutions of the general Hénon-Hardy equation, although we think that our results can be further developed in several directions involving extensions of what is known in the autonomous case. For instance, it would be interesting to investigate the existence and regularity of a extremal solution of the equation $-\Delta u= \lambda |x|^\alpha f(u)$ with the appropriate set of hypothesis over $\alpha$ and $f$ as in Section \ref{sec:brezis-vazquez-conjecture}. The content is developed through a series of lemmas and propositions, arriving at the end to the proof of the main theorem stated at the beginning of the next section.



\section{Statement of the main result and steps of the proof}
We directly proceed to state the main result of the chapter: a classification of the behavior of stable radial solutions of \eqref{non-autonomous-equation} depending on the dimension.
\begin{theorem}\label{principal}

Let $N\geq 2$, $f\in C^1(\mathbb{R})$, $\alpha>-2$ and $u\in H^1(B_1)$ be a
stable radial solution of (\ref{non-autonomous-equation}). Then there
exists a constant $C_{\alpha,N}$ depending only on $\alpha$ and $N$ such that:

\begin{enumerate}

\item[i)] If $N<10+4\alpha$, then $\Vert u\Vert_{L^\infty(B_1)}\leq C_{\alpha,N}
\Vert  u\Vert_{H^1 (B_1\setminus\overline{B_{1/2}})}$.

\

\item[ii)] If $N=10+4\alpha$, then $\vert u(r)\vert \leq C_{\alpha,10+4\alpha} \Vert
 u\Vert_{H^1 (B_1\setminus\overline{B_{1/2}})}\left( \vert \log r\vert +1\right) , \ \ \forall
r\in (0,1]$.

\

\item[iii)] If $N>10+4\alpha$, then $\displaystyle{\vert u(r)\vert \leq
C_{\alpha,N} \Vert  u\Vert_{H^1 (B_1\setminus\overline{B_{1/2}})}\, r^{\gamma(N,\alpha )}\, , \ \
\forall r\in (0,1]}$.
\end{enumerate}
where we have defined the exponent $\,\gamma(N,\alpha):=2-N/2+\alpha/2+\sqrt{(\alpha+2)(\alpha+2N-2)}/2$. 
\end{theorem}

\begin{remark}
    Observe that $\gamma(N,\alpha)=0$ if and only if $N=10+4\alpha$ and $\gamma(N,\alpha)<0$ if and only if $N>10+4\alpha$ so that the pointwise estimates $ii)$ and $iii)$ do not give rise to boundedness of $u$.
\end{remark}

We shall divide the proof in several lemmas and propositions.
\begin{lemma}\label{lemma:key-non-autonomous} Let $N\geq 2$, $f\in C^1(\mathbb{R})$, $\alpha>-2$ and $u\in H^1(B_1)$ be a
stable radial solution of \eqref{non-autonomous-equation}. Let $v\in C^{0,1}(0,1]$ such that $v(1)=0$. Then 
\begin{equation}
    \int_{r_0}^1 t^{N-1} u_r^2(t) \left(v'(t)^2 +\alpha \frac{v'(t)v(t)}{t}+\left(1-N-\alpha\frac{N}{2}\right)\frac{v(t)^2}{t^2} \right)dt\geq 0,
\end{equation}

\noindent for every $r_0\in (0,1)$.

\end{lemma}

\

\noindent {\bf Proof.} Assume, as a first step, that $v\in C^\infty (0,1)$ has compact support and $v\equiv 0$ in $(0,r_0]$ (later we will prove it for any $v\in C^{0,1}(0,1]$ with $v(1)=0$). Differentiating (\ref{non-autonomous-equation}) with respect to $r$ we obtain

$$(-\Delta u)_r=-\Delta u_r+\frac{N-1}{r^2}u_r=\alpha r^{\alpha -1}f(u)+r^\alpha f'(u)u_r.$$

Multiplying by the radial function $u_r v^2$ and integrating by parts yields

$$\int_{B_1}\left( \nabla u_r \nabla (u_r v^2)+\frac{N-1}{r^2}u_r^2 v^2 \right) dx=\int_{B_1}\left( \alpha r^{\alpha -1}f(u)+r^\alpha f'(u)u_r \right) u_r v^2 dx.$$

Therefore

\begin{equation}\label{hola}
\int_{B_1}\left( v^2 \vert\nabla u_r\vert^2+2u_r v\nabla u_r\nabla v+\frac{N-1}{r^2}u_r^2 v^2\right) dx= \int_{B_1}\left( \alpha r^{\alpha -1}f(u)u_r v^2+r^\alpha f'(u)u_r^2 v^2\right) dx.
\end{equation}

On the other hand, since $u$ is stable, we can consider the radial function $u_r v$ (which is a $C^1 (B_1)$ function with compact
support in $B_1\setminus\{ 0\}$) obtaining

\begin{equation}\label{laotra} \int_{B_1}  \vert \nabla (u_r v)\vert^2 dx\geq \int_{B_1} r^\alpha f'(u)(u_r v)^2 dx .
\end{equation}

Subtracting \eqref{hola} from \eqref{laotra} we can assert that

\begin{equation}\label{dfdf}
\int_{B_1}\left( u_r^2 v_r^2-\frac{N-1}{r^2}u_r^2 v^2\right)\geq-\alpha\int_{B_1}r^{\alpha -1}f(u)u_r v^2 dx.
\end{equation}

Let us expand this last term:

$$\int_{B_1}r^{\alpha -1}f(u)u_r v^2 dx=\int_{B_1}(-\Delta u)\frac{u_r v^2}{r}dx=\omega_N\int_{r_0}^1 t^{N-1}\left (-u_{rr}-\frac{N-1}{t}u_r \right)\frac{u_r v^2}{t}dt=$$

$$-\frac{\omega_N}{2}\int_{r_0}^1 \left(t^{2N-2}u_r^2\right)' t^{-N} v^2 dt=\frac{\omega_N}{2}\int_{r_0}^1 t^{2N-2}u_r^2 (t^{-N} v^2)' dt=$$

$$\omega_N\int_{r_0}^1 t^{N-1}\left( \frac{-N}{2}\frac{u_r^2 v^2}{t^2}+\frac{u_r^2 v v'}{t}\right)dt.$$

Combining this and \eqref{dfdf} we obtain 

$$\omega_N \int_{r_0}^1 t^{N-1}\left( u_r^2 v'^2-\frac{N-1}{t^2}u_r^2 v^2\right) dt\geq -\alpha\ \omega_N \int_{r_0}^1 t^{N-1}\left( \frac{-N}{2}\frac{u_r^2 v^2}{t^2}+\frac{u_r^2 v v'}{t}\right)dt,$$

\noindent which is the desired conclusion.

In fact, by standard density arguments, the above formula is also true if we consider a $C^{0,1}(0,1]$ function vanishing in $(0,r_0] \cup \{ 1\}$.

Now, for the sake of clarity of the exposition let us denote
\begin{equation}
   I(a,b;v):=\int_{a}^b t^{N-1} u_r^2(t) \left(v'(t)^2 +\alpha \frac{v'(t)v(t)}{t}+\left(1-N-\alpha\frac{N}{2}\right)\frac{v(t)^2}{t^2} \right)dt.
\end{equation}
For an arbitrary $v\in C^{0,1}(0,1]$ vanishing at $t=1$, we define a radial truncated function
\begin{equation}
    \bar{v}_\epsilon(t)=
\begin{cases}
0&\text{if }0< t<\epsilon,\\
\frac{v(r_0)}{r_0-\epsilon}(t-\epsilon)&\text{if }\epsilon\leq t\leq r_0,\\
v(t)&\text{if }r_0<t\leq 1,
\end{cases}
\end{equation}
for any $0<\epsilon<r_0$. Now, applying the first step of the lemma to the function $\bar{v}_\epsilon$ we have that $I(\epsilon,1;\bar{v}_\epsilon)\geq 0$ and thus
\begin{align*}
    I(r_0,1; v)&\geq -I(\epsilon,r_0;\bar{v}_\epsilon)\\
    &=-\left(\frac{v(r_0)}{r_0-\epsilon}\right)^2\int_{\epsilon}^{r_0}t^{N-1}u_r(t)^2\left[1+\alpha\frac{t-\epsilon}{t}+\left(1-N-\alpha\frac{N}{2}\right)\frac{(t-\epsilon)^2}{t^2}\right]dt.
\end{align*}
Observe that since $u\in H^1(B_1)$ the whole right hand side is uniformly bounded for $\epsilon\in(0,r_0)$. Therefore 
\begin{align*}
    I(r_0,1; v)\geq -\lim_{\epsilon\to 0}I(\epsilon,r_0;\bar{v}_\epsilon)
    =-\left(\frac{v(r_0)}{r_0}\right)^2(2+\alpha)\left(1-\frac{N}{2}\right)\int_{0}^{r_0}t^{N-1}u_r(t)^2dt
    \geq 0,
\end{align*}
as we claimed. \qed

\begin{remark}
    In the last step, it is shown that the assumption $u\in H^1(B_1)$ is essential. Indeed, it is well known (see for instance \cite{brezis-vazquez}), that there exists stable weak solutions to $-\Delta u =C_{N,q}(1+u)^{\frac{q-2}{q}}$ in $B_1$ of the form $u(x)=|x|^q-1$ such that $u\notin H^1(B_1)$ in the range $q \in \left( -\frac{N}{2} + 2 - \sqrt{N - 1},\; -\frac{N}{2} + 1 \right]$, $N\geq 3$.
\end{remark}

\begin{proposition}\label{prop:non-zero-u_r}
    Let $N\geq 2$, $f\in C^1(\mathbb{R})$, $\alpha>-2$ and $u\in H^1(B_1)$ be a non constant stable radial solution of \eqref{non-autonomous-equation}. Then $u_r\not=0$ in $(0,1]$.
\end{proposition}

\noindent {\bf Proof.} Assume by contradiction that there exists $r_1\in (0,1]$ 
 such that $u_r(r_1)=0$. From $u\in H^1(B_1)$, $f\in C^1(\mathbb{R})$ and the radial symmetry we know that $u(r)\in C^3((0,1])$. Thus, in particular there is a constant $C_{r_1}>0$ such that $|u_r(r)|\leq C_{r_1}|r-r_1|$ for any $r\in[r_1/2,r_1]$.
On the other hand, we take the following test function

\begin{equation*}
    v(t)=
\begin{cases}
\frac{t}{r_1-\epsilon}&\text{if }0< t<r_1-\epsilon,\\
\frac{r_1-t}{\epsilon}&\text{if } r_1-\epsilon\leq t\leq r_1,\\
0&\text{if }r_1<t\leq 1,
\end{cases}
\end{equation*}
with $0<\epsilon<r_1/2$. Applying Lemma \ref{lemma:key-non-autonomous} we have that $-I(r_0,r_1-\epsilon;v)\leq I(r_1-\epsilon,r_1;v)$ for any $0<r_0<r_1-\epsilon$. Therefore
\begin{align*}
    -\frac{1}{(r_1-\epsilon)^2}&\int_{r_0}^{r_1-\epsilon}t^{N-1} u_r^2(t)(2+\alpha)\left(1-\frac{N}{2}\right) dt\\ 
    &\leq\int_{r_1-\epsilon}^{r_1}\frac{t^{N-1} u_r^2(t)}{\epsilon^2}\left[1-\alpha  \frac{(r_1-t)}{t}+\left(1-N-\alpha\frac{N}{2}\right)\frac{\left(r_1-t\right)^2}{t^2}\right]dt\\
    &\leq K \int_{r_1-\epsilon}^{r_1}\frac{ u_r^2(t)}{\epsilon^2}dt\leq K C_{r_1}^{\,2}\int_{r_1-\epsilon}^{r_1}\frac{(t-r_1)^2}{\epsilon^2}dt=\frac{K C_{r_1}^2}{3}\,\epsilon.
\end{align*}
Taking limit as $\epsilon\to 0$ we have that
\begin{equation*}
    -\frac{1}{r_1^2}(2+\alpha)\left(1-\frac{N}{2}\right)\int_{r_0}^{r_1}t^{N-1} u_r^2(t)  dt\leq0,
\end{equation*}
that implies that $u_r\equiv0$ in $[r_0,r_1]$ if $N>2$. By the uniqueness of the corresponding Cauchy problem, $u_r\equiv0$ in $(0,1]$ which gives the contradiction for $N>2$. Finally, if $N=2$ we change the test function by
\begin{equation*}
    v(t)=
\begin{cases}
\left(\frac{t}{r_1-\epsilon}\right)^\beta&\text{if }0< t<r_1-\epsilon,\\
\frac{r_1-t}{\epsilon}&\text{if } r_1-\epsilon\leq t\leq r_1,\\
0&\text{if }r_1<t\leq 1,
\end{cases}
\end{equation*}
with $\beta\in\mathbb{R}$ to be chosen later. With similar arguments we can conclude
\begin{equation*}
    -\frac{1}{r_1^{2\beta}}(\beta-1)\left(\beta+\alpha+1\right)\int_{r_0}^{r_1}t^{N-1+2\beta-2} u_r^2(t)  dt\leq0.
\end{equation*}
Choosing any $\beta\in(-1-\alpha,1)$, we obtain $u_r\equiv0$ in $[r_0,r_1]$ which is again a contradiction.\qed
\begin{lemma}\label{lemma:bound-of-u_r}
    Let $N\geq 2$, $f\in C^1(\mathbb{R})$, $\alpha>-2$ and $u\in H^1(B_1)$ be a non constant stable radial solution of \eqref{non-autonomous-equation}. Then
    \begin{equation*}
        \int_{r/2}^ru_r^2\,dt\leq K_{\alpha,N}||\nabla u||^2_{L^2(B_1\setminus \overline{B_{1/2}})}\,r^{3-N+\alpha+\sqrt{(2+\alpha)(2N-2+\alpha)}},
    \end{equation*}
    for any $r\in (0,1]$.
\end{lemma}
\noindent {\bf Proof.} We first consider the case $0<r\leq 1/2$ and define
\begin{equation*}
    v(t)=
\begin{cases}
r^{s-1}t&\text{if }0< t<r,\\
t^s&\text{if } r\leq t\leq 1/2,\\
2^{1-s}(1-t)&\text{if }1/2<t\leq 1
\end{cases}
\end{equation*}
with $s\in\mathbb{R}$ to be determined. By Lemma \ref{lemma:key-non-autonomous} applied for $r_0=r/2$, we have
\begin{align*}
    0\leq r^{2s-2}&\int_{r/2}^r t^{N-1}u_r^2(t)\,(2+\alpha)\left(1-\frac{N}{2}\right)dt\\
    &+\int_{r}^{1/2}t^{N-3+2s}u_r^2(t)\left(s^2+\alpha s+1-N-\alpha\frac{N}{2}\right)dt\\
    &+\int_{1/2}^1t^{N-1}u_r^2(t)\,2^{2-2s}\left[1-\alpha \frac{1-t}{t}+\left(1-N-\alpha\frac{N}{2}\right)\left(\frac{1-t}{t}\right)^2\right]dt.
\end{align*}
Choosing $s_\alpha=-\alpha/2-\sqrt{(2+\alpha)(2N-2+\alpha)}/2$, the second term is indeed zero. Therefore, we deduce
\begin{align*}
     -r^{2s_\alpha-2}&\int_{r/2}^r t^{N-1}u_r^2(t)(2+\alpha)\left(1-\frac{N}{2}\right)dt\\
     &\leq \int_{1/2}^1t^{N-1}u_r^2(t)\,2^{2-2s_\alpha}\left[1-\alpha \frac{1-t}{t}+\left(1-N-\alpha\frac{N}{2}\right)\left(\frac{1-t}{t}\right)^2\right]dt.
\end{align*}
As a result
\begin{align*}
    - r^{2s_\alpha-2}\int_{r/2}^r & \left(\frac{r}{2}\right)^{N-1}u_r^2(t)(2+\alpha)\left(1-\frac{N}{2}\right)dt\\
    &\leq r^{2s_\alpha-2}\int_{r/2}^r t^{N-1}u_r^2(t)(2+\alpha)\left(1-\frac{N}{2}\right)dt\\
     &\leq \int_{1/2}^1t^{N-1}u_r^2(t)\,2^{2-2s_\alpha}\left[1-\alpha \frac{1-t}{t}+\left(1-N-\alpha\frac{N}{2}\right)\left(\frac{1-t}{t}\right)^2\right]dt\\
     &\leq C\,\int_{1/2}^1t^{N-1}u_r^2(t)dt,
\end{align*}
Finally, we obtain
\begin{align*}
    -\,  \frac{(2+\alpha)\left(1-\frac{N}{2}\right)}{2^{N-1}}r^{2s_\alpha-3+N}\int_{r/2}^r & u_r^2(t)dt\leq \frac{C}{\omega_N}||\nabla u||^2_{L^2(B_1\setminus\overline{B_{1/2}})}
\end{align*}
which is the desired conclusion if $N>2$.  If $N=2$, we change the test function by  
\begin{equation*}
    v(t)=
\begin{cases}
r^{-1-\alpha-\beta}t^\beta&\text{if }0< t<r,\\
t^{-1-\alpha}&\text{if } r\leq t\leq 1/2,\\
2^{2+\alpha}(1-t)&\text{if }1/2<t\leq 1.
\end{cases}
\end{equation*}
with $\beta\in (-1-\alpha,1)$ as in the previous proposition. We have
\begin{align*}
    -{r^{-2-2\alpha-2\beta}}\int_{r/2}^r & t^{-1+2\beta}\,u_r^2(t)(\beta-1)(\beta+\alpha+1)\,dt\\
    &\leq\int_{1/2}^1t\,u_r^2(t)\,2^{4+2\alpha}\left[1-\alpha \frac{1-t}{t}-\left(1+\alpha\right)\left(\frac{1-t}{t}\right)^2\right]dt,
    \end{align*}
that gives the desired conclusion once the same argument as above is employed for the last integral. 

On the other hand, let us call $R_{\alpha,N}$ a constant such that $0<R_{\alpha,N}<r^{3-N-2s_\alpha}$  for all $r\in[1/2,1]$. Hence, applying the above conclusion for $r=1/2$
\begin{align*}
    \int_{r/2}^r u_r^2(t)dt & \leq \int_{1/4}^{1/2} u_r^2(t)dt+\int_{1/2}^1u_r^2(t)dt\\
    &\leq K_{\alpha,N} \left(\frac{1}{2}\right)^{3-N-2s_\alpha}||\nabla u||^2_{L^2(B_1\setminus \overline{B_{1/2}})}+\frac{1}{\omega_N}||\nabla u||^2_{L^2(B_1\setminus \overline{B_{1/2}})}\\
    &= \left( K_{\alpha,N} \left(\frac{1}{2}\right)^{3-N-2s_\alpha}+\frac{1}{\omega_N}\right)\,\,||\nabla u||^2_{L^2(B_1\setminus \overline{B_{1/2}})}\\
    &\leq \frac{1}{R_{\alpha, N}}\left[K_{\alpha,N}\left(\frac{1}{2}\right)^{3-N-2s_\alpha}+\frac{1}{\omega_N}\right]||\nabla u||^2_{L^2(B_1\setminus \overline{B_{1/2}})}r^{3-N-2s_\alpha},
\end{align*}
giving rise to the same conclusion again.\qed

\begin{proposition}
    Let $N\geq 2$, $f\in C^1(\mathbb{R})$, $\alpha>-2$ and $u\in H^1(B_1)$ be a non constant stable radial solution of \eqref{non-autonomous-equation}. Then
    \begin{align*}
        |u(r)-u(r/2)|\leq K'_{\alpha,N}||\nabla u||_{L^2(B_1\setminus \overline{B_{1/2}})}\,r^{2-N/2+\alpha/2+\sqrt{(2+\alpha)(2N-2+\alpha)}/2},
    \end{align*}
    for all $r\in(0,1]$.
\end{proposition}

\noindent {\bf Proof.} 
Fix $r\in(0,1/2)$, using Cauchy-Schwartz inequality and Lemma \ref{lemma:bound-of-u_r}, we have the following chain of inequalities:
\begin{align*}
    |u(r)-u(r/2)|= &\left|\int_{r/2}^ru_r(t)\,dt\right|\leq \int_{r/2}^r\left|u_r(t)\right|\,dt \leq\left(\int_{r/2}^ru_r^2(t)dt\right)^{1/2}\left(\int_{r/2}^rdt\right)^{1/2}\\
     &\leq K_{\alpha,N}||\nabla u||_{L^2(B_1\setminus \overline{B_{1/2}})}r^\frac{{3-N+\alpha+\sqrt{(2+\alpha)(2N-2+\alpha)}}}{2}\left(\frac{r}{2}\right)^{\frac{1}{2}}\\
     &\leq K'_{\alpha,N}||\nabla u||^2_{L^2(B_1\setminus \overline{B_{1/2}})}r^\frac{{4-N+\alpha+\sqrt{(2+\alpha)(2N-2+\alpha)}}}{2}.\ \ \ \ \ \ \ \ \ \ \ \ \ \ \ \ \ \ \ \ \ \ \ \ \ \   \qed
\end{align*}
\textbf{Proof of Theorem \ref{principal}.} 
Following the reasoning of \cite{villegas}, let $r\in(0,1]$, there exists $r_1\in(1/2,1]$ such that $r=r_1/2^{m-1}$ for some $m\in\mathbb{N}$. By Sobolev embedding in dimension one we observe that $u(r_1)\leq||u||_{L^\infty(B_1\setminus B_{1/2})}\leq \mathfrak{s}_N||u||_{H^1{(B_1\setminus \overline{B_{1/2}})}}$. Therefore
\begin{align}\label{eq:u(r)-bound}\nonumber
    |u(r)|\leq & |u(r_1)-u(r)|+|u(r_1)|\leq \sum_{k=1}^{m-1}\left|u\left(\frac{r_1}{2^{k-1}}\right)-u\left(\frac{r_1}{2^{k}}\right)\right|+|u(r_1)|\\\nonumber
    \leq & K'_{\alpha,N}||\nabla u||_{L^2(B_1\setminus \overline{B_{1/2}})}\sum_{k=1}^{m-1}\left(\frac{r_1}{2^{k-1}}\right)^{{2-\frac{N}{2}+\frac{\alpha+\sqrt{(2+\alpha)(2N-2+\alpha)}}{2}}}+\mathfrak{s}_{N}||u||_{H^1{(B_1\setminus \overline{B_{1/2}})}}\\
    \leq & \left(K''_{\alpha,N}\sum_{k=1}^{m-1}\left(\frac{r_1}{2^{k-1}}\right)^{{2-\frac{N}{2}+\frac{\alpha+\sqrt{(2+\alpha)(2N-2+\alpha)}}{2}}}+\mathfrak{s}_{N}\right)||u||_{H^1{(B_1\setminus \overline{B_{1/2}})}}.
\end{align}

$\bullet$ If $2\leq N<10+4\alpha$, we have $2-N/2+\alpha/2+\sqrt{(2+\alpha)(2N-2+\alpha)}/2>0$. Then
\begin{align*}
    \sum_{i=1}^{m-1}\left(\frac{r_1}{2^{i-1}}\right)^{2-\frac{N}{2}+\alpha/2+\sqrt{(2+\alpha)(2N-2+\alpha)}/2}\leq
\sum_{i=1}^{\infty}\left(\frac{1}{2^{i-1}}\right)^{2-\frac{N}{2}+\alpha/2+\sqrt{(2+\alpha)(2N-2+\alpha)}/2},
\end{align*}

\noindent which is a convergent series. Hence, equation \eqref{eq:u(r)-bound} implies statement i) of the theorem.

\

$\bullet$ If $N=10+4\alpha$, we have $2-N/2+\alpha/2+\sqrt{(2+\alpha)(2N-2+\alpha)}/2=0$. From \eqref{eq:u(r)-bound} we obtain
\begin{align*}
 \vert u(r)\vert &\leq \left( K'_{\alpha,N}
(m-1)+\mathfrak{s}_N \right) \Vert u\Vert_{H^1(B_1\setminus
\overline{B_{1/2}})}\\
&=\left( K'_{\alpha,N} \left( \frac{\log r_1-\log
r}{\log 2}\right) + \mathfrak{s}_N \right) \Vert
u\Vert_{H^1(B_1\setminus \overline{B_{1/2}})}\\ 
&\leq \left( \frac{K'_{\alpha,N}}{\log 2}+ \mathfrak{s}_N \right)
\left( \vert \log r\vert +1\right) \Vert u\Vert_{H^1(B_1\setminus
\overline{B_{1/2}})},
\end{align*}
\noindent which proves statement ii).

\

$\bullet$ Finally, if $N>10+4\alpha$, we have $2-N/2+\alpha/2+\sqrt{(2+\alpha)(2N-2+\alpha)}/2<0$. Then
\begin{align*}
    \sum_{i=1}^{m-1}\left(\frac{r_1}{2^{i-1}}\right)^{2-\frac{N}{2}+\frac{\alpha}{2}+\frac{\sqrt{(2+\alpha)(2N-2+\alpha)}}{2}}
    =\frac{r^{2-\frac{N}{2}+\frac{\alpha}{2}+\frac{\sqrt{(2+\alpha)(2N-2+\alpha)}}{2}}-r_1^{2-\frac{N}{2}+\frac{\alpha}{2}+\frac{\sqrt{(2+\alpha)(2N-2+\alpha)}}{2}}}{(1/2)^{2-\frac{N}{2}+\frac{\alpha}{2}+\frac{\sqrt{(2+\alpha)(2N-2+\alpha)}}{2}}-1}.
\end{align*}
From this and again \eqref{eq:u(r)-bound}, we deduce
\begin{align*}
    \vert u(r)\vert \leq \left(
\frac{K'_{\alpha,N}}{(1/2)^{2-\frac{N}{2}+\frac{\alpha}{2}+\frac{\sqrt{(2+\alpha)(2N-2+\alpha)}}{2}}-1}+\mathfrak{s}_N \right)
r^{2-\frac{N}{2}+\frac{\alpha}{2}+\frac{\sqrt{(2+\alpha)(2N-2+\alpha)}}{2}}\Vert u\Vert_{H^1(B_1\setminus
\overline{B_{1/2}})}\, ,
\end{align*}
which completes the proof. \qed
\section{Optimality of the Theorem \ref{principal}}
    These a priori estimates are indeed optimal as one can see with the following examples:
\begin{itemize}
    
     \item If $N=10+4\alpha$, then the function $u(x)=|\log |x||$ is a solution to the Hénon-Gelfand equation
    \begin{equation*}
        -\Delta u = (N-2)|x|^{\alpha}e^{u(2+\alpha)}\qquad\text{in }B_1\setminus\lbrace0\rbrace,
    \end{equation*}
    such that $|x|^\alpha f'(u)=\frac{(N-2)^2}{4|x|^2}$. This implies that the stability inequality \eqref{Hénon-stability-inequality} attaining exactly the optimal Hardy constant and hence it is a stable solution.
    \\
    \item If $N>10+4\alpha$, recall that we have that for any $\alpha>-2$, $\gamma(\alpha,N)<0$. We consider the function $u(|x|)=|x|^{\gamma} -1$ for any 
    \begin{align}\label{eq:gamma-range-optimality-10+4alpha}
        \gamma(\alpha,N)\leq \gamma<0, 
    \end{align}
    which is the solution to the power-type equation
    \begin{equation*}
        -\Delta u =|x|^\alpha(-\gamma)(\gamma +N-2)(1+u)^{1+\frac{2+\alpha}{-\gamma}}\qquad\text{in }B_1\setminus\lbrace0\rbrace.
    \end{equation*}
    In this case 
    \begin{align}\label{stability-term-N-bigger-than-10+4alpha}
        |x|^\alpha f'(u)=\frac{(-\gamma+\alpha+2)(\gamma +N-2)}{|x|^2}.
    \end{align}
    Taking into account \eqref{eq:gamma-range-optimality-10+4alpha}, \eqref{stability-term-N-bigger-than-10+4alpha} and comparing again the stability inequality \eqref{Hénon-stability-inequality} with the hardy inequality, the constants satisfies:
    $$(-\gamma+\alpha+2)(\gamma +N-2)\leq \frac{(N-2)^2}{4},$$ 
    implying the stability of $u$.

\end{itemize}

\begin{savequote}[0.75\linewidth]
	This chapter contains the paper ``Existence and multiplicity of solutions of semilinear equations with degenerate nonlinearity." \textit{Preprint}
	\qauthor{} 
\end{savequote}
\chapter{Existence and multiplicity of equations with degenerate nonlinearity}\label{ch-5}
In this Chapter we consider the  semilinear elliptic problem
 \begin{equation*}
\left\{ 
\begin{array}{rcll}
-\Delta u & = & {g(u)} - h(x)f(u) & \mbox{in } \Omega, \\
u & = & 0, & \mbox{on } \partial\Omega,\\
\end{array}
\right.
 \end{equation*}
where $\Omega$ is a smoothly bounded domain of $\mathbb{R}^N$, $h\in L^\infty(\Omega)$ is nonnegative and nontrivial but vanishing in a positive measure set $\Omega_0\subset\Omega$. $g$ is asymptotically linear, $f$ is superlinear and ${g(0)}=f(0)=0$. Including some other asymptotic properties and growth conditions of $g$ and $f$ we obtain the existence of a pair of a positive and a negative solution and up to four solutions with some further assumptions.
\\
\\
\textbf{Keywords:} Multiplicity, critical group, degenerate nonlinearity.
\\
\textbf{MSC2020:} 35J20, 35J61, 58E05

\section{Introduction}

We consider the semilinear elliptic problem
\begin{equation}\label{problem-multiplicity}
\left \{ 
\begin{array}{rcll}
-\Delta u & = & {g(u)} - h(x)f(u), & \mbox{in } \Omega, \\
u & = & 0, & \mbox{on } \partial\Omega,\\
\end{array}
\right.\tag{$P3$}
\end{equation}
where $\Omega\subset\mathbb{R}^N$ is a bounded smooth domain,    $h\in L^\infty(\Omega)$ is nonnegative, nonzero  and satisfies
\begin{equation*}\label{vanishing-property-h}
h>0  \mbox{ a.e. in } \Omega \backslash \Omega_0  \mbox{ with } \Omega_0=\mbox{interior}\,\{x\in \Omega\, /\, h(x)=0\},  \ {\rm and} \  \, |\Omega_0|>0,
\end{equation*} 
and $g$  and  $f$ are  continuous functions.

With respect to the function $g$, we assume that it is asymptotically linear at $t=0$ and at infinity, i.e., there exists $g'(0)$ {and} 
	\begin{equation*}\label{g}
	\lim_{|t|\to\infty}\frac{g(t)}{t}=\lambda.
	\end{equation*}	
We assume also the following growth condition
\begin{equation*}\label{eq-growth-G}
g(t)t-2G(t)\leq C|t|+C,\quad \forall t\in\mathbb R;
	\end{equation*}
where $G(t)=\int_0^tg(s)ds$.
	
In addition, {with respect to the nonlinearity $f$, if we denote $F(t)=\int_0^tf(s)ds$,} we assume that {$f(0)=0$ and there exists $p>1$  for which the following list of hypotheses holds true:} 
	\begin{equation*}\label{f1}
	\exists R>0 \mbox{ such that } (p+1) F(t)\leq f(t)t \ \ \ \mbox{for}  \ \ |t|\geq R, 
	\end{equation*}
	\begin{equation*}\label{f2}
	F(t)>0 , \quad \forall t\in\mathbb R,
	\end{equation*}
	\begin{equation*}\label{f3}
	|f(t)| \leq {C  (1+|t|^{p})}, \quad \forall t\in \mathbb R
	\end{equation*}	
	and
	\begin{equation}\label{f4}
	 (f(t)-f(s)) (t-s)\geq C|t-s|^{p+1},  \forall\ t,s\in\mathbb{R}.
	\end{equation}
In particular, the  {condition \eqref{f4} implies that %
$f$ is non-decreasing and, together to $f(0)=0$, that }
$$
\frac{f(t)}{t}\geq C|t|^{p-1}, \mbox{ for } t\neq 0.
$$

Our interest is to investigate existence and multiplicity of solutions with definite sign. We take the variational approach of the problem. If we define the energy functional
\begin{equation*}\nonumber
E(u)=\int_\Omega \left\lbrace \frac{1}{2}|\nabla u|^2-G(u)+h(x)F(u)\right\rbrace dx,
\end{equation*}
where $G(u)=\int_0^u g(s)ds$ and $F(u):=\int_0^u f(s)ds$. Equation \eqref{problem-multiplicity} is nothing but the Euler-Lagrange equation of the energy above. Weak solutions $u$ are the critical points satisfying:
\begin{equation*}\label{dE}
0=dE_u(\varphi)=\int_\Omega \left\lbrace\nabla u\nabla \varphi -g(u)\varphi+h(x)f(u)\varphi\right\rbrace dx, \ \  \ \ \ \ \ \ \ \forall \varphi\in C^\infty_c(\Omega).
\end{equation*}

We
introduce the operator $H_\infty$ defined as the unique {bounded} self-adjoint operator associated to the quadratic form $a(u) = \int_\Omega |\nabla u |^2\,dx$ with domain
$$
H_D^1(\Omega_0):=\{u\in H_0^1(\Omega)\, : \, u(x)=0\mbox{ a.e. }x\in \Omega \backslash\Omega_0 \}.
$${
More precisely, for each
$v\in H_D^1(\Omega_0)$,  $H_\infty(v)$ is defined as the unique function in  $H_D^1(\Omega_0)$ satisfying 
\[
\int_\Omega \nabla H_\infty (v) \nabla w = \int_\Omega v \, w, \quad \forall w\in H_D^1(\Omega_0).
\]
}

We denote by  $\{\lambda_i(\Omega_0)\}_{i\in \mathbb{N}}$ the spectrum of $H_\infty$, ordered by the min-max principle with eigenvalues repeated according to their multiplicity, and by $\psi_i$ the associated eigenfunctions to $\lambda_i(\Omega_0)$, normalized so that $\int_{\Omega_0} \psi_i\psi_j\, dx=\delta_{i,j}$.  
{Similarly,   $\{\lambda_i\}_{i\in\mathbb{N}}$ denotes the spectrum of $(-\Delta, H_0^1(\Omega))$ with $\varphi_i$ the normalized eigenfunction associated to $\lambda_i$.} 

\subsection{Statement of the main results}

\begin{theorem}\label{main-theorem-multiplicity}
Assume that $\lambda_1<\lambda<\lambda_1(\Omega_0)$,
and that $g'(0)\notin \{\lambda_i\}_{i\in\mathbb{N}}$, then 
\begin{enumerate}
\item \eqref{problem-multiplicity} has two nontrivial solutions if $g'(0)>\lambda_1$, one solution is positive, while the other is negative;
\item \eqref{problem-multiplicity} has  three nontrivial solutions if $p<2^*-1$ and $g'(0)>\lambda_2$;
\item \eqref{problem-multiplicity} has  four nontrivial solutions solution if $p<2^*-1$, $\lambda_k<g'(0)\leq \lambda_{k+1}$, $k\geq 3$, and $G(t)\leq \lambda_{k+1}t^2/2$ for all $t\in\mathbb{R}$.
\end{enumerate}
\end{theorem}

\begin{proposition}\label{prop:multiplicity-non-existence}
Suppose that $g(t)/t\geq \lambda_1(\Omega_0)$, for all $t\neq 0$, in a subset of $\Omega_0$ with positive measure, then \eqref{problem-multiplicity}  has no positive solution. 
\end{proposition}

\subsection{Proof of Theorem \ref{main-theorem-multiplicity}}

\begin{enumerate}
\item  The method of sub-super solutions is employed. Following \cite{at1} we consider the eigenvalue problem \eqref{eq:eigenvalue-multiplicity-problem}. Remember that we denote $\mu_0^+:=\lim_{u\to0^+}\frac{f(u)}{u}$. Take $\mu\geq\mu_0^+$ so that $\lambda_1(\mu)<g'(0)$
We claim that if $t>0$ is small enough, then $\underline{v}:=te_1(\mu_0^+)$ is a sub-solution of \eqref{problem-multiplicity}. Indeed from \eqref{dE},
\begin{multline}
\int_\Omega\left\lbrace \nabla (te_1(\mu_0^+)) \nabla\varphi-g(te_1(\mu_0^+))\varphi+h(x)f(te_1(\mu_0^+))\varphi\right\rbrace dx \\
=-t\int_{\Omega}\left\lbrace -\left(\frac{f(te_1(\mu_0^+))}{te_1(\mu_0^+))}-\mu_0^+\right) h(x)+\frac{g(te_1(\mu_0^+))}{te_1(\mu_0^+)}-g'(0)    \right\rbrace \varphi  e_1(\mu_0^+) dx \\
-t\int_{\Omega}\left\lbrace (\mu-\mu_0^+)h(x)+(g'(0)-\lambda_1(\mu))   \right\rbrace \varphi  e_1(\mu_0^+) dx \leq 0,
\end{multline}
since in the first integral of the right hand side both terms are arbitrarily small as $t\rightarrow 0^+$. Analogously, if we define $\mu_0^-:=\lim_{u\rightarrow 0-}\frac{f(u)}{u}\geq 0$ we get a negative supersolution of the form $\overline{w}:=te_1(\mu_0^-)$ for some $t<0$ sufficiently close to $0$. On the other hand if we fix $\mu>0$ large enough to ensure $\lambda_1(\mu)>\lambda$ and consider the boundary problem,
\begin{equation*}\label{eq-aux2}
\left \{ 
\begin{array}{rcll}
-\Delta \psi - \mu h(x)\psi & = & \lambda \psi , & \mbox{in } \Omega ,\\
\psi & = & 1, & \mbox{on } \partial\Omega,\\
\end{array}
\right.
\end{equation*}
This problem has a unique solution satisfying $\psi(x)\geq c>0$ for all $x\in\overline{\Omega}$. Taking $\overline{v}:=t\psi$ with $t>>0$ sufficiently large to ensure $\overline{v}>\underline{v}$ and $f(t\psi)>\mu \,t\psi$,
\begin{multline}
\int_\Omega\left\lbrace \nabla (t\psi) \nabla\varphi-g(t\psi)\varphi+h(x)f(t\psi)\varphi\right\rbrace dx \\
=t\int_{\Omega}\left\lbrace \left(\frac{f(t\psi)}{t\psi}-\mu\right) h(x)+\left(\lambda-\frac{g(t\psi)}{t\psi}\right)   \right\rbrace \varphi  \psi dx  \geq 0.
\end{multline}
So we get a positive supersolution $\overline{v}$ and working with the large negative values of $t<<0$ another negative subsolution $\underline{w}=te_1(\mu)$ for certain $\mu<0$. The two solutions $u^+$ and $u^-$ arise as: $0<\underline{v}<u^+<\overline{v}$ and $\underline{w}<u^-<\overline{w}<0$. 

\item 
 According to \cite[Lemma 2.1, Ch.III]{chang} and following the same notation, if $u^\pm$ are isolated solutions, they have critical groups (with coefficients group $G$):
\begin{equation*}
C_k(\widetilde{E},u^\pm)=\left \{ 
\begin{array}{rcll}
G & if & k=0, \\
0 & if & k\not= 0,\\
\end{array}
\right.
\end{equation*}
where $\widetilde{E}=E|_{C_0^1(\overline{\Omega})}$. From \cite[Theorem 3.5, Ch.III]{chang} the weak version of the Mountain pass lemma can be employed to obtain another solution $u_3$. This solution is not the trivial one since $g'(0)>\lambda_2$ implies $j=\mathsf{ind}(d^2 E (0))\geq 2$ and then
\begin{equation*}
C_k(E,0)=C_{k-j}(\widetilde{E},0),
\end{equation*}
but this would be a contradiction with the well known fact (recall Theorem \ref{th:chang-mountain-pass}) that any mountain pass critical point satisfies,
\begin{equation*}
C_k(E,u_3)=\left \{ 
\begin{array}{rcll}
G & if & k=1 ,\\
0 & if & k\not= 1.\\
\end{array}
\right.
\end{equation*}
\item  For $p<2^*-1$, the result is a consequence of a contradiction in the Morse inequalities like in \cite[Theorem 3.5, part (iii), Ch.III]{chang}.
\end{enumerate}

\section{Proof of Proposition \ref{prop:multiplicity-non-existence}}

\begin{proof}
By contradiction, let $u$ be a positive solution of \eqref{problem-multiplicity}, the $u$ solves
$$
-\Delta u+w(x)u=p(x)u, \ \ u\in H_0^1,
$$
where $w(x)=h(x)f(u)/u$ and $p(x)=g(u(x))/u(x)$.  It means that $u$ is the first eigenfunction of the eigenvalue problem
\begin{equation*}\label{ep}
\left \{ 
\begin{array}{rcll}
-\Delta v+ w(x)v& = &\lambda p(x)v & \mbox{in } \Omega,\\
v & = & 0 & \mbox{on } \partial\Omega.\\
\end{array}
\right.
\end{equation*}
Moreover, the associated eigenvalue $\lambda_1(-\Delta +w, p,\Omega)=1$.  As a consequence of \cite[Proposition 2.3]{fht}, we have that
$$
1=\lambda_1(-\Delta +w, p,\Omega)\leq \lambda_1(-\Delta, p,\Omega_0).
$$
Thus, denoting by $\psi_p$ the eigenfunction associated to $\lambda_1(-\Delta, p,\Omega_0)$, we get
$$
\lambda_1(-\Delta, g,\Omega_0)\int_{\Omega_0} p(x)\psi_q\psi_1=\int_{\Omega_0} \nabla \psi_g\nabla\psi_1=\lambda_1(\Omega_0)\int_{\Omega_0} \psi_q\psi_1.
$$
It follows that
$$
\int_{\Omega_0} \left(\frac{g(x,u)}{u}-\lambda_1(\Omega_0)\right)\psi_q\psi_1\leq 0,
$$
this is a contradiction provided $g(x,u)/u> \lambda_1(\Omega_0)$ in a subset of $\Omega_0$ with positive measure.

\end{proof}

 

\chapter{Upcoming Research}\label{ch-upcoming-research}

There are many possible different avenues of research from the results obtained in this thesis which are coherent with the wide range of applicability of the selected topic. Some of the most likely interesting ones are in the direction of extending regularity results to the largest possible class of (possibly non autonomous) nonlinearities $f(x,u)$ and (possibly nonlinear) differential operators. 

\section{Possible extensions of regularity results}

\subsection{Completion of the counterexample to the extended Brezis-Vázquez conjecture}

The power of the counterexample in Chapter \ref{ch-counterexample} lies in the fact that the control on the non singular behavior of solutions only by the control on its radial Morse index is an unattainable wish. This is telling us that we could need some control of the full Morse index or additional assumptions on the properties of the solution in order to keep singularities at bay. 

There exists some results bounding the possible relation between the radial Morse index and the full Morse index. For instance in \cite[Theorem 1.1]{radial-vs-full-Morse-index} the authors obtain for $1<p<p_S=2^*-1$ and the problem:
\begin{equation}\label{eq:pacella-problem}
\left\lbrace
\begin{array}{rcll}
-\Delta u &=&  u|u|^{p-1} &\text{ in }B_1\\
u&=&0&\text{ on }\partial B_1
\end{array}\right.
\end{equation}
the following result:
\begin{theorem}[Theorem 1.1, \cite{radial-vs-full-Morse-index}]\label{teoPrincipaleMorse}
Let $N\geq3$ and $u_p$ be a radial  solution to \eqref{eq:pacella-problem} with $m\in\mathbb{N}^+$ nodal regions. Then
\begin{equation}\label{formulina}
\mathsf{ind}(u_p)= m+ N(m-1), \qquad \mbox{ for $p$ sufficiently close to $p_S$.
}
\end{equation}
\end{theorem}
\noindent That, together with:
\begin{theorem}[\cite{HarrabiRebhibSelmi}]\label{teoMorseRadiale} 
Let $\Omega$ be a ball and $f(u)=|u|^{p-1}u$, $p\in (1,p_S)$, $p_S=\frac{N+2}{N-2}$ if $N\geq 3$, $p_S=+\infty$ if $N=2$.
Let $u$ be a radial solution to \eqref{main_problem} with $m\in\mathbb{N}^+$ nodal regions. Then
\[\mathsf{ind_r}(u)=m.\]
\end{theorem}
\noindent give us that, in the above set of hypotheses,
\begin{equation}\label{radialmorseindex-morseindex-relation}
\mathsf{ind}(u_p)= \mathsf{ind_r}(u_p)+ N(\mathsf{ind_r}(u_p)-1), \qquad\mbox{ for $p$ sufficiently close to $p_S$.
}
\end{equation}
This is the kind of result that would make our counterexample in Chapter \ref{ch-counterexample} true for the full or complete Morse index setting. It would lead to stronger implications that should be compared with \cite{bahri-lions} and \cite{figalli-zhang} to see what the implications are on them and vise-versa.
Of course, \eqref{radialmorseindex-morseindex-relation} is not directly applicable to our case, as long as, we do not have full control of the nonlinearities  \eqref{eq:nonlinearity-counterexample} in the sequence of Theorem \ref{main_theorem} and it is not clear if the results of \cite{marchis-ianni-pacella} can be generalized to a larger class of nonlinearities compatible with our counterexample.
However, as De Marchis et al. point out in their paper, there are reasons to think that it could be possible to obtain some bounds on the nonradial negative eigenvalues from the radial ones.

It is observed that the Morse index $\mathsf{ind}(u_p)$ increases proportionally with the number $m$ of nodal domains. This number $m$ also coincides with the count of negative radial eigenvalues of the linearized operator $L_{p}: =-\Delta-p|u_p|^{p-1}$ (\emph{cf.} \cite{HarrabiRebhibSelmi}). This behavior is somewhat unexpected, as it is typically anticipated that the nonradial negative eigenvalues would greatly outnumber the radial ones. In fact, if we examine how the radial and nonradial eigenvalues of the Laplacian $(-\Delta)$ are distributed in the space $H^1_0(B_1)$, we find:

\begin{itemize} 

\item[(i)] On the one hand, from \cite{Bruning2}, we have that
\[\lambda_{r,m}\sim C m^2\quad \mbox{ as }m\rightarrow +\infty\]
where $\lambda_{r,m}$ denotes the $m$-th radial eigenvalue of $(-\Delta)$. This implies that the number $n_r(m^2)$ of radial eigenvalues of $(-\Delta)$ less than or equal to $m^2$ satisfies $n_r(m^2)\sim m \quad \mbox{ as} m\rightarrow +\infty$

\item[(ii)] On the other hand, the classical Weyl law (see e.g. \cite{StraussBook}) provides the estimate: 
\[n(m^2)\sim C m^N\quad \mbox{ as $m\rightarrow +\infty\ $ ($N$ is the dimension)}\]
where $n(m^2)$ is the total number of eigenvalues of $(-\Delta)$ in $H^1_0(B_1)$ not exceeding $m^2$. \end{itemize}

An alternative perspective is to consider a radial eigenfunction of $(-\Delta)$ in $H^1_0(B_1)$ with $m$ nodal zones, i.e., associated with the eigenvalue $\lambda_{r,m}$. The \emph{Morse index} of such a function equals the number of eigenvalues less than $\lambda_{r,m}$, and according to (i) and (ii), this quantity grows like $m^N$, which increases faster than $m$ (for $N\geq 2$) as $m\rightarrow +\infty$.\ Thus, $L_{p}$ provides an example of a Schrödinger-type linear operator with potential $V_p(x)=p|u_p(x)|^{p-1}$, for which  the typical spectral behavior described in (i) and (ii) for $p$ approaching $p_S$, is not reflected in the distribution of its negative eigenvalues.

A further possible implication of this discussion is related to the multiplicity of solutions to \eqref{eq:pacella-problem}, which would offer another illustration of how the Morse index can help to establish the multiplicity results discussed in Chapter \ref{ch-5}. The authors of \cite{radial-vs-full-Morse-index} state a conjecture that could be inferred by comparing the asymptotic behavior of \eqref{eq:pacella-problem} as $p\rightarrow 1$ and as $p\rightarrow p_S$: the Morse index $\mathsf{ind}(u_p)$ converges to the Morse index of the Dirichlet radial eigenfunction of $(-\Delta)$ with $m$ nodal regions (i.e., the one corresponding to $\lambda_{r,m}$), potentially increased by the multiplicity of $\lambda_{r,m}$, which is $1$. Indeed, with appropriate normalizations, solutions of \eqref{eq:pacella-problem} converge to eigenfunctions of the Laplacian as $p\rightarrow 1$ (see \cite{Bonheure}).
Consequently, the above observations suggest that for large $m$, the Morse index $\mathsf{ind}(u_p)$ for $p$ close to $1$ is of order $m^N$, which is significantly greater than $m+ N(m-1)$, which is the value given by \eqref{formulina} for $p$ near $p_S$. This discrepancy indicates the likely occurrence of bifurcations from $u_p$ as $p$ varies between $1$ and $p_S$, thereby implying that the set of solutions to \eqref{eq:pacella-problem} possesses a richer structure than might initially be expected.

%

\

Lastly, an additional remarkable observation is that the formula \eqref{formulina} \emph{does not apply when the spatial dimension is $N=2$ and $p\rightarrow p_S=+\infty$}. In fact, as demonstrated in the paper \cite{DeMarchisIanniPacellaMathAnn}, the following result holds:

\begin{theorem}[\cite{DeMarchisIanniPacellaMathAnn}]\label{teoPrincipaleMorseN2} Let $u_p$ be a radial sign-changing solution to \eqref{eq:pacella-problem} with $2$ nodal regions, with $B\subset\mathbb{R}^2$ and $p_S=+\infty$. Then
\[\mathsf{ind}(u_p)= 12 \qquad\mbox{ for $p$ sufficiently large. }\]
\end{theorem}
In which clearly, $12\neq m+N(m-1)=4$ for $N=2$ and $m=2$. It should be noted that, in this setting, the value of $\mathsf{ind}(u_p)$ turns out to be related with the Morse index of certain radial solutions to the singular Liouville equation in $\mathbb{R}^2$ (see \cite{ChenLin}).

So it would be interesting to research in the following

\textbf{Challenge}: To obtain a more general relation between the radial and nonradial Morse index of solutions to problems of the type \eqref{main_problem} and try to extend our counterexample to the bounded full Morse index scenario and/or to sign-changing solutions.

\subsection{Extensions of the non autonomous case}\label{non-autonomus-future-research}

\subsubsection{General convex domains}
It is worth mentioning the following observation about the non autonomous study of Chapter \ref{ch-4}. Remember that for our non autonomous problem:
\begin{equation*}
\left\lbrace
\begin{array}{rcll}
-\Delta u &=&  |x|^\alpha f(u) &\text{ in }\Omega\\
u&=&0&\text{ on }\partial \Omega.
\end{array}\right.
\end{equation*}
The stability condition read:
\begin{equation*}\label{eq:non-autonomous-second-variation-energy}
    Q_{u,\Omega}(\varphi):=d^2E(\varphi,\varphi)=\int_\Omega \left\lbrace\left|\nabla\varphi\right|^2-|x|^\alpha f'(u)\varphi^2\right\rbrace dx\ \ \ \ \ \forall\varphi\in C^1_0(\Omega)
\end{equation*}
Now looking at a test function of the form $\varphi(x)=(x\cdot \nabla u)|x|^\frac{a}{2}$ as in \eqref{eq:test-function-CFRS} one can realize that when taking the same test function in the stability condition above the two exponents of $|x|$ would add up and the change of test function $\varphi_\alpha=|x|^{-\alpha}\varphi$ merely redefines the exponent $a$ avoiding the effect of the non autonomous factor. In fact, if this is done one gets:
\begin{align}
\label{eq:2n10+4alpha}
(N-2)(10+4\alpha-N)\int_{B_\rho}|x|^{-N}|x\cdot \nabla u|^2\,dx \leq C\rho^{2-N}\int_{B_{3\rho/2}\setminus B_\rho}|\nabla u|^2\,dx &\\\nonumber
\forall\  0<\rho<{\textstyle\frac23}&
\end{align}
for the non autonomous case, attaining the optimal dimension restriction of our Theorem \ref{principal}. If one were able to bound the energy estimate on level sets of $u$ like in \eqref{eq:level-set-energy-estimate}, it is not difficult to see that the moving plane method can be used even in the non autonomous setting \eqref{non-autonomous-equation} to achieve the regularity result in any smoothly bounded convex domain $\Omega\subset\mathbb{R}^N$.
However, the difficulty lies in the fact that when the equation in the stability condition is used deriving with respect to the $x_i$-coordinate, one extra term arises
\[
-\Delta u_{x_i}=\alpha x_i|x|^{\alpha-2}f(u)+|x|^\alpha f'(u)u_{x_i}.
\]
eliminating the possibility of making the stability inequality $f$-independent; the term that depends on the nonlinearity without the derivative can not be canceled with any test function depending only on $u$ since the stability condition depends only on the derivative $f'$ and not on $f$. If it is possible to find a test function $\varphi_f$ depending on the nonlinearity $f$ chosen in such a miraculous way that the new terms cancel exactly the previous extra term and leave the inequality free of $f$-terms seems to be a very difficult task. In general, due to the rapid increasing of complexity of the stability condition with the complexity of test function only a few shapes of test function are useful in practice.

\textbf{Challenge}:
Find out if it is possible to construct a well-shaped test function $\varphi_f$ depending on the nonlinearity $f$ chosen in such a way that the new terms in the stability inequality cancel exactly the previous extra term and leave the inequality free of $f$-terms. If not possible, try to find another test function different from $|\nabla u|\eta$ with $\eta$ a cut-off function (see \cite[Proof of Proposition 1.4]{cabre-figalli-rosoton-serra}) that allows to bound directly 
\[\int_{B_1}\frac{|\nabla u|(-\Delta u)}{|x|}\,\eta^2 dx\leq C\]
which is the only estimate that differs with the autonomous case in the $1/|x|$ factor in the integrand.
\subsubsection{Non zero boundary condition and another non autonomous problem.}

Up to now, the study has been concerned with problems with zero Dirichlet boundary condition. If one imposes a more general non-zero Dirichlet boundary condition on the solution, say
\begin{equation*}\label{eq:non-zero-Dirichlet}
\left\lbrace
\begin{array}{rcll}
-\Delta u &=&  f(u) &\text{ in }\Omega\\
u&=&w&\text{ on }\partial \Omega.
\end{array}\right.
\end{equation*}
then one can redefine $v:=u-w$ and then the new function $v$ satisfies
\begin{equation*}\label{eq:redefined-zero-dirichlet}
\left\lbrace
\begin{array}{rcll}
-\Delta v &=&  f(v+w) -\Delta w&\text{ in }\Omega\\
v&=&0&\text{ on }\partial \Omega.
\end{array}\right.
\end{equation*}
Now, the problem have been converted into a non autonomous one of the form
\begin{equation}\label{eq:redefined-zero-dirichlet-non-autonomous}
\left\lbrace
\begin{array}{rcll}
-\Delta u &=&  f(u)+h(x)&\text{ in }\Omega\\
u&=&0&\text{ on }\partial \Omega.
\end{array}\right.
\end{equation}
One could consider as a first step the case where $h(x)=h(|x|)\equiv h(r)$ is radially symmetric and look also for radial solutions. Again, in the spirit of \cite[Lemma 2.1]{cabre-capella}, the stability inequality for the problem \eqref{eq:redefined-zero-dirichlet-non-autonomous} in the radial case, posed in the unit ball $\Omega= B_1\subset\mathbb{R}^N$, reads
\begin{equation}\label{radial_stability_h}
\int_0^1\left\lbrace \left((N-1)\varphi^2 - (t\varphi)'^2\right)u_r^2  \right\rbrace t^{N-1}dt\leq \int_0^1 h'u_r(t\varphi)^2 t^{N-1}dt
\end{equation}
Unfortunately, this time it is not possible to take advantage of any constant near the origin test function as in \cite{villegas} to bound the $L^2$-norm of radial derivative near the origin by a power or $r$. This is mainly due to the lack of homogeneity of the stability inequality.

\textbf{Challenge:} Try to find another test function that to do the job or a completely different approach to achieve the regularity of radial solutions to this second type of non autonomous equation.

\section{Formalizing the notion of strongly coupled limits in field theories}

A key goal is to systematically characterize when and how certain instabilities arise in field theories. This helps clarify what occurs when the kinematics of some perturbations break down due to vanishing effective terms in their kinetic energy. A degree of freedom becomes \textit{strongly coupled} and hence, does not propagate, when the kinetic matrix is singular, that is, when a perturbation's kinetic term vanishes for certain initial conditions or a specific background solution.

Following \cite{Delhom:2022vae}, we take a theory $\mathcal{S}(\Phi)$ depending on a set of fields $\{\Phi_i\}_{i\in I}$. As we have seen, in order to analyze stability of a solution $\Phi$ of the equations of motion we have to perform perturbations of it and study consequences order by order. Assume now that we introduce a multiplicative bifurcation parameter $\lambda\in\mathbb{R}$ in a neighborhood of $U$ of $1$. The solution fields and its perturbations would be denoted by $\{\Phi_{\lambda}\}$, $\{\varphi_{\lambda}\}$ representing its $\lambda$-dependence and the perturbation is $\Phi_{\lambda}+\epsilon\varphi_\lambda$ with $\varphi_\lambda\in C^\infty (\mathcal{M},V)$. In order to illustrate why vanishing kinetic term should be interpret as strongly coupled degree of freedom, assume that
\begin{itemise}
    \item $\Phi=\lim_{\lambda\to 1}\{\Phi_\lambda\}$, in for example, in $||\cdot||_\infty$,
    \item $\partial_t\varphi_\lambda\to 0$ when $\lambda\to 1$ in $||\cdot||_\infty$, so that the perturbation of one of the fields, denoted by $\varphi$, does not propagate; and
    \item $\det\,\mathcal{K}_{ab}(\lambda)\not=0$ for all $\lambda\in U\setminus \lbrace1\rbrace$.
\end{itemise}
Then, the kinetic matrix of the perturbations would become singular only when $\lambda\to1$. If we expand around any point with $\lambda \not= 1$, the kinetic term for $\varphi$ in the Lagrangian takes the form:
\begin{equation*}
    \mathcal{L}_\lambda=(\partial_t{\varphi}_\lambda)^2 + \mathcal{L}_\text{rest}\,.
\end{equation*}
After canonical normalization,
\begin{equation*}
    \mathcal{L}_\lambda|_\text{canon.}=\frac{1}{2}(\partial_t\varphi_{\lambda})_\text{canon.}^2 + \frac{1}{||\partial_t\varphi_\lambda||_\infty}\mathcal{L}_\text{rest} \,,
\end{equation*}
all the interactions terms of $\varphi$ will become arbitrary large around $\Phi\equiv\Phi_1$
\begin{equation*}
    \lim_{\lambda\to 1}\frac{1}{||\partial_t\varphi||_\infty}||\mathcal{L}_\text{rest}||_\infty= +\infty\,,
\end{equation*}
showing that the process of canonical normalization and taking the limit $\lambda\to 1$ are not interchangeable. Specifically, if a theory typically propagates $n$ modes but only $
m<n$ appear at linear order around a certain background solution, the missing modes may reappear at higher orders. This discrepancy signals a breakdown of the perturbative approach. Such backgrounds often correspond to geometric singularities in phase space. Unless one can compare with other backgrounds or track changes across perturbative orders, these instabilities are hard to detect.
The following example was analyzed in \cite{beltran-jimenez2021} and illustrates the problems with certain solutions in higher-order theories of gravity.

\subsubsection{Strongly coupled cosmological backgrounds in Einsteinian Cubic Gravity} \label{sec: analysis ECG}

We consider the theory described by the action
\begin{equation*} \label{eq: action}
    S = \int d^4x \sqrt{|g|} \left(- \Lambda + \frac{M_{Pl}^2}{2} R + \frac{\beta}{M_{Pl}^2} \mathcal{R}_{(3)} \right)\,,
\end{equation*}
where the first two terms correspond to the General Relativity sector with cosmological constant $\Lambda$, $M_{Pl}$ is the Planck mass, $\beta$ is some parameter and  $\mathcal{R}_{(3)}$, known as the Einsteinian Cubic Gravity (ECG) correction, is given by:
\begin{align*}
    \mathcal{R}_{(3)} & = - \frac{1}{8} \Big(
    12 R_\mu{}^\rho{}_\nu{}^\sigma R_\rho{}^\tau{}_\sigma{}^\eta R_\tau{}^\mu{}_\eta{}^\nu  
    +R_{\mu\nu}{}^{\rho\sigma} R_{\rho\sigma}{}^{\tau\eta} R_{\tau\eta}{}^{\mu\nu}
    +2 R R_{\mu\nu\rho\sigma}R^{\mu\nu\rho\sigma}\nonumber \\
    & \qquad \qquad
    -8 R^{\mu\nu}R_\mu{}^{\rho\sigma \tau}R_{\nu\rho\sigma\tau}
    +4 R^{\mu\nu}R^{\rho\sigma}R_{\mu\rho\nu\sigma}
    -4 R R_{\mu\nu}R^{\mu\nu}
    +8 R_\mu{}^\nu R_\nu{}^\rho R_\rho{}^\mu
    \Big)\,.
\end{align*}
The relative coefficients are precisely chosen to ensure that the linear spectrum of the theory, when expanded around maximally symmetric and cosmological backgrounds, matches that of General Relativity.

In order to study the strong coupling problem of isotropic (Friedmann-Lema\^itre-Robertson-Walker, FLRW) solutions in ECG, in  \cite{beltran-jimenez2021}, the authors analyze the linearized stability of the equations of motion of \eqref{eq: action} around a specific anisotropic deformation of these spacetimes. In particular, they consider the Bianchi I family of metrics described by the line element
\begin{equation*}
d s^2 = - d t^2+a^2(t) d x^2+b^2(t) d y^2+c^2(t) d z^2\,, \label{eq: BianchiI}
\end{equation*}
where $a$, $b$ and $c$ are functions of the cosmic time $t$. To make a more direct contact with the isotropic FLRW solutions, it is convenient to introduce the isotropic expansion rate
\begin{equation*}\label{eq: defH}
H(t):=\frac{1}{3}\left(\frac{\partial_t{a}}{a}+\frac{\partial_t{b}}{b}+\frac{\partial_t{c}}{c}\right)\,.
\end{equation*}
In addition, the anisotropic part is encapsulated in two functions, $\sigma_1(t)$ and $\sigma_2(t)$, defined implicitly by 
\begin{equation*} \label{eq: def sigmas}
\frac{\partial_t{a}}{a}= H + \epsilon_\sigma (2\sigma_1-\sigma_2)\,, \qquad 
\frac{\partial_t{b}}{b}= H - \epsilon_\sigma (\sigma_1-2\sigma_2)\,, \qquad 
\frac{\partial_t{c}}{c}= H - \epsilon_\sigma (\sigma_1+\sigma_2)\,.
\end{equation*}
where $\epsilon_\sigma$ is a certain (not necessarily small) fixed parameter that leads to the isotropic solution in the limit $\epsilon_\sigma\to 0$.

If we represent for an arbitrary function of time $f$ the derivative with respect to the number of e-folds as $f':= H^{-1}\partial_t f$, then the equations of motion can be recast as
\begin{equation*} \label{eq: systemeqs}
\mathcal K _{XY}\begin{pmatrix}H'' \\ X'''\\ Y'''\end{pmatrix} + \mathcal{V}=0\,,
\end{equation*}
where $X(t):= \sigma_1/H(t)$, $Y(t):=\sigma_2/ H(t)$, $\mathcal{K}$ is the matrix of coefficients for the principal part of the equation and $\mathcal{V}$ is a column vector depending on $\{X,Y,X',Y',X'',Y'',H,H'\}$. By imposing $\det \mathcal{K}=0$, where
\begin{align*}
     \det \mathcal{K}_{XY}= - \frac{1458 \beta^3 H^{15}}{M_{Pl}^{12}} &\big(X^2-2Y^2+2Y-X+2XY\big)\big(Y^2-2X^2+2X-Y+2XY\big)\times \nonumber\\
     &\qquad\times\big(Y^2+X^2-4XY-X-Y\big)\,,\label{eq: detHess}
\end{align*}
one finds that this becomes singular for the isotropic solution indicating the presence of a strong coupling issue. For a visualization of this behavior see the figures in the original publication, in which the singular separatrices, i.e., the curves in $(X,Y)$-space in which $\mathcal{K}$ becomes a singular matrix, are also shown.

\textbf{Challenge:} To find a characterization of the occurrence of strong coupling instabilities in any field theory and of the convergence of the perturbative expansion, and in affirmative cases the order in the perturbation from which the instability makes manifest and compare with numerical computations in some new examples.

\singlespacing

\bibliographystyle{apalike}  
\bibliography{bibliography}

\clearpage

\appendix
\pagenumbering{roman}  

\chapter{Standard auxiliary results}

\section{Standard mathematical analysis results}
The reader can check in general the reference \cite[Sec. 8.3]{gilbarg-trudinger} for almost any standard result in elliptic equations.

\begin{theorem}[Sobolev embeddings \cite{evans}]\label{th:sobolev-embeddings}
    Let $N,k\in\mathbb{N}$, $1\leq p<+\infty$, $\Omega\subset\mathbb{R}^N$ bounded domain, such that $\partial\Omega\in C^1$. Then, the following are continuous embeddings of Banach spaces:
    \begin{enumerate}
        \item If $k<N/p$, then $W^{k,p}(\Omega)\hookrightarrow L^q(\Omega)\qquad\forall q\in \left[p,\frac{Np}{N-kp}\right]$,
        \item if $k=N/p$, then $W^{k,p}(\Omega)\hookrightarrow L^q(\Omega)\qquad\forall p\leq q<\infty$ and
        \item if $k>N/p$, then $W^{k,p}(\Omega)\hookrightarrow C^{0,\alpha}(\overline\Omega)$
        where     
        \begin{equation*}\alpha=        \left\{
        \begin{aligned}
        k-N/p && \text{if } k-N/p<1,\\
        \text{any number in } [0,1)       && \text{if } k-N/p=1,\\
        1   && \text{if } k-N/p>1.
        \end{aligned}
        \right.
    \end{equation*}
    \end{enumerate}
    Furthermore, if $\Omega$ is bounded, every previous embedding above is compact except for $q=\frac{Np}{N-kp}$.
\end{theorem}

\begin{lemma}[Lemma A.1, \cite{cabre-figalli-rosoton-serra}]\label{strongconvergence}
Let $v\in L^1(B_R)$ be superharmonic in a ball $B_R\subset\ \mathbb{R}^N$, and let $r\in (0,R)$.
Then:
\begin{itemize}
\item[(a)] The distribution $-\Delta v=|\Delta v|$ is a nonnegative measure in $B_R$, $v\in W^{1,1}_{\rm loc}(B_R)$, 
$$
\int_{B_r}|\Delta v|\leq \frac{C}{(R-r)^2} \int_{B_R}|v|\,dx,
\qquad \mbox{and}
\qquad 
\int_{B_r}|\Delta v|\leq \frac{C}{R-r} \int_{B_R}|\nabla v|\,dx,
$$ 
where $C>0$ is a dimensional constant.
In addition, 
$$\int_{B_r}|\nabla v|\,dx\leq C(N,r,R) \int_{B_R}|v|\,dx
$$ 
for some constant $C(N,r,R)$ depending only on $N$, $r$, and $R$.
\end{itemize}

Assume now that $v_k\in L^1(B_R)$, $k=1,2,...$, is a sequence of superharmonic functions with $\sup_k \|v_k\|_{L^1(B_R)}<\infty$.
Then:
\begin{itemize}
\item[(b1)]  Up to a subsequence, $v_k\to v$ strongly in $W^{1,1}(B_r)$ to some superharmonic function~$v$.

\item[(b2)]
In addition, if for some $\gamma>0$ we have  $\sup_k \|v_k\|_{W^{1,2+\gamma}(B_r)}<\infty$, then $v_k\to v$  strongly in $W^{1,2}(B_r)$.
\end{itemize}
\end{lemma}

\section{Standard theoretical physics results}

\subsection{Noether's Theorem}\label{app:noether}
We introduce the setting of the theorem and some notation:
\begin{itemize}
    \item \( \Phi: \,\mathcal{M} \to V \) be a field defined on spacetime manifold \( \mathcal{M} \). Although the following results can be stated in a local chart of a general spacetime, for the sake of simplicity we keep the exposition restricted to de Minkowski spacetime. Every result remains essentially the same changing $\partial_\mu$ by $D_\mu$. Let us maintain the notation of $\mathcal{M}$ in order not to forget the generality of the theorem.
    \item \( G \) a Lie group of symmetries acting on the fields, with Lie algebra \( \mathfrak{g} \) and let us label its $\dim G$ components by $l\in\lbrace 1,\dots,\dim G\rbrace$,
    \item \( \{ g(\epsilon) \} = \{ \exp({\epsilon X}) \} \subset G \) a one-parameter subgroup with generator \( X \in \mathfrak{g} \) where $\exp$ is the \textit{exponential map} of $G$.
\end{itemize}
The symmetry group acts on fields via representation homomorphism $\rho:G\to GL(V)$:
\[
\Phi(x) \mapsto \Phi_\epsilon(x) = \rho({g(\epsilon)})_A^B\, \Phi^A(x).
\]
Is important to note the component-wise notation is adequate to follow the contractions for finite dimensional $V$, so if $V$ is infinite dimensional we take it as a formal notation. The transformation symmetry changes the action in the following way

\[
S(\Phi_\epsilon) := \int_M \mathcal{L}(\Phi_\epsilon(x), \partial_\mu \Phi_\epsilon(x), x) \, d^n x.
\]
The variation of the action would be computed as follows:
\[
\left. \frac{d}{d\epsilon} \right|_{\epsilon=0}  S(\Phi_\epsilon)= \int_M \left( \frac{\partial \mathcal{L}}{\partial \Phi} \cdot \left. \frac{d \Phi_\epsilon}{d\epsilon} \right|_{\epsilon=0}
+ \frac{\partial \mathcal{L}}{\partial (\partial_\mu \Phi)} \cdot \left. \frac{d}{d\epsilon} (\partial_\mu \Phi_\epsilon) \right|_{\epsilon=0} \right) d^n x.
\]
It is not difficult to see that
\[
\left. \frac{d \Phi_\epsilon}{d\epsilon} \right|_{\epsilon=0} = \rho_*(X)\, \Phi(x),
\quad
\left. \frac{d}{d\epsilon} \partial_\mu \Phi_\epsilon \right|_{\epsilon=0} = \partial_\mu (\rho_*(X)\, \Phi(x)).
\]
Thus:
\[
\left. \frac{d}{d\epsilon} S[\Phi_\epsilon] \right|_{\epsilon=0}
= \int_M \left( \frac{\partial \mathcal{L}}{\partial \Phi^A} \rho_*(X)_B^A\, \Phi^B(x) 
+ \frac{\partial \mathcal{L}}{\partial (\partial_\mu \Phi^A)}\partial_\mu (\rho_*(X)_B^A\, \Phi^B(x)) \right) d^n x.
\]
Now integrate the second term by parts:
\[
\int_M \frac{\partial \mathcal{L}}{\partial (\partial_\mu \Phi^A)}\partial_\mu (\rho_*(X)_B^A\, \Phi^B(x))  \, d^n x
= - \int_M \partial_\mu\frac{\partial \mathcal{L}}{\partial (\partial_\mu \Phi^A)} \rho_*(X)_B^A\, \Phi^B(x)  \, d^n x.
\]
where we have assumed that the boundary terms after the integration by parts vanish (e.g., fields with compact support), the total variation becomes:
\[
 S = \int_M \left\lbrace\left( \frac{\partial \mathcal{L}}{\partial \Phi^A} - \partial_\mu  \frac{\partial \mathcal{L}}{\partial (\partial_\mu \Phi^A)}  \right)  \rho_*(X)_B^A\, \Phi^B(x) +\partial_\mu \left(\frac{\partial \mathcal{L}}{\partial (\partial_\mu \Phi^A)} \rho_*(X)_B^A\, \Phi^B(x)\right)\right\rbrace\, d^n x.
\]

If \( \Phi \) satisfies the Euler–Lagrange equations:
\[
\frac{\partial \mathcal{L}}{\partial \Phi} - \partial_\mu  \frac{\partial \mathcal{L}}{\partial (\partial_\mu \Phi)} = 0,
\]

\begin{theorem}[\cite{noether}]\label{th:noether-theorem}
    If the transformation is a symmetry, i.e., $S (\Phi_\epsilon)= S (\Phi)$ for all $t$, then the following is a conserved current:
\[
j^\mu = \frac{\partial \mathcal{L}}{\partial (\partial_\mu \Phi^A)} (\rho_*(X)_B^A\, \Phi^B(x)), \quad \text{with } \partial_\mu j^\mu = 0.
\]
\end{theorem}
As we already said, when $\mathcal{M}=\mathbb{R}^{1,3}$ and $G=\lbrace T_v:\mathbb{R}^{1,3}\to\mathbb{R}^{1,3},\, T_v(x)^\mu=x^\mu +v^\mu :v\in \mathbb{R}^{1,3}\rbrace$, i.e, the group of (constant) translations, the conserved current $j^\mu$ is just the canonical stress-energy-momentum tensor \eqref{def:canonical-stress-energy-tensor} and the time component lead to the conservation of energy.

Imagine now that we allow for spacetime dependence on the parameter $\epsilon=\epsilon(x^\mu)$ (this is very naturally understood, under the construction of a principal fiber bundle $\mathcal{P}_G$ over $\mathcal{M}$, as taking a non constant section of the bundle) or in other words, given a basis $\lbrace X_1,\dots,X_{\dim G}\rbrace$ of the Lie algebra $T_eG$, we take the variation 
\begin{align*}
    \Phi_{\epsilon^lX_l}=\rho_{\exp(\epsilon^l(x)X_l)}\Phi
\end{align*}
for $\epsilon^l(x)$ any compactly supported smooth functions over $\mathcal{M}$. In this case we have to take care since now there will appear terms involving derivatives $(\partial_\mu \epsilon^l(x))X_l$. After these derivatives over the of the local parameters are integrated by parts we impose the invariance of the action and that the Euler-Lagrange equations are satisfied. We end up obtaining
\begin{align*}
    \int_\mathcal{M}\left\lbrace \partial_\mu \, \tilde j^\mu(\epsilon^l(x))+\epsilon^l(x)\mathcal{D}_l\left( \frac{\partial \mathcal{L}}{\partial \Phi} - \partial_\mu \left( \frac{\partial \mathcal{L}}{\partial (\partial_\mu \Phi)} \right) \right)\right\rbrace \,d^4x
\end{align*}
By divergence theorem we eliminate the term with $\partial_\mu\tilde j^\mu$ and then as $\epsilon^l(x)$ are arbitrary independent test functions we get $\dim G$ identities $\mathcal{D}_l (\text{Euler-Lagrange-equations})=0$ that are true regardless $\Phi$ is a solution of Euler-Lagrange equations and a posteriori the conservation of the current $\partial_\mu\tilde j^\mu=0$. All of this is the second version (and most used in particle physics) of Noether's Theorem:

\begin{theorem}[\cite{noether}]
If the action is invariant under local transformations (gauge transformations):
\[
S\left(\Phi^{\epsilon^l(x)X_l}\right) = S(\Phi)\qquad\left(\Rightarrow\left. \frac{d}{d\epsilon} S(\Phi^{\epsilon^lX_l}) \right|_{\epsilon=0} = 0,
\quad \text{for all } \epsilon^l \in C^\infty_c(M, \mathfrak{g}).\right).
\]
Then, there are $\dim G$ different restriction identities among :
\[\mathcal{D}_l
\left( \frac{\partial \mathcal{L}}{\partial \Phi} - \partial_\mu \left( \frac{\partial \mathcal{L}}{\partial (\partial_\mu \Phi)} \right) \right)  = 0\qquad l=1,\dots,\dim G.
\]
where $\mathcal{D}_l=\mathcal{D}_l(\partial_{\mu})$ is a differential operator. Furthermore, there is a conserved current
\begin{align*}
    \partial_\mu\, \tilde j^\mu=0
\end{align*}
where $\tilde j^\mu=\tilde j ^\mu(\epsilon^l)$ differs from $j^\mu$ only by terms proportional to the Euler-Lagrange equations.
\end{theorem}

This kind of identities are called Noether or Bianchi identities. An example of this identity are the homogeneous Maxwell equations:
\[
\partial_\lambda F_{\mu\nu} + \partial_\mu F_{\nu\lambda} + \partial_\nu F_{\lambda\mu} = 0\qquad(\Leftrightarrow  dF=0 \text{\, in the language of forms}).
\]

\begin{remark}\ 
    \begin{itemize}
    \item If the symmetry depends on arbitrary functions \( \epsilon^l(x) \), then there are dependencies among the Euler–Lagrange equations.
    \item This reflects a redundancy in the description of the degrees of freedom, as occurs in gauge theories.
    \item These identities are essential in the Batalin–Vilkovisky formalism \cite{BATALIN198127}, BRST quantization \cite{figueroa-kimura1991}, and constraint analysis.
\end{itemize}

In summary:
\begin{itemize}
    \item \textbf{Noether's First Theorem}: \emph{Global (rigid) symmetries} of the action functional imply the existence of \emph{conserved currents}.
    \[
    \text{If } \mathcal{S}(\Phi) \text{ is invariant under } \Phi \mapsto \Phi_\epsilon = \rho(\exp(\epsilon X)) \Phi, \text{ then } \partial_\mu j^\mu = 0,
    \]
    for some current \( j^\mu \) constructed from the symmetry generator \( X \).

    \item \textbf{Noether's Second Theorem}: \emph{Local (gauge) symmetries} depending on arbitrary functions (e.g., \( \epsilon^l(x) \)) lead to \emph{differential identities among the Euler--Lagrange equations}.
    
    In other words, if  $\mathcal{S}\left(\Phi^{\epsilon^lX_l}\right) = \mathcal{S}(\Phi)$ for any  $\epsilon^l(x)$, then there exist Noether identities and a identically  conserved current of the form
    
\[\mathcal{D}_l
\left( \frac{\partial \mathcal{L}}{\partial \Phi} - \partial_\mu \left( \frac{\partial \mathcal{L}}{\partial (\partial_\mu \Phi)} \right) \right)  = 0,\qquad \partial_\mu \tilde j^\mu=0.
\]
\end{itemize}

\end{remark}

\subsection{Ostrogradsky’s Theorem}\label{app:ostrogradsky}

Let $\mathcal{L} = \mathcal{L}(\Phi, \partial_t\Phi,\dots,\partial^{(n)}_t\Phi)$, define the $k$-th-order kinetic matrix as 
\begin{align*}
    \mathcal{K}^{(k)}_{ab}:=\frac{\partial \mathcal{L}}{\partial\left(\partial^{(k)}_t\Phi^a\right)\partial\left(\partial^{(k)}_t\Phi^B\right)}.
\end{align*}
We say that a Lagrangian is not degenerate if $\det \mathcal{K}^{(n)}_{ab}(x)\not= 0 $ for all $x\in\mathcal{M}$.
\begin{theorem}[\cite{Ostrogradsky:1850fid}]\label{ostrogradsky}
 Let $n\in\mathbb{N}$, and let $\mathcal{L} = \mathcal{L}(\Phi, \partial_t\Phi,\partial_t^2\Phi\dots,\partial_t^{(n)}\Phi)$ for $n>2$ be a non-degenerate Lagrangian. 
Then the associated Hamiltonian is linear in at least one canonical momentum.
\end{theorem}
\begin{remark}
    The Ostrogradsky's Theorem tell us that the corresponding Hamiltonian (and thus the energy of the system) becomes unbounded from below and above. Consequently, the system exhibits a runaway or ghost type instability.
\end{remark}

\end{document}